%% file: pointsource0.tex
\pdfoutput=1
\documentclass[a4paper,english]{jnsao}
\usepackage[utf8]{inputenc}
\usepackage[T1]{fontenc}
\usepackage{csquotes}
\usepackage{babel}
\usepackage[inline]{enumitem}
\usepackage{mathtools}
\usepackage[indLines=false,noEnd=false]{algpseudocodex}
\usepackage[nameinlink,capitalize]{cleveref}
\usepackage{float}

\input{symbols}

\input{symbols-texonly}

\manuscriptsubmitted{2022-12-07}
\manuscriptaccepted{2023-09-20}
\manuscriptvolume{4}
\manuscriptnumber{10433}
\manuscriptyear{2023}
\manuscriptdoi{10.46298/jnsao-2023-10433}
\manuscriptlicense{CC-BY-SA 4.0}

\theoremstyle{definition}
\newtheorem{assumption}[theorem]{Assumption}
\crefname{assumption}{Assumption}{Assumptions}

\floatstyle{ruled}
\floatname{algorithm}{Algorithm}
\newfloat{algorithm}{tbp!}{loa}
\crefname{algorithm}{Algorithm}{Algorithms}
\numberwithin{algorithm}{section}

\usepackage{tikz}
\usetikzlibrary{external}
\tikzset{external/system call={lualatex \tikzexternalcheckshellescape -halt-on-error -interaction=batchmode -jobname "\image" "\texsource"}}
\tikzexternaldisable
\input{plotting.tex}

\author{
    Tuomo Valkonen\thanks{ModeMat, Escuela Politécnica Nacional, Quito, Ecuador \emph{and} Department of Mathematics and Statistics, University of Helsinki, Finland, \email{tuomo.valkonen@iki.fi}, \orcid{0000-0001-6683-3572}}
    }
\shortauthor{Valkonen}
\title{Proximal methods for point source localisation}

\acknowledgements{%
    This research has been supported by the Escuela Politécnica Nacional internal grant PIIF-20-01 and the Academy of Finland grant 314701.
}

\begin{document}

\maketitle

\begin{abstract}
    Point source localisation is generally modelled as a Lasso-type problem on measures.
    However, optimisation methods in non-Hilbert spaces, such as the space of Radon measures, are much less developed than in Hilbert spaces.
    Most numerical algorithms for point source localisation are based on the Frank–Wolfe conditional gradient method, for which ad hoc convergence theory is developed.
    We develop extensions of proximal-type methods to spaces of measures.
    This includes forward-backward splitting, its inertial version, and primal-dual proximal splitting.
    Their convergence proofs follow standard patterns.
    We demonstrate their numerical efficacy.
\end{abstract}

\section{Introduction}

The point source localisation problem \cite{candes2014towards,lindberg2012mathematical} reads
\begin{equation}
    \label{eq:intro:problem}
    \min_{0 \le \mu \in \Meas(\Omega)}\, \frac{1}{2}\norm{A\mu-b}^2 + \alpha \norm{\mu}_{\Meas(\Omega)}
\end{equation}
for a regularisation parameter $\alpha>0$, measurement data $b \in \R^n$, and a forward operator $A \in \linear(\Meas(\Omega); \R^n)$ in the space $\Meas(\Omega)$ of Radon measures on $\Omega \subset \R^m$.
Commonly $A$ consists of several convolution operators over the sensors of a sensor grid.
Most iterative algorithms in the literature for \eqref{eq:intro:problem} are based on the Frank–Wolfe conditional gradient method \cite{frank1956algorithm}.
Forward-backward splitting and other proximal-type methods are notably absent.
We want to understand \emph{why are true proximal-type methods not amenable to \eqref{eq:intro:problem}---or are they?}

The basic scheme of conditional gradient methods \cite{brediespikkarainen2013inverse,denoyelle2019sliding,duval2017sparse,walter2019linear,blank2017extension,courbot2021fast} for \eqref{eq:intro:problem} is to add a single Dirac measure (or spike) to the discrete measure $\mu$, and then optimise the weights of the spikes so far inserted. Repeat.
The location of the spike is discovered by maximising $\abs{A_*(A\mu-b)}$, where $A_* \in \linear(\R^n; \FullPredual)$ is a pre-adjoint of $A$.
This is a difficult non-convex optimisation problem.
Some also include spike sliding \cite{denoyelle2019sliding}.
The approach has been extended to curve discovery \cite{bredies2022generalized}.
Semismooth Newton approaches have also been proposed \cite{casas2012approximation,casas2013parabolic}, and recently, an approach based on semi-infinite programming \cite{flinth2020linear}.
In some cases particle gradient descent can be used after lifting the measures to a higher-dimensional space \cite{chizat2021sparse,chizat2018global}. The number of particles (or point sources) is, however, fixed, although many-particle limits are also analysed. In \cite{chizat2023convergence} standard Bregman-proximal methods are applied to point source reconstruction by working with densities with respect to a fixed reference measure. In practical application, the reference measure is discretised, so the method is not \emph{grid-free}: the sources cannot be located at arbitrary points of the domain $\Omega$.

Forward-backward splitting is commonly used for $\ell_1$-regularised regression (Lasso), which is a discrete, finite-dimensional variant of \eqref{eq:intro:problem}. With $F(\mu)=\frac{1}{2}\norm{A\mu-b}^2$, we can also sketch forward-backward splitting for \eqref{eq:intro:problem} as solving on each step $k \in \N$ the surrogate problem
\begin{equation}
    \label{eq:intro:fb}
    \mu^{k+1} \in \argmin_{0 \le \mu \in \Meas(\Omega)}\, F(\this\mu) + \dualprod{F'(\this\mu)}{\mu-\this\mu} + \alpha \norm{\mu}_{\Meas}
    + \text{proximal penalty}.
\end{equation}
In Hilbert spaces the proximal penalty would be $\frac{1}{2\tau}\norm{\mu-\this\mu}^2$ for a step length parameter $\tau>0$ satisfying $\tau L < 1$ for a Lipschitz factor $L$ of $F'$.
This choice is not practical in non-Hilbert spaces, as a crucial (Pythagoras') three-point identity does not hold.
Such an identity, however, holds for Bregman divergences \cite{bregman1967relaxation}; see \cite{tuomov-firstorder}. The Frank–Wolfe method can, moreover, be seen a forward-backward method for a modified but equivalent problem\footnote{$\norm{\mu}_{\Masses}$ is replaced by $\phi(\norm{\mu}_{\Masses})$ for $\phi$ that quadratically penalises values greater than $F(0)/\alpha$.} with zero as the proximal penalty. However, as $F'$ cannot be Lipschitz with respect to such a degenerate proximal penalty, a form of line search is employed to ensure descent. In the case of \eqref{eq:intro:problem}, this amounts to the finite-dimensional weight optimisation subproblems.
These considerations raise the question whether we could develop a distance on $\Meas(\Omega)$ that would give a practical forward-backward method?

As the first stage of our endeavour to understand \eqref{eq:intro:problem}, documented in this manuscript, we were able to do so simply by constructing on $\Meas(\Omega)$ the semi-inner product
\[
    \iprod{\mu}{\nu}_{\Wave} \defeq \dualprod{\Wave \mu}{\nu}
    \quad
    (\mu, \nu \in \Meas(\Omega))
\]
for a  “particle-to-wave” operator $\Wave: \linear(\Meas(\Omega); \FullPredual)$, and then basing the proximal penalty in \eqref{eq:intro:fb} on the square of the induced semi-norm $\norm{\freevar}_\Wave$.
We analyse an exemplifying choice of $\Wave$ based on convolution in \cref{sec:distances}. In particular we relate it to the weak-$*$ convergence of measures. To develop in \cref{sec:fb} a forward-backward method for \eqref{eq:intro:problem}, we will need $F'$ to be Lipschitz with respect to $\norm{\freevar}_\Wave$. For the squared data term, this reduces to $A_*A \le L \Wave$ for some factor $L \ge 0$. Based on approximation theory \cite{cheney2000approximation}, Bochner's theorem and Fourier transforms, we study this relationship in \cref{sec:sensorgrids}.

Our convergence proofs in \cref{sec:fb} are then standard and simple \cite{tuomov-proxtest,clasonvalkonen2020nonsmooth}.
However, we will need to account for inexact computation of each step of the method. The resulting approach bears more resemblance to the semi-infinite approach of \cite{flinth2020linear} than to conditional gradient methods in that more than one point may need to be inserted into the support of $\mu$ on each iteration of the method.
Inertial extensions, i.e., FISTA \cite{beck2009fista} for \eqref{eq:intro:problem}, are immediate.
We treat one in \cref{sec:fista}.

A forward-backward or conditional gradient method can only handle differentiable data terms with nonsmooth (Radon norm) regularisation.
However, our overall approach readily extends to other proximal-type type methods. We therefore sketch in \cref{sec:pdps} an extension of the primal dual proximal splitting (PDPS) of \cite{chambolle2010first} to problems of the form
\[
    \min_{0 \le \mu \in \Meas(\Omega)}\, F_0(A\mu) + \alpha \norm{\mu}_{\Meas},
\]
where $F_0$ may be nonsmooth.
Thus our overall approach is applicable to problems that are (so far) out of the reach of conditional gradient methods. As a further advantage, thanks to standard convergence proofs, our methods can readily be extended to product spaces $\Masses \times X$ for $X$ a Hilbert space.

We finish with numerical demonstrations in \cref{sec:numerical}.

\subsection*{Notation}

We denote the extended reals by $\extR \defeq [-\infty,\infty]$, and the space of finite Radon measures on a locally compact Borel measurable set $\Omega \subset \R^n$ by $\Meas(\Omega)$.
We write $\norm{\freevar}_{\Masses}$ or, for short, $\norm{\freevar}_\Meas$ for the Radon norm.
The subspace $\DiscreteMasses \subset \Meas(\Omega)$ consists of \term{discrete measures} $\mu=\sum_{k=1}^n \alpha_k \delta_{x_k}$ for any $n \in \N$, where the weights $\alpha_k \in \R$, and locations $x_k \in \Omega$.
Here, $\delta_x$ is the Dirac measure with mass one at a point $x \in \Omega$.
For $\mu \in \Meas(\Omega)$ and a Borel set $A \subset \Omega$, we write $\mu \llcorner A$ for the measure defined by $(\mu \llcorner A)(B) \defeq \mu(A \isect B)$ for Borel sets $B$.
We denote by $\B(x, r)$ the closed ball of radius $r>0$ and centre $x$, and by $\#K$ the cardinality of a finite set $K$.

For Fréchet differentiable $F:X \to R$ on a normed space $X$, we write $F'(x) \in X^*$ for the Fréchet derivative at $x \in X$.
Here $X^*$ is the dual space to $X$.
We call $F$ \term{pre-differentiable} if $F'(x) \in X_*$ for $X_*$ a designated \term{predual space} of $X$, satisfying $X=(X_*)^*$.
We have $(X_*)^{**}=X^*$, so $X_*$ canonically injects into $X^*$.
A predual of $\Meas(\Omega)$ is the space $C_c(\Omega)$ of continuous functions with compact support.
We write $C_0(\Omega) \defeq \closure C_c(\Omega)$ for continuous functions on $\Omega$ that vanish at infinity, which is also a predual space of $\Meas(\Omega)$ \cite[Theorem 1.200]{fonseca2007mmc}.

For a convex function $F: X \to \extR$, we write $\subdiff F: X \setto X_*$ for the (set-valued) pre-subdifferential map or, defined as $\subdiff F(x) = \{ x_*  \in X_* \mid F(\tilde x) - F(x) \ge \dualprod{x_*}{\tilde x - x} \text{ for all } \tilde x \in X\}$.
We write $\delta_C: X \to \extR$ for the $\{0,\infty\}$-valued indicator function of a set $C \subset X$.
Observe that the notation has similarities to the Dirac measure.

We write $\iprod{x}{x'}$ for the inner product between two elements $x$ and $x'$ of a Hilbert space $X$, and $\dualprod{x^*}{x} \defeq x^*(x)$ for the dual product in a Banach space.
We write $\linear(X; Y)$ for the space of bounded linear operators between two vector spaces $X$ and $Y$.
For $A \in \linear(X; Y)$, we write $\range A \subset Y$ for the range.
We let $\Id \in \linear(X; X)$ stand for the identity operator.
With $Y$ a Hilbert space for simplicity, we call $A \in \linear(X; Y)$ \term{pre-adjointable} if there exists a \term{pre-adjoint} $A_*: \linear(X_*; Y)$ whose mixed Hilbert–Banach adjoint $(A_*)^* = A$.
In other words $\dualprod{x}{A_*z} = \iprod{Ax}{z}$ for all $z \in Y$ and $x \in X$.

With $\imag=\sqrt{-1}$ the imaginary unit, we write $\fourier u$ or $\hat u$ for the Fourier transform, defined for of $u \in L^{1}(\R^n)$ by
\[
    [\fourier u](\xi) \defeq \int_{\R^n} u(x) e^{-2\pi\imag\iprod{x}{\xi}} \d x
    \quad (\xi \in \R^n),
\]
and for $u \in L^2(\R^n)$ by continuous extension; see \cite{rudin2006functional}.
Then the Hilbert-space adjoint $\fourier^*=\inv\fourier$.
We occasionally also take Fourier transforms of measures as tempered distributions.
We write $u * v$ for the convolution and $\autoconvolution u = u * u$ for the autoconvolution of functions $u, v$, recalling that $\fourier[u]^2 = \fourier[\autoconvolution[u]]$.

\section{Distances of measures}
\label{sec:distances}

As we have mentioned in the introduction, the norm or the squared norm is not particularly useful for deriving proximal-type methods in non-Hilbert spaces.
Nevertheless, let $M(\mu) \defeq \frac{1}{2}\norm{\mu}_{\Masses}^2$. The \term{preduality map} $\mu \mapsto \subdiff M(\mu)$, $\Meas(\Omega) \setto \FullPredual$ associates with every measure $\mu \in \Masses)$ a set of corresponding “wave functions”
\[
    \omega_\mu \in \subdiff M(\mu) = \{
        \omega\norm{\mu}_{\Masses}
        \mid
        \omega \in \FullPredual,\,
        -1 \le \omega \le 1,\,
        \dualprod{\omega}{\mu}=\norm{\mu}_{\Masses}
    \}.
\]
If this relationship were linear and single-valued, as it is in Hilbert spaces, $M(\mu-\this\mu)$ could be a practical proximal penalty term in \cref{eq:intro:fb}. Unfortunately, this is not the case in $\Masses$.
However, the concept of the preduality map suggests to define a distance on measures using some \term{particle-to-wave operator} $\Wave \in \linear(\Meas(\Omega); \FullPredual)$ that explicitly and linearly associates every measure, in the manuscript typically a discrete measure, with a corresponding wave function, an element of the predual.
In \cref{sec:distances:convolution} we will study such operators based on convolution.
Before this, in \cref{sec:distances:duality} we study general properties of seminorms and semi-inner products constructed with linear operators $\Wave$ to a predual space.

\subsection{Wave-particle duality}
\label{sec:distances:duality}

Let $X$ be a (real) Banach space with a predual space $X_*$, i.e., $X=(X_*)^*$.
Pick $\Wave \in \linear(X; X_*)$, which is \term{self-adjoint} and \term{positive semi-definite}, that is,
\begin{gather}
    \label{eq:wave:self-adjoint}
    \dualprod{\Wave x}{y}_{X_*,X} = \dualprod{x}{\Wave y}_{X,X_*}
    \quad (x, y \in X)
\shortintertext{and}
    \label{eq:wave:pos-def}
    \dualprod{\Wave x}{x}_{X_*,X} \ge 0
    \quad (x, x \in X).
\end{gather}
If the latter inequality is strict, we call $\Wave$ (strictly) \term{positive definite}.
Define
\[
    \iprod{x}{x}_{\Wave} \defeq \dualprod{\Wave x}{x}_{X_*,X}
    \quad\text{and}\quad
    \norm{x}_{\Wave} \defeq \sqrt{\iprod{x}{x}_{\Wave}}.
\]
Then it is easy to see that $\iprod{\freevar}{\freevar}_{\Wave}$ is a pseudo-inner product on $X$, and $\norm{\freevar}_{\Wave}$ a seminorm, i.e., non-negative and satisfies the triangle inequality.
Consequently, by a standard argument,\footnote{$\iprod{\mu}{\nu}_\Wave = \frac12\norm{\mu-\nu}_\Wave^2 - \frac12\norm{\mu}_\Wave^2 - \frac12\norm{\nu}_\Wave^2 \le \frac12\left(\norm{\mu}_\Wave+\norm{\nu}_\Wave\right)^2 - \frac12\norm{\mu}_\Wave^2 - \frac12\norm{\nu}_\Wave^2 = \norm{\mu}_\Wave\norm{\nu}_\Wave$.} it also satisfies the Cauchy–Schwarz inequality.
Moreover, the three-point identity holds:
\begin{equation}
    \label{eq:distances:wave-3-point}
    \begin{aligned}[t]
    \frac{1}{2}\norm{x-y}_{\Wave}^2
    &
    =
    \frac{1}{2}\dualprod{\Wave(x-y)}{x-y}
    \\
    &
    =
    \frac{1}{2}\dualprod{\Wave((x-z)+(z-y))}{(x-z)+(z-y)}
    \\
    &
    =
    \frac{1}{2}\norm{x-z}_{\Wave}^2
    +
    \frac{1}{2}\norm{z-y}_{\Wave}^2
    +
    \iprod{x-z}{z-y}_{\Wave}.
    \end{aligned}
\end{equation}

If $\Wave$ is injective, this construction turns $X$ into an inner product space. The (semi-)norm $\norm{\freevar}_{\Wave}$ may not, however, be equivalent to the original norm $\norm{\freevar}$ on $X$. Since $\Wave$ is assumed bounded, we always have $\norm{x}_{\Wave} \le \norm{\Wave} \norm{x}$, so the induced topology is weaker than the original (strong) topology. If the opposite inequality also holds, the norms are equivalent. For our constructions, this will, however, not be the case.

\subsection{Particle-to-wave operators based on convolution}
\label{sec:distances:convolution}

We now take $\Wave \mu = \rho * \mu$ for a convolution kernel $\rho \in C_0(\R^n)$.
We call $\rho$ \term{symmetric} if $\rho(-x)=\rho(x)$ for all $x \in \R^n$, and \term{positive semi-definite} if $0 \le \fourier[\rho]$.
We first show that $\Wave$ is well-defined.
Throughout, we take $\FullPredual$ as our designated predual space of $\Masses$.
When  $\Omega$ is compact, $\FullPredual=C(\Omega)$.

\begin{lemma}
    \label{lemma:wave:constr:basic}
    Let $\rho \in C_0(\R^n)$ be symmetric.
    On a closed domain $\Omega \subset \R^n$, let $\Wave \mu \defeq \rho * \mu$ for $\mu \in \Masses$.
    Then $\Wave \in \linear(\Meas(\Omega); \FullPredual)$ and is self-adjoint, i.e., satisfies \eqref{eq:wave:self-adjoint}.
\end{lemma}

\begin{proof}
    The definition of the convolution readily shows that $\Wave\mu \in C(\Omega)$.
    Let $\epsilon>0$.
    Since $\rho$ is continuous and vanishes at infinity, $\norm{\rho}_\infty < \infty$, and there exist $r>0$ such that $\abs{\rho(x)}\le \epsilon$ for $x \in \R^n \setminus \B(0, r)$.
    Consequently
    \[
        \abs{[\Wave\mu](y)}
        \le
        \adaptabs{\int_{\B(y, r)} \rho(y-x) \d\mu(x)}
        + \adaptabs{\int_{\Omega \setminus \B(y, r)} \rho(y-x) \d\mu(x)}
        \le
        \norm{\rho}_\infty \abs{\mu}(\B(y, r))
        + \epsilon \norm{\mu}_{\Meas(\Omega)}.
    \]
    Since $\norm{\mu}_{\Meas(\Omega)} < \infty$, for any $\epsilon>0$, we can find $R>0$ such that $\mu(\R^n \setminus \B(0, R)) \le \epsilon$.
    Thus $\abs{[\Wave\mu](y)} \le \epsilon(\norm{\rho}_\infty + \norm{\mu}_{\Meas(\Omega)})$ for $\norm{y} > r + R$.
    Since $\epsilon>0$ was arbitrary, this shows that $\Wave\mu \in C_0(\R^n)$.
    Since $\Omega$ is closed, it follows that  $\Wave\mu \in C_0(\Omega)$.
    Indeed, by the definition of $C_0(\Omega)$ as the closure of $C_c(\Omega)$, there exist $\phi_k \in C_c(\R^n)$ with $\phi_k \to \Wave\mu$ uniformly in $\R^n$.
    But for a closed $\Omega$, $\phi_k|\Omega \in C_c(\Omega)$ with $\phi_k|\Omega \to (\Wave\mu)|\Omega$ uniformly in $\Omega$.

    We then expand and rearrange using Fubini's theorem for any $\nu \in \Masses$ that
    \begin{equation}
        \label{eq:wave:constr:expand-dualprod}
        \begin{aligned}[t]
        \dualprod{\Wave\mu}{\nu}_{\FullPredual,\Meas(\Omega)}
        &
        =
        \int_{\Omega} [\Wave\mu](x) \d\nu(x)
        =
        \int_{\Omega} \int_\Omega \rho(x-y) \d\mu(y) \d\nu(x)
        \\
        &
        =
        \int_{\Omega} \int_\Omega \rho(x-y) \d\nu(x) \d\mu(y)
        =
        \int_{\Omega} \int_\Omega \rho(y-x) \d\nu(x) \d\mu(y)
        \\
        &
        =
        \int_{\Omega} [\rho * \nu](y) \mu(y)
        =
         \dualprod{\Wave\nu}{\mu}_{\FullPredual,\Meas(\Omega)}.
        \end{aligned}
    \end{equation}
    
    Hence $\norm{\rho * \nu}_\infty \le \norm{\rho}_\infty \norm{\nu}_{\Meas}$; see, e.g., \cite[§2.1]{ambrosio2000fbv}.
    Consequently \eqref{eq:wave:constr:expand-dualprod} establishes
    $
        \dualprod{\Wave\mu}{\nu}_{\FullPredual,\Meas(\Omega)}
        \le  \norm{\rho}_\infty \norm{\nu}_{\Meas}\norm{\mu}_{\Meas}.
    $
    This shows that $\Wave$ is bounded. It is evidently also linear so that $\Wave \in \linear(\Meas(\Omega); \FullPredual)$ as claimed.
    From \eqref{eq:wave:constr:expand-dualprod} we also observe that $\Wave$ is self-adjoint.
\end{proof}

We now characterise weak-$*$ convergence in terms of $\Wave$-seminorms.
Following \cite[Chapter 18]{cheney2000approximation} we call a subset $V$ of a normed space $X$ \term{fundamental} if $\closure \Span V$ is dense in $X$. We are interested in the fundamentality in $\FullPredual$ of
\[
    V_\rho(\Omega) \defeq \{ x \mapsto \rho(x - y) \mid y \in \Omega \}.
\]
According to \cite[Theorem 18.1]{cheney2000approximation}, this is the case when $\fourier[\rho]$ is a finite-valued real non-negative measure and satisfies for all non-empty open sets $U \subset \R^n$ the strict lower bound
\begin{equation}
    \label{eq:wave:fundamentality-condition}
    \int_U \fourier[\rho](\xi) \d \xi > 0.
\end{equation}

\begin{example}
    Since the Fourier transform of a Gaussian $\rho$ is a Gaussian, $V_\rho(\Omega)$ is fundamental in $C(\Omega)$ for any compact $\Omega \subset \R^n$.
\end{example}

We will require the following lemma on a few occasions.

\begin{lemma}
    \label{lemma:distances:uniform-equicontinuity}
    Suppose $\phi \in C_0(\R^n)$,
    the domain $\Omega \subset \R^n$,
    and $\{\this\mu\}_{k \in \N} \subset \Masses$ is bounded.
    Then the collection $\{x \mapsto [\phi * \this\mu](x)\}_{k \in \N}$ is uniformly equicontinuous on $\closure \Omega$.
\end{lemma}

\begin{proof}
    
    Let $\epsilon>0$ and
    $C \defeq \sup_{k \in \N} \norm{\this\mu}_{\Masses}$.
    By assumption $C < \infty$.
    Since $\phi \in C_0(\R^n)$, it is uniformly continuous. Therefore, there exists $\delta>0$ such that $\abs{\phi(x)-\phi(y)} \le \epsilon$ whenever $\norm{x-y} \le \delta$.
    Now
    \[
        \abs{[\phi * \this\mu](x) - [\phi * \this\mu](y)}
        \le
        \int_{\Omega} \abs{\phi(x-z)-\phi(y-z)} \d\abs{\this\mu}(z)
        \le
        \int_{\Omega} \epsilon \d\abs{\this\mu}(z)
        \le
        C \epsilon
    \]
    for any $k \in \N$ and $x, y \in \Omega$ with $\norm{x-y} \le \delta$. Since the right hand side $C \epsilon$ does not depend on $k$, and can be made arbitrarily small for sufficiently small $\epsilon$ and corresponding $\delta$, this proves the claim.
\end{proof}

Now we can state and prove our main result on particle-to-wave operators based on convolution.
Recall that we write $\autoconvolution u \defeq u * u$ for the autoconvolution.

\begin{theorem}
    \label{thm:wave:constr:weak-new}
    Let $0 \not\equiv \rho \in C_0(\R^n) \isect L^2(\R^n)$ be symmetric and positive semi-definite.
    On a closed domain $\Omega \subset \R^n$, let $\Wave \in \linear(\Masses; \FullPredual)$ be defined by $\Wave \mu = \rho * \mu$ for $\mu \in \Meas(\Omega)$. Then
    \begin{enumerate}[label=(\roman*)]
        \item\label{item:wave:constr:pos-def}
        $\Wave$ is self-adjoint and positive semi-definite.
        \item\label{item:wave:constr:strict-pos-def}
        $\Wave$ is (strictly) positive definite on the discrete measures $\DiscreteMasses$.
    \end{enumerate}
    Assume further that $\Omega$ is compact, and that there exists $\rho^{1/2} \in L^2(\R^n) \isect C_0(\R^n)$ such that $\rho=\autoconvolution[\rho^{1/2}]$. Then
    \begin{enumerate}[resume*]
        \item\label{item:wave:constr:weak-to-wave-new}
         $\this\mu \weaktostar \mu$ weakly-$*$ in $\Meas(\Omega)$ implies $\norm{\this\mu-\mu}_\Wave \to 0$.
        \item\label{item:wave:constr:wave-to-weak-new}
        If  $V_{\rho^{1/2}}$ is fundamental for $\FullPredual$, then $\norm{\this\mu-\mu}_\Wave \to 0$ with $\{\this\mu\}_{k \in \N} \subset \Meas(\Omega)$ bounded implies $\this\mu \weaktostar \mu$ weakly-$*$ in $\Meas(\Omega)$.
    \end{enumerate}
\end{theorem}

\begin{proof}
    We know from \cref{lemma:wave:constr:basic} that $\Wave \in \linear(\Masses; \FullPredual)$, and is self-adjoint.
    Bochner's theorem  proves \cref{item:wave:constr:pos-def}, and \cref{item:wave:constr:strict-pos-def} follows from \cite[Theorem 13.3]{cheney2000approximation}. However, due to slightly different assumptions for the latter, we find it easiest to provide the proofs.
    We write $\hat\rho=\fourier[\rho]$.
    Since $\rho \in  L^2(\R^n)$ is symmetric, also $\hat\rho \in L^2(\R^n)$ is symmetric.
    Therefore, by Fubini's theorem and the Fourier inverse transform, we have
    \begin{equation}
        \label{eq:wave:constr:posdef}
        \begin{aligned}[t]
        \dualprod{\Wave \mu}{\mu}
        &
        =
        \int_{\Omega} \int_{\Omega} \rho(x-y) \d\mu(x) \d\mu(y)
        \\
        &
        =
        \int_{\Omega} \int_{\Omega} \int_{\R^n} \hat\rho(\xi) e^{2\pi\imag\iprod{\xi}{x-y}} \d\xi \d\mu(x) \d\mu(y)
        \\
        &
        =
        \int_{\R^n} \hat\rho(\xi) \left| \int_{\Omega} e^{-2\pi\imag\iprod{\xi}{x}} \d\mu(x) \right|^2 \d\xi
        \\
        &
        =
        \int_{\R^n} \hat\rho(\xi) |\fourier[\mu](\xi)|^2 \d\xi.
        \end{aligned}
    \end{equation}
    Since $\hat\rho \ge 0$ by our definition of the positive semi-definiteness of $\rho$, this shows \cref{item:wave:constr:pos-def}.

    To prove \ref{item:wave:constr:strict-pos-def}, we start by observing that since $\rho \in L^2(\R^n)$ is not identically zero and has $\hat\rho \ge 0$, necessarily $\int \hat\rho \d\xi > 0$. But then there exists $\epsilon>0$ and a Borel set $E$ with nonzero Lebesgue measure such that $\hat\rho \ge \epsilon$ on $E$. When $\mu \in \DiscreteMasses$, the zero set of the continuous function $\fourier[\mu]$, a finite sum of complex exponentials, has Lebesgue measure zero by \cite[Lemma 13.6]{cheney2000approximation}.
    Therefore $\int_{\R^n} \hat\rho(\xi) |\fourier[\mu](\xi)|^2 \d\xi \ge \epsilon \int_E |\fourier[\mu](\xi)|^2 \d\xi > 0$.
    The claim now follows from \eqref{eq:wave:constr:posdef}.
    
    For \cref{item:wave:constr:weak-to-wave-new} we observe that $\Wave = \WaveSqr_*\WaveSqr$ for $\WaveSqr \in \linear(\Meas(\Omega); L^2(\R^n))$ defined by $\WaveSqr \mu \defeq \rho^{1/2} * \mu$ and $\WaveSqr_* f \defeq \rho^{1/2} * f$ for $f \in L^2(\R^n)$ and $\mu \in \Masses$.
    We therefore need to prove that $\rho^{1/2} * (\this\mu-\mu) \to 0$ in $L^2(\R^n)$.
    Observe that $\rho_y^{1/2}(x) \defeq \rho^{1/2}(x - y)$ satisfies $\rho_y^{1/2} \in \FullPredual$ for all $y \in \R^n$.
    (For this we need the closedness of $\Omega$, but not yet the boundedness; compare the proof of \cref{lemma:wave:constr:basic}.)
    Thus $\dualprod{\rho_x^{1/2}}{\this\mu-\mu} \to 0$  for all $x \in \R^n$ by the assumed weak-$*$ convergence.
    %
    Since by assumption $\rho^{1/2} \in L^2(\R^n)$, we can for each $\epsilon>0$ find a bounded $U_\epsilon \subset \R^n$ such that
    $\int_{\R^n \setminus U_\epsilon} \rho^{1/2}(x)^2 \d x \le \epsilon  \norm{\rho^{1/2}}_\infty^2$.
    Letting $\Theta_\epsilon \defeq U_\epsilon + \Omega$, for all $y \in \Omega$ then
    \[
        \int_{\R^n \setminus \Theta_\epsilon} \rho_y^{1/2}(x)^2 \d x
        =
        \int_{\R^n \setminus (\Theta_\epsilon - y)} \rho^{1/2}(x)^2 \d x
        \le
        \int_{\R^n \setminus U_\epsilon} \rho^{1/2}(x)^2 \d x \le \epsilon  \norm{\rho^{1/2}}_\infty^2.
    \]
    Using Jensen's inequality and Fubini's theorem, then
    \[
        \begin{aligned}[t]
        \int_{\R^n \setminus \Theta_{\epsilon}} [\rho^{1/2} * (\this\mu-\mu)](x)^2 \d x
        & =
        \int_{\R^n \setminus \Theta_\epsilon}\left(\int_{\Omega} \rho_y^{1/2}(x) \d \abs{\this\mu-\mu}(y)\right)^2 \d x
        \\
        &
        \le
        \norm{\this\mu-\mu}_{\Meas} \int_{\R^n \setminus \Theta_\epsilon}\int_{\Omega} \rho_y^{1/2}(x)^2 \d \abs{\this\mu-\mu}(y) \d x
        \\
        &
        =
        \norm{\this\mu-\mu}_{\Meas} \int_{\Omega} \int_{\R^n \setminus \Theta_\epsilon} \rho_y^{1/2}(x)^2  \d x \d \abs{\this\mu-\mu}(y)
        \\
        &
        \le
        \epsilon \norm{\rho^{1/2}}_\infty^2 \norm{\this\mu-\mu}_{\Meas}^2.
        \end{aligned}
    \]

    The set $\Theta_\epsilon$ is of finite measure, because $\Omega$ is bounded.
    By Egorov's theorem we can thus find a set $E_\epsilon \subset \Theta_{\epsilon}$ such that the Lebesgue measure $\Lebesgue(\Theta_{\epsilon} \setminus E_{\epsilon}) \le \epsilon$ and $x \mapsto \dualprod{\rho_x^{1/2}}{\this\mu-\mu}$ converges in $E_{\epsilon}$ uniformly to zero.
    Then
    \[
        \begin{aligned}[t]
        \int_{\R^n} [\rho^{1/2} * (\this\mu-\mu)](x)^2 \d x
        &
        =
        \int_{E_\epsilon} \dualprod{\rho_x^{1/2}}{\this\mu-\mu}^2 \d x
        +
        \int_{\Theta_{\epsilon} \setminus E_\epsilon} [\rho^{1/2} * (\this\mu-\mu)](x)^2 \d x
        \\
        \MoveEqLeft[-1]
        +
        \int_{\R^n \setminus \Theta_{\epsilon}} [\rho^{1/2} * (\this\mu-\mu)](x)^2 \d x
        \\
        &
        \le
        \int_{E_\epsilon} \dualprod{\rho_x^{1/2}}{\this\mu-\mu}^2 \d x
        + 2 \epsilon \norm{\rho^{1/2}}_\infty^2\norm{\this\mu-\mu}_{\Meas}^2.
        \end{aligned}
    \]
    Thus $\limsup_{k \to \infty} \norm{\rho^{1/2} * (\this\mu-\mu)}_{L^2(\R^n)} \le \sqrt{2\epsilon}\norm{\rho^{1/2}}_\infty\sup_k \norm{\this\mu-\mu}_{\Meas}$.
    Since $\epsilon>0$ was arbitrary, and $\{\this\mu-\mu\}_{k \in \N}$ is bounded by weak-$*$ convergence, it follows that $\rho^{1/2} * (\this\mu-\mu) \to 0$ in $L^2(\R^n)$.
    
    For \cref{item:wave:constr:wave-to-weak-new} we note that since $V_{\rho^{1/2}}$ is by assumption fundamental for $\FullPredual$, we can for any $\phi \in \FullPredual$ and $\epsilon>0$ find a finite set $U_\epsilon \subset \Omega$ and coefficients $\{\alpha_y\}_{y \in U_\epsilon} \subset \R$ such that
    \[
        \sup_{x \in \Omega} \Biggl|\phi(x) - \sum_{y \in U_\epsilon} \alpha_y \rho^{1/2}_y(x) \Biggr| < \epsilon.
    \]
    Then, for all $k \in \N$,
    \begin{equation}
        \label{eq:distance:convolution:weak-approximation}
        \dualprod{\phi}{\this\mu-\mu}
        =
        \sum_{y \in U_\epsilon} \dualprod{\alpha_y\rho^{1/2}_y}{\this\mu-\mu}
        +
        \Biggl\langle \phi - \sum_{y \in U_\epsilon} \alpha_y \rho^{1/2}_y \Biggm| \this\mu-\mu \Biggr\rangle.
    \end{equation}
    The assumption $\norm{\this\mu-\mu}_\Wave \to 0$ readily implies $\rho^{1/2} * (\this\mu-\mu) \to 0$ in $L^2(\R^n)$.
    Hence, there exists a subsequence $\{k_j\}_{j \in \N}$ of $\{k\}_{k \in \N}$ such that $\rho^{1/2} * (\mu^{k_j}-\mu) \to 0$ almost everywhere.
    Since, by the Heine--Borel theorem,  $\Omega$ is bounded, and $\sup_{k \in \N} \norm{\this\mu-\mu}_{\Masses} \le M$ for some $M > 0$,
    \cref{lemma:distances:uniform-equicontinuity} shows that the family $\{ \rho^{1/2} * (\this\mu-\mu)\}_{k \in \N}$ is uniformly equicontinuous on $\Omega$.
    By a compact covering argument, it follows that $\rho^{1/2} * (\mu^{k_j}-\mu) \to 0$ uniformly on $\Omega$.
    In particular, $\abs{\dualprod{\alpha_y\rho^{1/2}_y}{\mu^{k_j}-\mu}} < \epsilon/\#U_\epsilon$ for all $y \in U_\epsilon$ for $j \in \N$ large enough.
    It follows from \eqref{eq:distance:convolution:weak-approximation} that
    \begin{equation}
        \label{eq:distance:convolution:weak-epsilon}
        \abs{\dualprod{\phi}{\mu^{k_j}-\mu}} \le (1 + M) \epsilon
        \quad\text{for $j$ large enough.}
    \end{equation}
    Since $\epsilon>0$ and $\phi \in \FullPredual$ were arbitrary, we deduce $\mu^{k_j} \weaktostar \mu$.
    If there would exist a subsequence $\{\mu^{k_\ell}\}_{\ell \in \N}$ of $\{\mu^k\}_{k \in \N}$ not weakly-$*$ convergent to $\mu$, then we could find a further unrelabelled subsequence, $\epsilon>0$, and $\phi \in \FullPredual$ such that $\abs{\dualprod{\phi}{\mu^{k_\ell}-\mu}} \ge 2 (1 + M) \epsilon$ for all $\ell \in \N$. Repeating the above arguments, we would obtain \eqref{eq:distance:convolution:weak-epsilon}, hence a contradiction.
    Therefore $\mu^k \weaktostar \mu$.
\end{proof}

\begin{corollary}
    \label{cor:distances:norm}
    Let all the assumptions of \cref{thm:wave:constr:weak-new}\,\cref{item:wave:constr:wave-to-weak-new} hold.
    Then $\norm{\freevar}_\Wave$ is a norm on $\Masses$.
\end{corollary}

\begin{proof}
    Based on the discussion that in the beginning of \cref{sec:distances:duality}, it only remains to prove that $\norm{\mu-\nu}_\Wave = 0$ implies $\mu=\nu$.
    Indeed, setting $\mu^k=\mu$ for all $k \in \N$, we have $\norm{\mu^k-\nu}_\Wave \to 0$.
    Then \cref{thm:wave:constr:weak-new}\,\cref{item:wave:constr:wave-to-weak-new} implies $\mu = \mu^k \weaktostar \nu$. Thus $\mu=\nu$.
\end{proof}

\section{Convolution kernels and sensor grids}
\label{sec:sensorgrids}

We will need the data term $F(\mu) = \frac{1}{2}\norm{A\mu-b}^2$ to satisfy a smoothness property also known as the “descent lemma” with respect to the wave-particle norm $\norm{\freevar}_\Wave$. As shown in, e.g., \cite[Lemma 7.1]{clasonvalkonen2020nonsmooth}, this in general follows from $F'$ being Lipschitz. In our case the Lipschitz property needs to be with respect to the $\Wave$-seminorm.
For quadratic $F$ as in \eqref{eq:intro:problem}, this requirement reduces to
\[
    A_*A \le L \Wave \quad\text{for some}\quad L>0.
\]
This means that we have to find $\Wave$ such that  $L \Wave - A_*A$ is positive semi-definite.
We could take $\Wave=A_*A$ and $L=1$ and still obtain a theoretical algorithm with convergent function values, although in general not the weak-$*$ convergence of iterates.
However, the proximal step of the method that we develop in the next section consists of finding new points $x$ to insert into the support of $\nexxt\mu$ to approximately satisfy first-order optimality conditions for \eqref{eq:intro:fb} with $\frac{1}{2}\norm{\freevar}_\Wave$ as the proximal penalty.
To avoid such an insertion from perturbing optimality conditions globally, it will be practical if $\Wave\delta_x$ is localised around $x$. This can be achieved when $\Wave$ is a simple convolution operator.
However, even for physical processes that can be modelled as convolutions, $A_*A$ is in general not \emph{exactly} a convolution operator when the range of $A$ consists of measurements on a \emph{finite} sensor grid, and its effects are not localised: $A^*A\delta_x$ may have multiple peaks.

We will, therefore, first in \cref{sec:sensorgrid:sensorgrid} construct upper approximations by convolution operators $\Wave$ for forward operators $A$ modelling measurements of convolved signals on a finite sensor grid.
We provide specific examples in \cref{sec:sensorgrid:examples}.
Although we work with convolution operators in this section, in the following sections on our proposed algorithms, there is in general no need for $A$ or $\Wave$ to be based on convolution.

\subsection{Sensor grids}
\label{sec:sensorgrid:sensorgrid}

\def\grid{\mathscr{G}}
\def\shift{T}

We now construct some example forward operators $A$, based on a simple physical spreading process followed by sensing. More precisely, the \term{spread} $\psi \in C_c(\R^n)$ is a convolution kernel that models the physical spreading process of, e.g., light. If a source of unit intensity is at $x$, then $\psi(y-x)$ is the intensity that falls at $y$. This light is sensed by the \term{sensor} $\theta_z \in L^2(\R^n) \isect L^1(\R^n)$ placed at $z$ on a finite \term{grid} $\grid$ with the sensitivity $\theta_z(y)$. With this we define the forward operator $A$ as
\begin{equation}
    \label{eq:sensorgrid:sensorgrid}
    A\mu \defeq (\mu(\theta_z * \psi))_{z \in \grid},
    \quad\text{i.e.,}\quad
    [A\mu]_z = \int (\theta_z * \psi)(x) \d\mu(x)
    \quad\text{for}\quad
    z \in \grid.
\end{equation}
 The sensors or, more precisely, the field-of-views of the sensors, are a technical necessity for our proof. They can, nevertheless, also be physically justified through practical light sensors having positive area instead of just sensing the light falling at a singular point.
For $A$ of this form, the next lemma establishes an upper bound on $A_* A$ in terms of a simple convolution operator.

\begin{lemma}
    \label{lemma:sensorgrid:base-spread-construction}
    Let the \term{sensor grid} $\grid$ be an arbitrary finite index set with corresponding \term{sensors} $0 \le \theta_z \in L^2(\R^n) \isect L^1(\R^n)$, ($z \in \grid$), satisfying for some $L_0 \ge 0$ the bound
    \begin{equation}
        \label{eq:convolution:fov-assumption}
        L_0 \ge \sup_{x \in \R^n} \sum_{z \in \grid} \theta_z(x) \int_{\R^n} \theta_z(w) \d w.
    \end{equation}
    Pick a symmetric \term{spread} $\psi \in C_0(\R^n) \isect L^2(\R^n)$, and on a domain $\Omega \subset \R^n$, define $A \in \linear(\Meas(\Omega); \R^{\#\grid})$ by \eqref{eq:sensorgrid:sensorgrid}.
    Then the pre-adjoint $A_* \in \linear(\R^n; \FullPredual)$ exists and is given by $A_* y = \sum_{z \in \grid} y_z \theta_z * \psi$, and we have
    \[
        A_* A \le L_0 \dualprod{\freevar}{\autoconvolution[\psi] * \freevar}.
    \]
\end{lemma}

\begin{proof}
    We easily verify that $(A_*)^*=A$. $A_*$ is also bounded since $A^* = J A_*$ where $J$ is the canonical injection $\FullPredual \hookrightarrow \Meas(\Omega)^*$.
    We have
    \[
        \dualprod{\mu}{A_*A\mu}
        = \norm{A\mu}_2^2
        = \sum_{z \in \grid} \left( \int \int_{\R^n} \theta_z(w) \psi(x-w) \d w \d\mu(x) \right)^2
    \]
    Recalling that $\theta_z \ge 0$, defining $I_{\theta_z} \defeq \int_{\R^n} \theta_z(w) \d w \in [0, \infty)$, and using Hölder's inequality, we continue
    \[
        \begin{aligned}
        \dualprod{\mu}{A_*A\mu}
        &
        = \sum_{z \in \grid}  \left(\int_{\R^n} \int  \psi(x-w) \d\mu(x) \theta_z(w) \d w \right)^2
        \\
        &
        \le
        \sum_{z \in \grid}  \int_{\R^n} \theta_z(w) \d w \int_{\R^n} \left( \int \psi(x-w) \d\mu(x) \right)^2 \theta_z(w) \d w.
        \end{aligned}
    \]
    Using \eqref{eq:convolution:fov-assumption} and the symmetricity of $\psi$, it follows
    \[
        \begin{aligned}
            \dualprod{\mu}{A_*A\mu}
            &
            \le
            L_0 \int_{\R^n} \left( \int \psi(x-w) \d\mu(x) \right)^2 \d w
            \\
            &
            =
            L_0 \int_{\R^n} \int \psi(x-w) \d\mu(x) \int \psi(y-w) \d\mu(y) \d w
            \\
            &
            =
            L_0 \int \int \int_{\R^n} \psi(x-w)\psi(w-y) \d w \d\mu(y) \d\mu(x)
            \\
            &
            =
            L_0 \dualprod{\mu}{\autoconvolution[\psi] * \mu}.
        \end{aligned}
    \]
    This establishes the claim.
\end{proof}

\begin{example}[Rectangular sensors]
    \label{ex:convolution:rectangular-sensor}
    Let $\theta_z(x) = \chi_{[-b,b]^n}(x-z)$ for some $0 < b < \inf_{z, \tilde z \in \grid} \norm{z-\tilde z}_\infty$.
    Then \eqref{eq:convolution:fov-assumption} is satisfied by $L_0 = (2b)^n$.
\end{example}

Having bounded $A_* A$ by $\dualprod{\freevar}{\autoconvolution[\psi] * \freevar}$, we now need to bound $\dualprod{\freevar}{\autoconvolution[\psi] * \freevar}$ by $\Wave$ to establish the overall bound $A_*A \le L \Wave$. This is done in the next theorem.

\begin{theorem}[Putting it all together]
    \label{thm:convolution:together0}
    Suppose the “spread” and the “kernel” $\psi, \rho \in C_0(\R^n) \isect L^2(\R^n)$ satisfy
    \begin{equation}
        \label{eq:convolution:bochner-condition}
        \fourier[\psi]^2 \le L_1 \fourier[\rho]
    \end{equation}
    Also assume $0 \le \theta_z \in L^2(\R)$ for all $z \in \grid$ on a finite grid $\grid \subset \R$ with $L_0$ satisfying \eqref{eq:convolution:fov-assumption}. 
    For a closed $\Omega \subset \R^n$, define
    \begin{equation}
        \label{eq:convolution:together:defs}
        A\mu \defeq (\mu(\theta_z * \psi))_{z \in \grid}
        \quad\text{and}\quad
        \Wave\mu \defeq \rho * \mu
        \quad\text{for}\quad
        \mu \in \Meas(\Omega).
    \end{equation}
    Then a pre-adjoint $A_* \in \linear(\R^n; \FullPredual)$ is given by $A_* y = \sum_{z \in \grid} y_z \theta_z * \psi$, and we have
    \begin{equation}
        \label{eq:convolution:together0}
        A_*A \le L_0 L_1 \Wave.
    \end{equation}
\end{theorem}

\begin{proof}
    Due to \cref{lemma:sensorgrid:base-spread-construction}, we only need to show that $\dualprod{\freevar}{\autoconvolution[\psi] * \freevar} \le L_1 \Wave$.
    By Bochner's theorem (see \cite[Chapter 12]{cheney2000approximation} or the proof of \cref{thm:wave:constr:weak-new}), and the Fourier transform convolution-multiplication exchange rule (see, e.g., \cite[Theorem 7.19]{rudin2006functional}), this amounts to \eqref{eq:convolution:bochner-condition}.
\end{proof}

\subsection{Examples}
\label{sec:sensorgrid:examples}

One choice satisfying \cref{thm:convolution:together0} with $L_1=1$ would be  $\rho=\autoconvolution[\psi]$.
This could also easily be made to satisfy all the additional conditions of \cref{thm:wave:constr:weak-new} with $\rho^{1/2}=\psi$, so would be an option if we were able to calculate $\autoconvolution[\psi]$ efficiently.
For numerical reasons, both $\rho$ and $\psi$ should also have small support, which is not the case for a Gaussian $\psi$.
Therefore, we would need to cut the spread to ensure small compact support.
The next example, on the other hand, avoids even forming the autoconvolution by taking $\rho=\psi$. Why we call the choice “fast” will be apparent once stumble onto Gauss error functions when treating a cut Gaussian spread.

To proceed, we define for $b>0$ the hat function
\begin{equation}
    \label{eq:wave:tri}
    \tri_b(t) \defeq \frac{1}{b}\autoconvolution[\chi_{[-b/2,b/2]}](t)
    = \begin{cases}
        1-\abs{t}/b, & t \in (-b, -b),
        \\
        0, & \text{otherwise}.
    \end{cases}
\end{equation}
We also recall for $a>0$ that (see, e.g., \cite[Appendix 3]{kammler2008first})
\begin{equation}
    \label{eq:wave:indicator-fourier}
    \fourier[\chi_{[-a,a]}] = \fourier[\rect(\tfrac{1}{2a}\freevar)] = 2a\sinc(2\pi a\freevar).
\end{equation}
for the (unnormalised) $\sinc t \defeq \sin t/t$.

\begin{example}[“Fast” spread]
    \label{ex:convolution:fast}
    For some $\sigma>0$, let $\rho=\psi$ and for $x \in \R$,
    \[
        \psi(x)
        = \frac{4}{\sigma} \autoconvolution[\tri_{1/2}] (x/\sigma)
        =  \frac{4}{\sigma}
        \begin{cases}
            \frac{2}{3} (x/\sigma+1)^3 & -1 < x/\sigma \le -\frac{1}{2}, \\
            2 \abs{x/\sigma}^3-2 (x/\sigma)^2+\frac{1}{3} & -\frac{1}{2} < x/\sigma \le \frac{1}{2}, \\
            -\frac{2}{3} (x/\sigma-1)^3 & \frac{1}{2} \le x/\sigma < 1, \\
            0 & \text{otherwise}.
        \end{cases}
    \]
    Also let $\theta_z(x)=\chi_{[-b,b]}(x-z)$ for some $0 < b < \inf_{z, \tilde z \in \grid} \abs{z-\tilde z}$ on the finite grid $\grid \subset \R$, and define $A$ and $\Wave$ according to \eqref{eq:convolution:together:defs}.
    Clearly $\supp\psi=\supp\rho=[-\sigma,\sigma]$.
    It can also be verified that $\int \psi \d x=1$.
    Writing $\hat\psi \defeq \fourier[\psi]$, we expand using standard Fourier transform dilation and convolution rules (see, e.g., \cite[Appendix 3]{kammler2008first}), \eqref{eq:wave:tri}, and \eqref{eq:wave:indicator-fourier} that
    \[
        \begin{aligned}
        \fourier[\psi](\xi)
        &
        = 4 \fourier[\autoconvolution[\tri_{1/2}]](\sigma \xi)
        = 4 \fourier[\tri_{1/2}](\sigma \xi)^2
        \\
        &
        = 16 \fourier[\chi_{[-1/4, 1/4]}](\sigma \xi)^4
        = 16 \left(\frac{1}{2}\sinc(\pi\sigma\xi/2) \right)^4
        = \sinc(\pi\sigma\xi/2)^4
        \quad (\xi \in \R).
        \end{aligned}
    \]
    Since $\sinc$ achieves its maximum at $\xi=0$ with value 1, it follows that $\fourier[\psi]^2 \le \fourier[\rho]$.
    Thus we may apply \cref{thm:convolution:together0} with $L_1=1$ and $L_0$ given by \cref{ex:convolution:rectangular-sensor}. This gives the bound
    \[
        A_*A \le 2b \Wave.
    \]
    
    Obviously $\rho^{1/2} = \frac{2}{\sigma}\tri_{1/2}(\freevar/\sigma) \in L^2(\R^n) \isect C_c(\R^n)$.
    The above calculations establish that $\int_U \fourier[\rho^{1/2}](\xi) \d\xi = \int_U \sinc(\pi\sigma\xi/2)^2 \d\xi > 0$ for every open set $U$, so that \eqref{eq:wave:fundamentality-condition} holds for $\rho^{1/2}$.
    Therefore \cite[Theorem 18.1]{cheney2000approximation} proves that $V_{\rho^{1/2}}(\Omega)$ is fundamental on compact sets $\Omega \subset \R^n$.
    Now \cref{thm:wave:constr:weak-new} shows that $\norm{\freevar}_\Wave$ characterises weak-$*$ convergence on such sets.

    Finally, to apply $A$, we need to be able to calculate $\theta_z * \psi = \chi_{[-b, b]} * v$.
    This is a piecewise polynomial, so easily implemented in software, although sizeable to write down case-by-case.
\end{example}

To base the spread $\psi$ on a Gaussian, we need to cut it to have a compact support.
We therefore next construct $\psi = \chi v$ for a cut-off function $\chi$ and a base spread $v$.
Since $\rho = \autoconvolution[\psi]$ can be numerically unwieldy, we also take $\rho=\phi u$ for
$\phi=\autoconvolution[\chi]$ and a computationally simple base kernel $u$, in our numerical practise $u=v$.
The condition \eqref{eq:convolution:compare:assumption} roughly says that the base spread $v$ cannot be more localised in space than the base kernel $u$.

In the following, we use the abbreviations $\hat u = \fourier[u]$, $\hat v = \fourier[v]$, etc., without explicit mention.

\begin{lemma}
    \label{lemma:convolution:compare}
    Let $u,v,\chi \in L^2(\R^n)$.
    Then \eqref{eq:convolution:bochner-condition} holds if $\hat v \ge 0$ with $\int \hat v \d\xi < \infty$, and, for some $L_1 \ge 0$ we have
    \begin{equation}
        \label{eq:convolution:compare:assumption}
        I_{\hat v}
        \hat v
        \le
        L_1 \hat u.
    \end{equation}
\end{lemma}

\begin{proof}
    Since $\fourier[\autoconvolution(\chi)u](\xi) = [\hat\chi^2 * \hat u](\xi)$ by standard Fourier transform exchange rules for convolution and multiplication, we need to prove for all $\xi \in \R^n$ that
    $
        \fourier[\chi v](\xi)^2
        \le
        [\hat\chi^2 * \hat u](\xi).
    $
    If $I_{\hat v} \defeq \int \hat v \d\xi = 0$, also $v=0$, so this is immediate by the assumed non-negativity of $\hat u$. So suppose $I_{\hat v}>0$.
    Since, by assumption $I_{\hat v} < \infty$ and $\hat v \ge 0$, we have $\hat v \in L^1(\R^n)$.
    Now Hölder's inequality shows as claimed that
    \[
        (\hat\chi * \hat v)(\xi)^2
        = \left(
            \int_{\R^n} \hat\chi(z) \hat v(\xi-z) \d z
        \right)^2
        \le
            \left(\int_{\R^n} \hat v(\xi-z)\d z\right)
            \left(\int_{\R^n} \hat\chi(z)^2 \hat v(\xi-z) \d z\right)
        =
        I_{\hat v} [\hat\chi^2 * \hat v](\xi).
        \qedhere
    \]
\end{proof}

If we want $\norm{\freevar}_\Wave$ defined with the kernel $\rho=\phi u$ to characterise weak convergence we need to show for \cref{thm:wave:constr:weak-new} the existence of $\rho^{1/2} \in L^2(\R^n) \isect C_0(\R^n)$, and that $V_{\rho^{1/2}}(\Omega)$ is fundamental. Unfortunately, we have no simple characterisation of the requirement $\rho^{1/2} \in C_0(\R^n)$, although it would hold by properties of the Fourier transform if $\rho$ were a rapidly decreasing function (i.e., a test function for tempered distributions) \cite[Chapter 7]{rudin2006functional}.

\begin{lemma}
    \label{lemma:wave:cut-square-root}
    Suppose $\rho=\phi u$ for symmetric $\phi, u \in L^2(\R^n)$ such that $ \hat\phi, \hat u \in L^1(\R^n)$ and $\hat\phi, \hat u \ge 0$. Then there exists $\rho^{1/2} \in L^2(\R^n)$ such that $\rho=\autoconvolution[\rho^{1/2}]$.
    If $\hat\rho(\freevar)^{1/2} \in L^1(\R^n)$, then, moreover, $\rho^{1/2} \in C_0(\R^n)$.
    
\end{lemma}

\begin{proof}
    We have $\hat \rho = \hat\phi * \hat u \ge 0$ with all of the Fourier transforms real-valued and symmetric by the corresponding properties of $\phi$ and $u$.
    Since $\hat \phi, \hat u \in L^1(\R^n)$, standard properties of convolutions imply $\hat\rho \in L^1(\R^n)$.
    Therefore $\hat\rho(\freevar)^{1/2} \in L^2(\R^n)$.
    This quantity is non-negative, real-valued, and symmetric by the corresponding properties of $\hat\rho$.
    Now letting $\rho^{1/2} \defeq \fourier^*[\hat\rho(\freevar)^{1/2}]$, since the Fourier transform is by Plancherel's theorem an isometry of $L^2(\R^n)$, it follows that $\rho^{1/2} \in L^2(\R^n)$.
    Also $\rho^{1/2}$ is real-valued and symmetric by the corresponding properties of $\hat\rho(\freevar)^{1/2}$.
    Finally, if $\hat\rho(\freevar)^{1/2} \in L^1(\R^n)$, we have $\rho^{1/2} \in C_0(\R^n)$ by, e.g., \cite[Theorem 7.5]{rudin2006functional}.
\end{proof}

\begin{lemma}
    \label{lemma:wave:product-fundamental:u}
    Let $\phi, u \in L^2(\R^n) \isect L^1(\R^n)$ satisfy:
    \begin{enumerate}[label=(\roman*)]
        \item
        \label{item:wave:product-fundamental:u}
        $\fourier[u] \ge 0$ and for every non-empty open set $U \subset \R^n$ that $\int_{U} \fourier[u](\xi) \d \xi > 0$.
        \item
        \label{item:wave:product-fundamental:phi}
        $\fourier[\phi] \ge 0$ and for every $r>0$, $\int_{\B(0, r)} \fourier[\phi](\xi) d \xi > 0$.
    \end{enumerate}
    Then $V_{(\phi u)^{1/2}}(\Omega)$ is fundamental for $C(\Omega)$ for any compact $\Omega \subset \R^n$.
\end{lemma}

\begin{proof}
    Let $\rho \defeq \phi u$.
    We use \cite[Theorem 18.1]{cheney2000approximation}, which establishes the required fundamentality provided that $\hat\rho \ge 0$ and  $\int_U \hat\rho \d \xi>0$ for every non-empty open set $U$.
    We have $\hat\rho = \hat\phi * \hat u$, so
    clearly $\hat\rho \ge 0$ by the corresponding assumptions on $\hat\phi$ and $\hat u$.
    Let then $U \subset \R^n$ be open and non-empty.
    There then exists $\eta_0 \in U$ and $r>0$ such that $\B(\eta_0, r) \subset U$.
    Let $U_r \defeq \{\xi \in \R^n \mid \B(\xi, r) \subset U \}$.
    Then $U_r$ is open, non-empty, and $\Isect_{\eta \in \B(0, r)} (U - \eta) \supset U_r$.
    Write $I_{\hat\phi} \defeq \int_{\B(0, r)} \hat\phi(\eta) \d\eta$.
    By assumption \cref{item:wave:product-fundamental:phi}, $I_{\hat\phi}>0$ and $\hat\phi \ge 0$, we have
    \[
        \begin{aligned}
        \int_U \hat\rho(\xi)^{1/2} \d \xi
        &
        =
        \int_U \left( \int_{\R^n} \hat u(\xi - \eta) \hat\phi(\eta) \d\eta\right)^{1/2} \d\xi
        \\
        &
        \ge
        \int_U \left( \int_{\B(0, r)} \hat u(\xi - \eta)  \hat\phi(\eta) \d\eta\right)^{1/2} \d\xi
        \\
        &
        =
        I_{\hat\phi}^{1/2} \int_U \left(\frac{1}{I_{\hat\phi}} \int_{\B(0, r)} \hat u(\xi - \eta)  \hat\phi(\eta)  \d\eta\right)^{1/2} \d\xi.
        \end{aligned}
    \]
    Continuing with the reverse Jensen's inequality for concave functions, Fubini's theorem, and \cref{item:wave:product-fundamental:u}, we get
    \[
        \begin{aligned}
        \int_U \hat\rho(\xi)^{1/2} \d \xi
        &
        \ge
        I_{\hat\phi}^{-1/2} \int_U \int_{\B(0, r)} \hat u(\xi - \eta)^{1/2} \hat\phi(\eta) \d\eta \d\xi
        \\
        &
        \ge
        I_{\hat\phi}^{-1/2}\norm{\hat u}_\infty^{-1/2} \int_U \int_{\B(0, r)} \hat u(\xi - \eta) \hat\phi(\eta)  \d\eta \d\xi
        \\
        &
        \ge
        I_{\hat\phi}^{-1/2}\norm{\hat u}_\infty^{-1/2}  \int_{\B(0, r)} \hat\phi(\eta)\d\eta \int_{U_r} \hat u(\xi) \d\xi 
        \\
        &
        =
        I_{\hat\phi}^{1/2}\norm{\hat u}_\infty^{-1/2} \int_{U_r} \hat u(\xi) \d\xi
        > 0.
        \end{aligned}
    \]
    Since we assume $\Omega$ compact, the aforementioned \cite[Theorem 18.1]{cheney2000approximation} now establishes the claim.
\end{proof}

\begin{example}[Cut Gaussian spread with triangular--Gaussian kernel for $\Wave$]
    \label{ex:convolution:gaussian}
    On $\R$, let $u$ and $v$ be Gaussians of variances $\sigma_u^2$ and $\sigma_v^2$, i.e., $u(x)=C_u e^{-x^2/(2\sigma_u^2)}$ and $v(x)=C_v e^{-x^2/(2\sigma_v^2)}$ for the unit scaling factors $C_u$ and $C_v$.
    Also let $\theta_z(x)=\chi_{[-b,b]}(x-z)$ for some $0 < b < \inf_{z, \tilde z \in \grid} \abs{z-\tilde z}$ on the finite grid $\grid \subset \R$.
    Define $A$ and $\Wave$ according to \eqref{eq:convolution:together:defs} with $\psi=\chi_{[-a,a]}u$ and $\rho=\autoconvolution[\chi_{[-a,a]}]v$ for some $a>0$.
    We know from, e.g., \cite[Appendix 2 \& 3]{kammler2008first} that $\hat u(\xi) = e^{-2\pi^2\sigma_u^2\xi^2}$ and $\hat v(\xi) = e^{-2\pi^2\sigma_v^2\xi^2}$.
    Thus \eqref{eq:convolution:compare:assumption} reads
    \[
        I_{\hat v} e^{-2\pi^2\sigma_v^2\xi^2} \le L_1 e^{-2\pi^2\sigma_u^2\xi^2}
        \quad\text{where}\quad
        I_{\hat v}
        = \sqrt{2\pi}\sigma_v.
    \]
    This holds with
    $L_1 = \sqrt{2\pi}\sigma_v$ if $\sigma_v \ge \sigma_u$.
    By \cref{ex:convolution:rectangular-sensor}, the bound \eqref{eq:convolution:fov-assumption} is satisfied by $L_0 = 2b$.
    Now \cref{lemma:convolution:compare,thm:convolution:together0} establish
    \[
        A_*A \le 2 b \sqrt{2\pi}\sigma_v \Wave
        \quad\text{if}\quad
        \sigma_v \ge \sigma_u.
    \]

    By \cref{lemma:wave:cut-square-root,lemma:wave:product-fundamental:u}, $\rho^{1/2} \in L^2(\R^n)$ exists, and $V_{\rho^{1/2}}$ is fundamental for $C(\Omega)$ for any compact set $\Omega \subset \R^n$.
    However, $\rho^{1/2}$ does not appear to be $C_0(\R^n)$, so \cref{thm:wave:constr:weak-new} does not show that $\norm{\freevar}_\Wave$ would characterise weak-$*$ convergence unlike the “fast” spread of \cref{ex:convolution:fast}.

    The next lemma provides a formula for calculating $\theta_z * \psi$ in the definition of $A$ in \cref{ex:convolution:gaussian}.
    Numerically, the error functions $\erf$ in the lemma are expensive, which leads us to the moniker “fast” for \cref{ex:convolution:fast}. It only involves low-order piecewise polynomials.
\end{example}

\begin{lemma}
    \label{lemma:convolution:fov-spread}
    On $\R$, let $\psi(x) = \chi_{[-a,a]}(x) C e^{-\frac{1}{2\sigma^2}x^2}$ and $\theta_z(x)=\chi_{[-b,b]}(x-z)$ for some $a,b,\sigma>0$.
    Then
    \[
        [\theta_z * \psi](x) =
        \begin{cases}
            0,
            &
            c_1(x) \ge c_2(x),
            \\
            \frac{C\sigma\sqrt{8}}{\sqrt{\pi}}\left[
                \erf\left(\frac{c_2(x)}{\sqrt{2}\sigma}\right)
                -\erf\left(\frac{c_1(x)}{\sqrt{2}\sigma}\right)
            \right],
            &
            \text{otherwise},
        \end{cases}
    \]
    where $c_1(x) \defeq \max\{z-x-b, -a\}$ and $c_2(x) \defeq \min\{z-x+b,a\}$, and the \term{error function}
    \[
        \erf(s) = \frac{2}{\sqrt{\pi}} \int_0^s e^{-t^2}\d t.
    \]
\end{lemma}

\begin{proof}
    We have
    \[
        \begin{aligned}
        [\theta_z * \psi](x)
        &
        =
        \int_{-\infty}^\infty \chi_{[-b,b]}(y-z)\chi_{[-a,a]}(x-y) C e^{-\frac{1}{2\sigma^2}(x-y)^2} \d y
        \\
        &
        =
        \int_{-\infty}^{\infty} \chi_{[z-x-b,z-x+b]}(w)\chi_{[-a,a]}(w) C e^{-\frac{1}{2\sigma^2}w^2} \d w
        \end{aligned}
    \]
    If $c_1(x) \ge c_2(x)$, this gives $[\theta_z * \psi](x)=0$.
    So assume $c_1(x)<c_2(x)$. Then, with the convention that $\int_c^d = - \int_d^c$ if $c>d$, we have
    \[
        \begin{aligned}
        [\theta_z * \psi](x)
        &
        =
        \int_{c_1(x)}^{c_2(x)} C e^{-\frac{1}{2\sigma^2}w^2} \d w
        \\
        &
        =
        \int_{c_1(x)/(\sqrt{2}\sigma)}^{c_2(x)/(\sqrt{2}\sigma)} \sqrt{2}\sigma C e^{-t^2} \d t
        \\
        &
        =
        \int_{0}^{c_2(x)/(\sqrt{2}\sigma)} \sqrt{2}\sigma C e^{-t^2} \d t
        -\int_{0}^{c_1(x)/(\sqrt{2}\sigma)} \sqrt{2}\sigma C e^{-t^2} \d t
        \\
        &
        =
        \frac{C\sigma \sqrt{8}}{\sqrt{\pi}}\left[
            \erf\left(\frac{c_2(x)}{\sqrt{2}\sigma}\right)
            -\erf\left(\frac{c_1(x)}{\sqrt{2}\sigma}\right)
        \right].
        \end{aligned}
    \]
    This establishes the claim.
\end{proof}

\begin{remark}[Higher dimensions]
    In higher dimensions, we will work with uniform products
    \[
        u^{(n)}(x_1,\ldots,x_n) \defeq
        u(x_1) \cdots u(x_n).
    \]
    Thus we will generally replace $\chi = \chi_{[-a,a]}$ by $\chi = \chi_{[-a,a]^n}=\chi_{[-a,a]}^{(n)}$, and $\phi = \autoconvolution[\chi_{[-a,a]}]$ by $\phi = \autoconvolution[\chi_{[-a,a]^n}] = \autoconvolution[\chi_{[-a,a]}]^{(n)}$.
    Since $u^{(n)} * v^{(n)} = (u * v)^{(n)}$ and $\fourier[u^{(n)}] = \fourier[u]^{(n)}$, the above results readily extend to higher dimensions with product factors $L_1^n$ in place of $L_1$ in \eqref{eq:convolution:bochner-condition}.
    The one-dimensional factor $L_1$ can be calculated following \cref{ex:convolution:gaussian,ex:convolution:fast}.
    On a regular grid $\grid$ with $\theta_{(z_1,\ldots,z_n)}(x_1,\ldots,x_n) = \prod_{i=1}^n \theta_{z_i}(x_i)$, also $L_0$ can be replaced by $L_0^n$ for the one-dimensional $L_0$.
\end{remark}

\begin{remark}[Higher dimensions alternative]
    Alternatively, instead of cutting $\psi=\chi v$ with $\chi = \chi_{[-a,a]^n}$, it would be possible to cut with $\chi = \chi_{\B(0, a)}$.
    Then $\phi = \autoconvolution[\chi_{\B(0, a)}]$ in $\rho=\phi u$.
    In $\R^2$, it is possible to calculate using geometric arguments on the area of the intersection of two disks of equal radius that
    \[
        \autoconvolution[\chi_{\B(0, a)}](x)
        =
        \begin{cases}
            0,
            &
            d \ge a, \\
            2a^2*\cos^{-1}(d/a) - d\sqrt{a^2-d^2},
            &
            d < a,
        \end{cases}
    \]
    where $d = \norm{x}/2$ is the distance between the centres of the two disks.
    Thus $\rho = \phi u$ is computable when $u$ is.
    An extension of \cref{lemma:convolution:fov-spread} for the calculation of $\theta_z * \psi$ is more involved.
    The calculation of the factor $L_1$ for the “fast” spread as in \cref{ex:convolution:fast} requires working with Bessel functions in place of the $\sinc$.
\end{remark}

\begin{remark}[Microlocal analysis]
    Fourier analysis of products $\phi u$ is also central to microlocal analysis \cite{hormander2003analysis}.
    However, a family of cut-off functions $\phi$ is usually employed.
\end{remark}

\section{Forward-backward splitting}
\label{sec:fb}

We now develop a forward-backward approach for \eqref{eq:intro:problem}.
In fact, we do so for the general problem
\begin{equation}
    \label{eq:fb:problem}
    \min_{\mu \in \Meas(\Omega)} [F+G](\mu)
\end{equation}
for some $\alpha>0$, convex and pre-differentiable $F$, and
\begin{equation}
    \label{eq:fb:g}
    G(\mu) \defeq \alpha \norm{\mu}_{\Meas} + \delta_{\ge 0}(\mu).
\end{equation}
In the case of \eqref{eq:intro:problem}, $F(\mu) \defeq \frac{1}{2}\norm{A\mu-b}^2$.

In \cref{sec:fb:basic} we formulate the overall method and approximate optimality conditions that each step of our method tries to solve. Then in \cref{sec:fb:insertion} we formulate an algorithm for satisfying these approximate optimality conditions. Finally, in \cref{sec:fb:convergence,sec:fb:weak-convergence}, we prove function value convergence with rates and weak-$*$ convergence of the iterates.

\subsection{Optimality conditions and basic method}
\label{sec:fb:basic}

We choose a self-adjoint and positive semi-definite particle-to-wave operator $\Wave \in \linear(\Masses; \FullPredual)$.
On each step $k \in \N$, given a \term{base point} $\this{\breve\mu} \in \Masses$, we take the next iterate $\nexxt\mu \in \Masses$ as a solution of the (for convenience $\tau$-scaled) surrogate problem
\[
    \min_{\mu \in \Meas(\Omega)} E_k(\mu)
    \quad\text{where}\quad
    E_k(\mu)
    =
    \tau[F(\this{\breve\mu}) + \dualprod{F'(\this{\breve\mu})}{\mu-\this{\breve\mu}} + G(\mu)] + \frac{1}{2}\norm{\mu-\this{\breve\mu}}_{\Wave}^2.
\]
For the methods of this section, the base point $\this{\breve\mu}=\this\mu$, but in \cref{sec:fista} we will use an inertial base point.
The pre-subdifferential of $E_k$ is
\[
    \subdiff E_k(\mu)
    =
    \tau [ F'(\this{\breve\mu}) + \subdiff G(\mu)] + \Wave(\mu-\this{\breve\mu}),
\]
where the pre-subdifferential of $G$ is characterised by
\[
    \alpha \nexxt w \in \subdiff G(\nexxt\mu)
    \iff
    \nexxt w \le 1
    \quad\text{and}\quad
    \nexxt{\mu}(\nexxt w)=\norm{\mu}
    \quad\text{with}\quad
    \nexxt w \in \FullPredual.
\]
Thus the inclusion $\tilde\epsilon_{k+1} \in \subdiff E_k(\nexxt\mu)$ expands with $\this v \defeq F'(\this{\breve\mu})$ as
\begin{align}
    \label{eq:fb:implicit:expanded:convex}
    &
    \left\{\begin{aligned}
        \tilde\epsilon_{k+1} & = \tau[\this v + \alpha \nexxt{w}] + \Wave(\nexxt\mu-\this{\breve\mu}),
        \\
        \nexxt{w} & \le 1 \text{ and } \nexxt{\mu}(\nexxt{w})=\norm{\nexxt\mu}_\Meas.
    \end{aligned}\right.
\end{align}
In practise, we are unable to satisfy the Fermat principle  $0 \in \subdiff E_k(\nexxt\mu)$.
Instead, we pick a tolerance $\epsilon_{k+1}>0$ for $\tilde\epsilon_{k+1}$, and solve for $\nexxt\mu$,
\begin{subequations}
\label{eq:fb:step-problem}
\begin{equation}
    \label{eq:fb:step-problem:main}
    \left\{\begin{aligned}
        \tau[\this v + \alpha] + \Wave(\nexxt\mu - \this{\breve\mu}) & \ge -\epsilon_{k+1}
        &&\text{on } \Omega,
        \\
        \tau[\this v + \alpha] + \Wave(\nexxt\mu - \this{\breve\mu}) & \le \epsilon_{k+1}
        &&\text{on } \supp \nexxt\mu.
    \end{aligned}\right.
\end{equation}
To ensure that $\{\nexxt\mu\}_{k \in \N}$ stays bounded, we further require for some $\kappa > 0$ that
\begin{equation}
    \label{eq:fb:step-problem:boundedness}
    \dualprod{\tau[\this v + \alpha] + \Wave(\nexxt\mu - \this{\breve\mu})}{\nexxt\mu} \le \kappa\epsilon_{k+1}.
\end{equation}
\end{subequations}
Iteratively solving \eqref{eq:fb:step-problem} with $\this{\breve\mu}=\this\mu$ produces the overall structure of \cref{alg:fb:fb}, our proposed forward-backward splitting for \eqref{eq:fb:problem}.

The next lemma shows that for a solution $\nexxt\mu$ to \eqref{eq:fb:step-problem:main} there is some $\epsilon_{k+1}$-bounded $\tilde\epsilon_{k+1} \in \subdiff E_k(\nexxt\mu)$. In the next subsection, we verify the existence of solutions to the entire \eqref{eq:fb:step-problem}.

\begin{lemma}
    \label{lemma:fb:epsilon}
    Suppose $\Omega$ is closed and $\nexxt\mu \in \Meas(\Omega)$ has compact support and solves \eqref{eq:fb:step-problem} for a fixed $k \in \N$ and a given $\this v \in \FullPredual$.
    Then there exists $\tilde\epsilon_{k+1} \in \FullPredual$ satisfying \eqref{eq:fb:implicit:expanded:convex} with $-\epsilon_{k+1} \le \tilde\epsilon_{k+1} \le \epsilon_{k+1}$, and
    \begin{equation}
        \label{eq:fb:tilde-epsilon-estimate}
        \dualprod{\tilde\epsilon_{k+1}}{\nexxt\mu-\opt\mu}
        \le \epsilon_{k+1} (\kappa+\norm{\opt\mu}_{\Masses}).
    \end{equation}
    If $\this v=F'(\this{\breve\mu})$, then $\tilde\epsilon_{k+1} \in \subdiff E_k(\nexxt\mu)$.
\end{lemma}

\begin{proof}
    Pick $\phi \in C_c(\Omega)$ such that
    \[
        1 \ge \phi \ge \frac{1}{\tau\alpha}\max\{0, -\epsilon_{k+1} - \tau\this v - \Wave(\nexxt\mu-\this{\breve\mu})\}
        \quad\text{with}\quad
        \phi|\supp \nexxt\mu=1.
    \]
    Indeed, by the first inequality in \eqref{eq:fb:step-problem:main}, the lower bound on $\phi$ is at most the upper bound. Since $\tau\this v + \Wave(\nexxt\mu-\this{\breve\mu}) \in \FullPredual$ and $\epsilon_{k+1}>0$, the lower bound becomes zero outside a compact set.
    Since $\nexxt\mu$ has compact support, it therefore suffices to take $\phi$ as a mollified indicator of a sufficiently large ball.

    We then let
    
    \begin{gather*}
        \tilde\epsilon_{k+1}
        = \min\{\epsilon_{k+1}, \tau[\this v + \alpha \phi] + \Wave(\nexxt\mu-\this{\breve\mu})\}
        = \tau [\this v + \alpha \nexxt w] + \Wave(\nexxt\mu-\this{\breve\mu})
    \shortintertext{for}
        \nexxt w \defeq \min\left\{\phi, \frac{1}{\tau\alpha}\left(\epsilon_{k+1} - \tau \this v - \Wave(\nexxt\mu-\this{\breve\mu})\right)\right\}.
    \end{gather*}
    Then $\tilde\epsilon_{k+1} \in \FullPredual$, and by construction $-\epsilon_{k+1} \le \tilde\epsilon_{k+1} \le \epsilon_{k+1}$.
    We also have $\nexxt w \le 1$ while the construction of $\phi$ and the second inequality of \eqref{eq:fb:step-problem:main} ensure that $\nexxt w = 1$ on $\supp \nexxt\mu$.
    Therefore, \eqref{eq:fb:implicit:expanded:convex} holds.
    When $\this v=F'(\this{\breve\mu})$, this means that $\tilde\epsilon_{k+1} \in \subdiff E_k(\nexxt\mu)$.

    Since the second inequality of \eqref{eq:fb:step-problem:main} ensures that $\tilde\epsilon_{k+1} = \tau[\this v + \alpha] + \Wave(\nexxt\mu-\this{\breve\mu})$ on $\supp\nexxt\mu$, we now deduce from  \eqref{eq:fb:step-problem:boundedness} that $\dualprod{\tilde\epsilon_{k+1}}{\nexxt\mu} \le \kappa \epsilon_{k+1}$.
    This bound with $\tilde\epsilon_{k+1} \ge -\epsilon_{k+1}$ establishes \eqref{eq:fb:tilde-epsilon-estimate}.
\end{proof}

\begin{algorithm}[t]
    \caption{Forward-backward for Radon norm regularisation of non-negative measures ($\mu$FB)}
    \label{alg:fb:fb}
    \begin{algorithmic}[1]
        \Require Regularisation parameter $\alpha>0$; convex and pre-differentiable $F: \Masses \to \R$;
            self-adjoint particle-to-wave operator $\Wave \in \linear(\Meas(\Omega); \FullPredual)$.
        \State Choose tolerances $\{\epsilon_{k+1}\}_{k \in \N} \subset (0, \infty)$ and a step length $\tau>0$
            subject to \cref{thm:fb:convergence:function-ergodic}, \ref{thm:fb:convergence:function}, or \ref{thm:fb:convergence:weak}.
        \State Choose a fractional tolerance $\kappa \in (0, 1)$ for finite-dimensional subproblems.
        \State Pick an initial iterate $\mu^0 \in \DiscreteMasses$.
        \For{$k \in \N$}
            \State $\this v \defeq F'(\this\mu)$.
            \State\label{step:fb:fb:sub} $\nexxt\mu \defeq \textproc{insert\_and\_adjust}(\this{\mu}, \tau \this v - \Wave\this\mu, \tau\alpha, \epsilon_{k+1}, \kappa)$.
                \Comment{Solves \eqref{eq:fb:step-problem} with \cref{alg:fb:insert-and-remove}.}
            \State Prune zero weight spikes from $\mu^{k+1}$.
                \Comment{Optionally also apply a spike merging heuristic.}
        \EndFor
    \end{algorithmic}
\end{algorithm}

\subsection{Point insertion and reweighting}
\label{sec:fb:insertion}

\begin{algorithm}[t]
    \caption{Point insertion and weight adjustment}
    \label{alg:fb:insert-and-remove}
    \begin{algorithmic}[1]
        \Require $\mu \in \DiscreteMasses$, $\eta \in \FullPredual$, $\lambda,\epsilon > 0$, $\kappa \in (0, 1)$ on a domain $\Omega \subset \R^n$.
        \Function{insert\_and\_adjust}{$\mu$, $\eta$, $\lambda$, $\epsilon$, $\kappa$}
            \State Let $S \defeq \supp\mu$.
            \Repeat
                \State%
                \abovedisplayskip=5pt%
                \belowdisplayskip=5pt%
                Solve for $\vec\beta=(\beta_x)_{x \in S} \in \R^{\#S}$ the finite-dimensional subproblem
                \[
                    \min_{\vec \beta \ge 0}~
                    \frac{1}{2}\iprod{\vec \beta}{D \vec \beta}
                    + \iprod{\vec \eta}{\vec \beta} + \lambda \norm{\vec \beta}_1
                    \quad\text{with}\quad
                    \left\{\begin{array}{l}
                        \vec\eta \defeq (\eta(x))_{x \in S} \in \R^{\#S},
                        \\
                        D \defeq (\rho(y-x))_{x, y \in S} \in \R^{\#S \times \#S},
                    \end{array}\right.
                \]
                to the accuracy
                \belowdisplayskip=0pt%
                \[
                    \norm{ D\vec\beta + \vec\eta + \lambda \vec w}_\infty \le \frac{\kappa\epsilon}{1 + \norm{\vec \beta}_1}
                    \quad\text{for some}\quad \vec w \in \subdiff\norm{\freevar}_1(\vec\beta).
                \]
                \label{step:fb:insert-and-remove:findim}
                \State Let $\mu \defeq \sum_{x \in S} \beta_x \delta_x$
                \label{step:fb:insert-and-remove:mu-update}
                \State Find $\opt x$ (approximately) minimising $\Wave\mu + \eta + \lambda$.
                \label{step:fb:insert-and-remove:optx}
                \Comment{For example, branch-and-bound.}
                \State Let $S \defeq S \union \{\opt x\}$
                \Comment{$\optx$ will only be inserted into $\mu$ if the next bounds check fails.}
            \Until $(\Wave\mu)(\opt x) + \eta(\opt x ) + \lambda \ge -\epsilon$
            \label{step:fb:insert-and-remove:in-bounds}
            \State \Return $\mu$
        \EndFunction
    \end{algorithmic}
\end{algorithm}

We present in \cref{alg:fb:insert-and-remove} a scheme to approximately solve \eqref{eq:fb:step-problem} for $\nexxt\mu$ when starting with $\this{\breve\mu}$ a discrete measure. We denote the class of such measures by $\DiscreteMasses \subset \Masses$.
We stress that \cref{alg:fb:insert-and-remove} is \emph{just one possibility for satisfying \eqref{eq:fb:step-problem}}. In all of our theory it could be replaced with another algorithm that could, for example, incorporate heuristics such as point merging and the sliding of \cite{denoyelle2019sliding}.
The method is based on optimising the weights of the points already in the support of the discrete measure $\mu=\this{\breve\mu}$ and, if necessary, inserting a new point into the support, then repeating until \eqref{eq:fb:step-problem} is satisfied.
Compared to the corresponding step for conditional gradient methods from the literature, which always insert a single point, \cref{alg:fb:insert-and-remove} has one major advantage: it can skip point insertion, if mere weight optimisation is enough to satisfy \eqref{eq:fb:step-problem}. As we have observed in the numerical experiments of \cref{sec:numerical}, this makes the merging heuristics that are critical for conditional gradient methods, unnecessary. Our analysis also incorporates inexact solution of the non-convex point discovery subproblem.

The next lemma proves the finite termination of \cref{alg:fb:insert-and-remove} along with the solvability of \eqref{eq:fb:step-problem}.
It is the only point in our algorithmic theory in this and the following sections, where we need $\Wave$ to be a convolution operator, instead of an abstract operator. Indeed, $\Wave=\rho *$ for a simple single-peaked kernel $\rho$ appears beneficial for the easy solution of \eqref{eq:fb:step-problem}, but in principle other options are possible when the bound $A_*A \le L \Wave$ cannot be satisfied for $\Wave$ a convolution operator.

\begin{lemma}
    \label{lemma:fb:step-problem:alg}
    On a compact set $\Omega \subset \R^n$, let $\Wave \mu = \rho * \mu$ for $\mu \in \Meas(\Omega)$ with the symmetric and positive semi-definite $0 \not\equiv \rho \in C_0(\R^n) \isect L^2(\R^n)$.
    Let $\epsilon_{k+1}>0$ and $\kappa \in (0, 1)$.
    Then, for any $\mu \in \DiscreteMasses$, \eqref{eq:fb:step-problem} can be satisfied in finitely many steps by \textproc{insert\_and\_adjust}($\mu, \tau \this v - \Wave\this{\breve\mu}, \tau\alpha, \epsilon_{k+1}, \kappa$) defined in \cref{alg:fb:insert-and-remove}.
\end{lemma}

\begin{proof}
    Observe that by \cref{thm:wave:constr:weak-new}, $\Wave$ is (strictly) positive definite on $\DiscreteMasses$.
    Let $\eta_k \defeq \Wave\this{\breve\mu} - \tau \this v$.
    For clarity, index the iterations of the inner loop of \cref{alg:fb:insert-and-remove} by $\ell$, and correspondingly write $\mu_\ell=\mu$ and $S_\ell=S$ for the situation on \cref{step:fb:insert-and-remove:mu-update}.
    Since the matrix $D$ on \cref{step:fb:insert-and-remove:findim} is positive definite by the corresponding property of $\Wave$, there always exist solution to the subproblem.
    Moreover, any minimising sequence $\{\vec\beta_j\}_{j \in \N}$ is bounded, ensuring that the $\norm{\vec\beta}_1$-scaled the accuracy requirement can be satisfied by an inexact solution.

    Together with the accuracy requirement, the finite-dimensional subproblem amounts to finding non-negative factors $\{\beta_x\}_{x \in S_\ell}$ such that $\mu_\ell \defeq \sum_{x \in S_\ell} \beta_x \delta_x$ satisfies for some $1 \ge w \in \FullPredual$ with $w(y)=1$ if $\beta_y > 0$ the bounds
    \begin{equation}
        \label{eq:fb:step-problem:findim}
        -\frac{\kappa\epsilon_{k+1}}{1 + \norm{\mu_\ell}_{\Masses}} \le [\Wave\mu_\ell + \eta_k + \lambda w](y) \le \frac{\kappa\epsilon_{k+1}}{1 + \norm{\mu_\ell}_{\Masses}}
        \quad\text{for all}\quad y \in S_\ell.
    \end{equation}
    Since $\kappa \in (0, 1)$, $w=1$ on $\supp\mu_\ell$, and in the parametrisation of $\textproc{insert\_and\_adjust}$, $\lambda=\tau\alpha$ and $\eta_k = \tau \this v - \Wave\this{\breve\mu}$, it follows that $\nexxt\mu = \mu_\ell$ thus constructed satisfies the upper bound of \eqref{eq:fb:step-problem:main} as well as \eqref{eq:fb:step-problem:boundedness}.
    The function $\textproc{insert\_and\_adjust}$ only returns when the lower bound of \eqref{eq:fb:step-problem:main} is also satisfied.
    Therefore, it remains to prove that the lower bound of \eqref{eq:fb:step-problem:main} is satisfied for large enough $\ell$.

    Since $\Wave$ is strictly positive definite, $\rho(0) = \dualprod{\Wave\delta_x}{\delta_x} > 0$ for any $x \in \Omega$. Since $\rho$ is continuous, there exists $R > 0$ such that $\rho \ge \frac{1}{2}\rho(0)\chi_{\B(0, R)}$.
    The upper bound of \eqref{eq:fb:step-problem:findim} ensures for all $\ell \in \N$ and $y \in S_\ell^> \defeq \{y \in S_\ell \mid \beta_y > 0 \}$ that
    \[
        \begin{aligned}
        0
        &
        \le
        \frac{1}{2} \rho(0) \norm{\mu_\ell \llcorner \B(y, R)}_{\Masses}
        =
        \sum_{x \in S_\ell} \frac{1}{2}  \beta_x \rho(0) \chi_{\B(y, R)}(x)
        \\
        &
        \le
        \sum_{x \in S_\ell} \beta_x \rho(y-x)
        =
        \Wave \mu_\ell(y) \le -[\eta_k + \lambda](y)
        + \kappa\epsilon_{k+1}.
        \end{aligned}
    \]
    Since $\eta_k \in \FullPredual$, it is  bounded on the compact set $\Omega$.
    Therefore, by the above, there exists $M > 0$ such that $\norm{\mu_\ell \llcorner \B(y, R)}_{\Masses} \le M$ for all for all $\ell \in \N$ and $y \in S_\ell^>$.
    By the finite Vitali covering lemma, see, e.g., \cite[Proof of Theorem 2.1]{mattila1999geometry}, there exists a subcollection $\tilde S_\ell \subset S_\ell^>$ such that $S_\ell^> \subset \Union_{y \in S_\ell^>} \B(y, R/3) \subset \Union_{y \in \tilde S_\ell} \B(y, R)$, and the latter balls are disjoint.
    Thus
    \[
        \norm{\mu_\ell}_{\Masses}
        = \bigl\| \mu_\ell \llcorner S_\ell^> \bigr\|_{\Masses}
        = \sum_{y \in \tilde S_\ell} \norm{\mu_\ell \llcorner \B(y, R)}_{\Masses}
        \le \#\tilde S_\ell M.
    \]
    Since $\Omega$ is bounded, $\#\tilde S_\ell$ is bounded by a constant dependent on $R$ and $\Omega$ alone; see, e.g., \cite[Lemma 2.6]{mattila1999geometry}.
    Therefore $\{\norm{\mu_\ell}_{\Masses}\}_{\ell \ge 1}$ is bounded.
    \Cref{lemma:distances:uniform-equicontinuity} now shows that the family $\{\Wave\mu_\ell\}_{\ell \in \N}$ is uniformly equicontinuous.
    Since $w(y) \le 1$, the lower bound in \eqref{eq:fb:step-problem:findim} holds with $w$ replaced by $1$.
    Therefore, there must exist $r>0$ such that
    \begin{equation}
        \label{eq:fb:insertion-estimate}
        [\Wave \mu_\ell + \eta_k + \lambda](x) \ge -\frac{\kappa + 1}{2}\epsilon_{k+1} > -\epsilon_{k+1}
        \quad\text{for all}\quad
        x \in \B(y, r),
        \quad
        y \in S_\ell,
        \quad\text{and}\quad
        \ell \ge 1.
    \end{equation}
    It follows that $\B(y, r) \isect S_\ell = \{y\}$ for all $y \in S_\ell$, meaning that \cref{alg:fb:insert-and-remove} does not insert new points into $S$ in the $r$-neighbourhoods of existing points.
    Since a finite number of such neighbourhoods cover the compact set $\Omega$, \cref{alg:fb:insert-and-remove} must terminate in a finite number of steps.
\end{proof}

\begin{remark}[Insertion count]
    \label{rem:fb:insertion}
    The proof of \cref{lemma:fb:step-problem:alg} implies that the maximum number of insertions in \cref{alg:fb:insert-and-remove} is bounded by the number of open $r$-balls with non-intersecting centres, that can be packed into $\Omega$. Subject to further Lipschitz assumptions on $\rho$ and $\range A_*$, such an $r$ is by \eqref{eq:fb:insertion-estimate} proportional to $\epsilon_{k+1}$, so the count is in the order of $1/\epsilon_{k+1}^n$.
    If no pruning is performed in \cref{alg:fb:fb}, this becomes a bound on the cumulative insertions.
    Under a strict complementarity assumption, when $\this\mu$ is near a solution $\opt\mu$, it would be possible to further improve the proof to show that insertions can only happen in a subset of the domain near the spikes of $\opt\mu$.
\end{remark}

\begin{remark}[Complete reconstruction]
    According to \cref{lemma:fb:step-problem:alg}, it would be possible to replace $\this\mu$ by $0$ in the application of \textproc{insert\_and\_adjust} in \cref{alg:fb:fb}.
    This would force complete reconstruction of $\nexxt\mu$ on each iteration.
    In our numerical experience the resulting algorithm has much higher computational demands.
\end{remark}

\begin{remark}[Unconstrained problem]
    \label{rem:fb:unconstrained}
    If $G(\mu)=\alpha\norm{\mu}_{\Masses}$ without the positivity constraint, the inclusion $\tilde\epsilon_{k+1} \in \subdiff E_k(\nexxt\mu)$ expands with $\this v \defeq F'(\this{\breve\mu})$ as
    \begin{align*}
        &
        \left\{\begin{aligned}
            \tilde\epsilon_{k+1} & = \tau[\this v + \alpha \nexxt{w}] + \Wave(\nexxt\mu-\this{\breve\mu}),
            \\
            -1 & \le \nexxt{w} \le 1 \text{ and } \nexxt{\mu}(\nexxt{w})=\norm{\nexxt\mu}_\Meas.
        \end{aligned}\right.
    \end{align*}
    Thus \eqref{eq:fb:step-problem:main} can be replaced by
    \begin{subequations}
    \label{eq:fb:step-problem-unconstr}
    \begin{equation}
        \label{eq:fb:step-problem-unconstr:main}
        \left\{\begin{aligned}
            \tau[\this v + \alpha] + \Wave(\nexxt\mu - \this{\breve\mu}) & \ge -\epsilon_{k+1}
            &&\text{on } \Omega,
            \\
            \tau[\this v - \alpha] + \Wave(\nexxt\mu - \this{\breve\mu}) & \le \epsilon_{k+1}
            &&\text{on } \Omega,
            \\
            \tau[\this v + \alpha] + \Wave(\nexxt\mu - \this{\breve\mu}) & \le \epsilon_{k+1}
            &&\text{on } \supp^+ \nexxt\mu,
            \\
            \tau[\this v - \alpha] + \Wave(\nexxt\mu - \this{\breve\mu}) & \ge -\epsilon_{k+1}
            &&\text{on } \supp^- \nexxt\mu.
        \end{aligned}\right.
    \end{equation}
    where $\supp^\pm \mu \defeq \supp \mu^\pm$ in the minimal decomposition $\mu = \mu^+ - \mu^-$ with $\mu^+,\mu^- \ge 0$.
    To ensure that $\{\nexxt\mu\}_{k \in \N}$ stays bounded, we replace  \eqref{eq:fb:step-problem:boundedness} with
    \begin{equation}
        \label{eq:fb:step-problem-nconstr:boundedness}
        \dualprod{\tau[\this v + \nexxt{\tilde w} \alpha] + \Wave(\nexxt\mu - \this{\breve\mu})}{\nexxt\mu} \le \kappa\epsilon_{k+1},
    \end{equation}
    \end{subequations}
    for some Borel $\nexxt{\tilde w}: \supp\nexxt\mu \to \{1, -1\}$ with $\nexxt{\tilde w} \nexxt\mu = \abs{\nexxt\mu}$.
    For a discrete measure $\nexxt\mu=\sum_{i=1}^n \beta_i \delta_{x_i}$ this falls down to the pointwise values $\nexxt{\tilde w}(x_i)=\sign \beta_i$ for all $i=1,\ldots,n$.
    %

    In \cref{alg:fb:insert-and-remove} all we have to change is:
    \begin{enumerate}[nosep]
        \item On \cref{step:fb:insert-and-remove:findim} remove the constraint $\vec\beta \ge 0$.
        \item On \cref{step:fb:insert-and-remove:optx} find $\opt\xi$ maximising $\abs{\Wave\mu + \eta}$.
        \item On \cref{step:fb:insert-and-remove:in-bounds} exit the loop only if \emph{also} the upper bound check $(\Wave\mu)(\opt x) + \eta(\opt x ) - \lambda \le \epsilon$ succeeds.
    \end{enumerate}
    It is also not difficult to extend \cref{lemma:fb:epsilon,lemma:fb:step-problem:alg} to this setting. Then the theorems in the remainder of this manuscript also go through for the modified algorithm without the positivity constraint.
    %
\end{remark}

\subsection{Function value convergence}
\label{sec:fb:convergence}

We now embark on proving the convergence \cref{alg:fb:fb}, where we could tacitly replace \cref{alg:fb:insert-and-remove} by any other way to solve \cref{eq:fb:step-problem}.
We require some basic regularity and step length conditions.

\begin{assumption}
    \label{ass:fb:all}
    We assume that $\Omega \subset \R^n$ is compact, the regularisation parameter $\alpha>0$, and that
    \begin{enumerate}[label=(\roman*)]
        \item\label{item:fb:all:wave}
        The operator $\Wave \in \linear(\Meas(\Omega); \FullPredual)$ is defined by $\Wave\mu = \rho * \mu$ where $0 \not\equiv \rho \in C_0(\R^n) \isect L^2(\R^n)$ is symmetric and positive definite, i.e., $\rho(-y)=\rho(y)$ for all $y$, and $\fourier[\rho] \ge 0$.
        \item\label{item:fb:all:f}
        The convex function $F: \Masses \to \R$ is pre-differentiable and \emph{$L$-smooth} with respect to $\norm{\freevar}_{\Wave}$:
        \[
            F(\nu)-F(\mu) - \dualprod{F'(\mu)}{\nu-\mu}_{\FullPredual,\Masses} \le \frac{L}{2} \norm{\nu-\mu}_{\Wave}^2
            \quad (\nu, \mu \in \Masses).
        \]
        \item\label{item:fb:all:tau}
        The step length $\tau>0$ satisfies $1 \ge \tau L$ (strictly for weak-$*$ convergence).
    \end{enumerate}
\end{assumption}

\begin{remark}[Compactness]
    The compactness assumption on $\Omega$ is only due to \cref{lemma:fb:step-problem:alg} on the finite termination of \cref{alg:fb:insert-and-remove}. It could be removed if \eqref{eq:fb:step-problem} is otherwise ensured on \cref{step:fb:fb:sub} of the algorithm.
\end{remark}

\begin{example}
    \label{example:fb:convergence:standard-f-lipschitz}
    For $F(\mu)=\frac{1}{2}\norm{A\mu-b}_Y^2$ with $A \in \linear(\Masses; Y)$ and $Y$ a Hilbert space, \cref{ass:fb:all}\,\cref{item:fb:all:f} reduces to
    $A$ being pre-adjointable with $A_* A \le L \Wave$.
\end{example}

Our convergence proofs are based on the following simplified (scalar-tested, without strong convexity) version of \cite[Theorem 2.3]{tuomov-inertia} with the addition of inexact solutions ($\epsilon_{k+1} \ne 0$).
It includes inertial parameters $\lambda_k$ that we need in \cref{sec:fista}.
For non-inertial methods we take $\lambda_k \equiv 1$ and $\thisz \equiv \this{\breve x} \equiv \thisx$.

\begin{lemma}
    \label{lemma:fb:general-estimate}
    On a Banach space $X$ with predual $X_*$, let $F, G: X \to \extR$ with $F$ transport-subdifferentiable and $F$ 3-point smooth.
    For an initial $x^0 = \breve x^0 \in X$, let $\{\nextx\}_{k \in \N}$ and $\{\nexxt{\breve x}\}_{k \in \N}$ be defined for some $\lambda_k>0$ and $\tilde\epsilon_{k+1} \in X_*$ through
    \begin{gather}
        \nonumber
        \tilde\epsilon_{k+1} \in \tau_k[F'(\this{\breve x}) + \subdiff G(\nextx)] + \Wave(\nextx-\this{\breve x})
    \shortintertext{and}
        \label{eq:fb:inertial-relationship}
        \lambda_k(\nextz-\thisz) = \nextx-\this{\breve x}
        \quad\text{as well as}\quad
        \lambda_k(\nexxt z-\thisx) = \nextx - \thisx.
    \end{gather}
    Then for any $\optx \in X$ we have
    \begin{multline*}
        \frac{\lambda_k^2}{2}\norm{\nextz-\optx}_{\Wave}^2
        +\frac{\lambda_k^2(1-\tau_k L)}{2}\norm{\nextz-\thisz}_{\Wave}^2
        + \tau_k([F+G](\nextx) - [F+G](\optx))
        \\
        - (1-\lambda_k)\tau_k([F+G](\thisx)-[F+G](\optx))
        \le
        \frac{\lambda_k^2}{2}\norm{\thisz-\optx}_{\Wave}^2
        +\lambda_k\dualprod{\tilde\epsilon_{k+1}}{\thisz - \opt x}.
    \end{multline*}
    In particular, if $\lambda_k \equiv 1$ and $\thisz \equiv \this{\breve x} \equiv \thisx$, we have
    \begin{equation*}
        \frac{1}{2}\norm{\nextx-\optx}_{\Wave}^2
        +\frac{1-\tau_k L}{2}\norm{\nextx-\thisx}_{\Wave}^2
        + \tau_k([F+G](\nextx) - [F+G](\optx))
        \le
        \frac{1}{2}\norm{\thisx-\optx}_{\Wave}^2
        +\dualprod{\tilde\epsilon_{k+1}}{\thisx - \opt x}.
    \end{equation*}
\end{lemma}

\begin{proof}
    Let $\nexxt q \in \subdiff G(\nextx)$ satisfy $\tilde\epsilon_{k+1} = \tau_k[F'(\this{\breve x}) + \nexxt q] + \Wave(\nextx-\this{\breve x})$.
    Then
    \begin{equation}
        \label{eq:fb:tested}
        \lambda_k \dualprod{\tilde\epsilon_{k+1}}{\nexxt z - \opt x}
        =
        \lambda_k \dualprod{\tau_k[F'(\this{\breve x}) + \nexxt q] + \Wave(\nextx-\this{\breve x})}{\nextz - \opt x}.
    \end{equation}
    Due to the three-point identity \eqref{eq:distances:wave-3-point} and \eqref{eq:fb:inertial-relationship}, we have
    \[
        \begin{aligned}[t]
        \lambda_k\dualprod{\Wave(\nextx-\this{\breve x})}{\nexxt z - \opt x}
        &
        =
        \lambda_k^2\dualprod{\Wave(\nextz-\thisz)}{\nexxt z - \opt x}
        \\
        &
        =\frac{\lambda_k^2}{2}\norm{\nextz-\thisz}_{\Wave}^2
        +\frac{\lambda_k^2}{2}\norm{\nextz-\optx}_{\Wave}^2
        -\frac{\lambda_k^2}{2}\norm{\nextz-\optx}_{\Wave}^2.
        \end{aligned}
    \]
    By the definition of the convex pre-subdifferential and the $L$-smoothness of $F$ with respect to $\Wave$, using both relationships of \eqref{eq:fb:inertial-relationship}, we obtain
    \[
        \begin{aligned}[t]
        \lambda_k \dualprod{F'(\this{\breve x}) + \nexxt q}{\nextz - \opt x}
        &
        =
        \lambda_k \dualprod{F'(\this{\breve x}) + \nexxt q}{\nextx - \opt x}
        + (1-\lambda_k) \dualprod{F'(\this{\breve x}) + \nexxt q}{\nextx - \thisx}
        \\
        &
        \ge
        \lambda_k\left(
            [F+G](\nextx)-[F+G](\optx) - \frac{L}{2}\norm{\nextx-\this{\breve x}}_\Wave^2
        \right)
        \\
        \MoveEqLeft[-1]
        +(1-\lambda_k)\left(
            [F+G](\nextx)-[F+G](\thisx) - \frac{L}{2}\norm{\nextx-\this{\breve x}}_\Wave^2
        \right)
        \\
        &
        =
        ([F+G](\nextx)-[F+G](\optx))
        \\
        \MoveEqLeft[-1]
        - (1-\lambda_k)([F+G](\thisx)-[F+G](\optx))
        - \frac{L\lambda_k^2}{2}\norm{\nextz-\thisz}_\Wave^2.
        \end{aligned}
    \]
    Inserting these two expressions into \eqref{eq:fb:tested}, the claim follows.
\end{proof}

We now readily obtain ergodic function value convergence:

\begin{theorem}[Ergodic function value convergence]
    \label{thm:fb:convergence:function-ergodic}
    Suppose \cref{ass:fb:all} holds and let $\opt\mu \in \Masses$ satisfy $0 \in \subdiff[F + G](\opt \mu)$.
    Let $\{\this\mu\}_{k \ge 1}$ be generated by \cref{alg:fb:fb} for some $\mu^0 \in \DiscreteMasses$ with the tolerance sequence $\{\epsilon_{k+1}\}_{k \in \N} \subset (0, \infty)$ satisfying
    \begin{equation}
        \label{eq:fb:tilde-cn}
        \lim_{N \to \infty}  \frac{1}{N} \sum_{k=0}^{N-1} \epsilon_{k+1} = 0.
    \end{equation}
    Define the \term{ergodic iterate} $\tilde\mu^N \defeq \frac{1}{N} \sum_{k=0}^{N-1} \this\mu$.
    Then $[F+G](\tilde\mu^N) \to [F+G](\opt\mu)$, more precisely,
    \[
        [F+G](\tilde\mu^N)
        \le
        [F+G](\opt\mu)
        + \frac{1}{N\tau}\left(
            (\kappa+\norm{\opt\mu}_{\Masses}) \sum_{k=0}^{N-1} \epsilon_{k+1}
            + \frac{1}{2}\norm{\mu^0-\opt\mu}_{\Wave}^2
        \right).
    \]
\end{theorem}

\begin{proof}
    \Cref{alg:fb:fb,lemma:fb:step-problem:alg,ass:fb:all} ensure \eqref{eq:fb:step-problem} for all $k \in \N$.
    Therefore \cref{lemma:fb:epsilon} provides $\tilde\epsilon_{k+1} \in \subdiff E_k(\nexxt\mu)$ satisfying \eqref{eq:fb:tilde-epsilon-estimate}.
    The case $\lambda_k \equiv 1$ of \cref{lemma:fb:general-estimate} with $\tau_k \equiv \tau$ now establishes
    \begin{equation}
        \label{eq:convergence:descent-estimate}
        \frac{1}{2}\norm{\nexxt\mu - \opt\mu}_\Wave^2
        + \frac{1-\tau L}{2}\norm{\nexxt\mu-\this\mu}_\Wave^2
        + \tau [F+G](\nexxt\mu) - \tau [F+G](\opt\mu)
        \le
        \frac{1}{2}\norm{\this\mu - \opt\mu}_\Wave^2
        + \dualprod{\tilde\epsilon_{k+1}}{\nexxt\mu-\opt\mu}.
    \end{equation}
    Summing over $k=0,\ldots,N-1$ we obtain
    \begin{equation}
        \label{eq:convergence:ergodic-first-estimate}
        \frac{1}{2}\norm{\mu^N-\opt\mu}_{\Wave}^2
        + \sum_{k=0}^{N-1}\left(
            \tau[F+G](\nexxt\mu) - \tau[F+G](\opt\mu)
            + \frac{1-\tau L}{2}\norm{\nexxt\mu-\this\mu}_\Wave^2
            \right)
        \le
        \frac{1}{2}\norm{\mu^0-\opt\mu}_{\Wave}^2
        + \tilde C_N,
    \end{equation}
    where \eqref{eq:fb:tilde-epsilon-estimate} shows that
    $
        \tilde C_N \defeq \sum_{k=0}^{N-1} \dualprod{\tilde\epsilon_{k+1}}{\nexxt\mu-\opt\mu}
         \le (\kappa+\norm{\opt\mu}_{\Masses}) \sum_{k=0}^{N-1} \epsilon_{k+1}.
    $
    Since $\tau L \le 1$, dividing \eqref{eq:convergence:ergodic-first-estimate} by $N\tau$ and using Jensen's inequality, the claim follows.
\end{proof}

\begin{example}[Tolerance sequence]
    \label{ex:fb:tolerance:ergodic}
    Take $\epsilon_{k+1}=1/(k+1)^p$ for some $p > 1$.
    Then \eqref{eq:fb:tilde-cn} holds.
    In fact, $\sum_{k=0}^{N-1} \epsilon_{k+1} < \infty$, so \cref{thm:fb:convergence:function-ergodic} yields $O(1/N)$ ergodic function value convergence.
\end{example}

\begin{example}[Controlling the insertion count]
    Recall \cref{rem:fb:insertion}. To keep the growth of the (cumulative) maximum insertion count in \cref{alg:fb:insert-and-remove} linear, we could take $\epsilon_{k+1}=1/(k+1)^{1/n}$. This comes at the expense of reducing the $O(\frac{1}{N}\sum_{k=0}^{N-1}\epsilon_{k+1})$ convergence rate given by the theorem to $O(\log N/N)$ for $n=1$ and $O(N^{-1/n})$ for $n>1$.
    By our numerical experience in \cref{sec:numerical}, in practise even the option of \cref{ex:fb:tolerance:ergodic} requires very few insertions on each iteration.
\end{example}

With a tighter finite-dimensional subproblem solution accuracy, the method is nearly monotone, and we get non-ergodic convergence. If the tolerances are as in  \cref{ex:fb:tolerance:ergodic}, this convergence is $O(1/N)$.

\begin{theorem}[Function value convergence]
    \label{thm:fb:convergence:function}
    Suppose \cref{ass:fb:all} holds and let $\opt\mu \in \Masses$ satisfy $0 \in \subdiff[F + G](\opt \mu)$.
    Let $\{\this\mu\}_{k \ge 1}$ be generated by \cref{alg:fb:fb} for some $\mu^0 \in \DiscreteMasses$ with the tolerance sequence $\{\epsilon_{k+1}\}_{k \in \N} \subset (0, \infty)$ satisfying \eqref{eq:fb:tilde-cn}, and \eqref{eq:fb:step-problem} amended by
    \begin{equation}
        \label{eq:fb:step-problem:boundedness:extra}
        \dualprod{\min\{\epsilon_{k+1}, \tau[\this v + \alpha] + \Wave(\nexxt\mu - \this{\mu}) \}}{\nexxt\mu-\this\mu} \le \frac{\kappa\epsilon_{k+1}}{k}
        \quad\text{for}\quad k \ge 1.
    \end{equation}
    Then $[F+G](\mu^N) \to [F+G](\opt\mu)$, more precisely,
    \[
        [F+G](\mu^N)
        \le
        [F+G](\opt\mu)
        + \frac{1}{N\tau}\left( (2\kappa+\norm{\opt\mu}_{\Masses}) \sum_{k=0}^{N-1} \epsilon_{k+1} + \frac{1}{2}\norm{\mu^0-\opt\mu}_{\Wave}^2\right).
    \]
\end{theorem}

\begin{proof}
    \Cref{alg:fb:fb,lemma:fb:step-problem:alg,ass:fb:all} ensure \eqref{eq:fb:step-problem} for all $k \in \N$.
    Therefore \cref{lemma:fb:epsilon} provides $\tilde\epsilon_{k+1} \in \subdiff E_k(\nexxt\mu)$ satisfying \eqref{eq:fb:tilde-epsilon-estimate}. In fact, $\tilde\epsilon_{k+1}= \min\{\epsilon_{k+1}, \tau[\this v + \alpha] + \Wave(\nexxt\mu - \this{\mu}) \}$ in the proof. Therefore, \eqref{eq:fb:step-problem:boundedness:extra} reads
    $
        \dualprod{\tilde\epsilon_{k+1}}{\nexxt\mu-\this\mu} \le \frac{\kappa\epsilon_{k+1}}{k}.
    $
    Since $\tau L \le 1$, the case $\lambda_j \equiv 1$ of \cref{lemma:fb:general-estimate} with $\tau_k \equiv \tau$ and $\optx=\mu^k$ now establishes for all $j \ge 1$ the “monotonicity with error”
    \begin{equation*}
        \tau [F+G](\mu^{j+1})
        \le
        \tau [F+G](\mu^j)
        + \dualprod{\tilde\epsilon_{j+1}}{\mu^{j+1}-\mu^j}
        \le
        \tau [F+G](\mu^j)
        + \frac{\kappa\epsilon_{j+1}}{j}.
    \end{equation*}
    Hence, summing over all $j=k+1,\ldots,N-1$, we obtain for all $N \ge k$ and $k \in \N$ that
    \begin{equation*}
        \tau [F+G](\mu^N)
        \le
        \tau [F+G](\nexxt\mu)
        + \kappa \sum_{j=k+1}^{N-1} \frac{\epsilon_{j+1}}{j}.
    \end{equation*}
    Applying this estimate in \eqref{eq:convergence:ergodic-first-estimate} for all $k=0,\ldots,N-1$, yields
    \[
        \frac{1}{2}\norm{\mu^N - \opt\mu}_\Wave^2
        + N \tau \bigl( [F+G](\mu^N) - [F+G](\opt\mu) \bigr)
        \le
        \frac{1}{2}\norm{\mu^0 - \opt\mu}_\Wave^2
        + \kappa \sum_{k=0}^{N-1} \sum_{j=k+1}^{N-1} \frac{\epsilon_{j+1}}{j}
        + \sum_{k=0}^{N-1}\dualprod{\tilde\epsilon_{k+1}}{\nexxt\mu-\opt\mu}.
    \]
    Since $\sum_{k=0}^{N-1} \sum_{j=k+1}^{N-1} \frac{\epsilon_{j+1}}{j} = \sum_{k=1}^{N-1} \epsilon_{k+1}$, using \eqref{eq:fb:tilde-epsilon-estimate} on the final term and dividing by $N\tau$, establishes the claim.
\end{proof}

\begin{remark}[Finite-dimensional subproblems]
    The accuracy requirement \eqref{eq:fb:step-problem:boundedness:extra} holds if
    \[
        \abs{\tau[\this v + \alpha \this w] + \Wave(\nexxt\mu - \this{\mu})}(y) \le \frac{\kappa\epsilon_{k+1}}{k(\norm{\nexxt\mu}_\Meas+\norm{\this\mu}_\Meas)}
        \quad\text{for}\quad
        y \in \supp \nexxt\mu \union \supp\this\mu,
    \]
    where $1 \ge \this w \in \FullPredual$ with $\this w = 1$ on $\supp\nexxt\mu$.
    Since the set $S$ only grows in \cref{alg:fb:insert-and-remove}, it follows that \eqref{eq:fb:step-problem:boundedness:extra} is satisfied by tightening the right-hand-side of the accuracy requirement on \cref{step:fb:insert-and-remove:findim} to $\kappa\epsilon_{k+1}/(1 + k(\norm{\vec\beta}_1 + \norm{\vec\beta_0}_1))$, where $\vec\beta_0$ is the vector of weights of the initial $\mu$.
\end{remark}

\begin{remark}[Linear convergence]
    If $F$ is strongly convex with respect to the $\Wave$-seminorm, it is not difficult to prove linear convergence along the lines of the Hilbert space proofs in \cite{clasonvalkonen2020nonsmooth,tuomov-proxtest}.
    For $F(\mu)=\frac{1}{2}\norm{A\mu-b}^2$ the strong convexity requirement reduces to $A_*A \ge \gamma \Wave$ for some $\gamma>0$.
    In practise, we do not expect this condition to hold, so do not pursue the precise proof here.
    It would be possible to relax the strong convexity requirement to the (strong or non-strong) \term{metric subregularity} of $\subdiff[F+G]$ near a solution;
    see \cite{tuomov-subreg,clasonvalkonen2020nonsmooth} for the convergence theory and \cite{tuomov-regtheory} for examples of satisfaction of such a condition in other contexts.
    Proving metric subregularity for the problem \eqref{eq:intro:problem} is outside the scope of the present manuscript.
\end{remark}

\subsection{Weak-$*$ convergence}
\label{sec:fb:weak-convergence}

We now prove weak-$*$ convergence of the iterates.
We recall the following deterministic version of the results of \cite{robbins1971convergence}:

\begin{lemma}
    \label{lemma:convergence-inequality}
    Let $\{a_k\}_{k \in \N}$, $\{b_k\}_{k \in \N}$, $\{c_k\}_{k \in \N}$, and $\{d_k\}_{k \in \N}$ be non-negative sequences with
    \[
        a_{k+1} \le a_k(1+b_k) + c_k - d_k
        \quad\text{for all}\quad k \in \N.
    \]
    If $\sum_{k=0}^\infty b_k < \infty$ and $\sum_{k=0}^\infty c_k < \infty$, then
    \begin{enumerate*}[label=(\roman*)]
        \item $\lim_{k \to \infty} a_k$ exists and is finite; and
        \item $\sum_{k=0}^\infty d_k < \infty$.
    \end{enumerate*}
\end{lemma}

The next generalisation of Opial's lemma \cite{opial1967weak} follows the proof from \cite{clasonvalkonen2020nonsmooth} in Hilbert spaces.
We prove it for Bregman divergences
\[
    B_M(x, z) \defeq M(z) - M(x) - \dualprod{M'(x)}{z-x},
\]
as they add no extra difficulties. Recall that in Hilbert spaces $\frac{1}{2}\norm{z-x}^2=B_M(x, z)$ for $M=\frac{1}{2}\norm{\freevar}^2$.

\begin{lemma}
    \label{lemma:opial}
    Let $X=(X_*)^*$ be the dual space of a corresponding separable normed space $X_*$.
    Also let $M : X \to \R$ be convex, proper, and Gateaux differentiable with $M': X \to X_*$ weak-$*$-to-weak continuous.
    Finally, let $\hat X \subset X$ be non-empty and $\{e_{k+1}\}_{k \in \N} \in X_*$. If
    \begin{enumerate}[label=(\roman*)]
        \item\label{item:opial-limit} all weak-$*$ limit points of $\{\thisx\}_{k \in \N}$ belong $\hat X$.
        \item\label{item:opial-epsilon-bound} $\sum_{k=0}^\infty \max\{0, \dualprod{e_{k+1}}{\nextx-\bar x}\} <\infty$ for all $\bar x \in \hat X$; and
        \item\label{item:opial-nonincreasing} $B_M(\nextx, \optx) \le B_M(\thisx, \optx) + \dualprod{e_{k+1}}{\nextx-\optx}$ for all $\optx \in \hat X$ and $k \in \N$;
    \end{enumerate}
    then all weak-$*$ limit points of $\{\this x\}_{k \in \N}$ satisfy $\hat x, \opt x \in \hat X$ and
    \begin{equation}
        \label{eq:opial:mprime-zero}
        \dualprod{M'(\hat x) - M'(\bar x)}{\hat x - \bar x} = 0.
    \end{equation}
    If $\{\thisx\}_{k \in \N} \subset X$ is bounded, then such a limit point exists.
    If, in addition to all the previous assumptions, \eqref{eq:opial:mprime-zero} implies $\hat x = \bar x$ (such as when $M$ is strongly monotone), then $\thisx \weaktostar \hat x$ in $X$ for some $\hat x \in \hat X$.
\end{lemma}
\begin{proof}
    Let $\bar x$ and $\hat x$ be weak-$*$ limit points of $\{\thisx\}_{k \in \N}$.
    Since Bregman divergences $B_M \ge 0$ for convex $M$, the conditions \cref{item:opial-epsilon-bound,item:opial-nonincreasing} establish the assumptions of \cref{lemma:convergence-inequality} for $a_k = B_M(\thisx; \optx)$, $b_k=0$, $c_k=\max\{0, \dualprod{e_{k+1}}{\nextx-\bar x}\}$, and $d_k=0$. It follows that $\{B_M(x^k; \bar x)\}_{k\in\N}$ is convergent.
    Likewise we establish that $\{B_M(x^k; \hat x)\}_{k\in\N}$ is convergent.
    Therefore
    \begin{equation*}
         \dualprod{M'(x^k)}{\bar x - \hat x}
            = B_M(x^k; \hat x) - B_M(x^k; \bar x) + M(\bar x) - M(\hat x)
            \to c \in \R.
    \end{equation*}
    Since $\bar x$ and $\hat x$ are a weak-$*$ limit point, there exist subsequence $\{x^{k_n}\}_{n\in \N}$ and $\{x^{k_m}\}_{m\in \N}$ with $x^{k_n}\weakto \bar x$ and $x^{k_m}\weakto \hat x$.
    Using the weak-$*$-to-weak continuity of $M': X \to X_*$, \eqref{eq:opial:mprime-zero} follows from
    \begin{equation*}
        \dualprod{M'(\bar x)}{\bar x - \hat x}
        = \lim_{n\to \infty} \iprod{M'(x^{k_n})}{\bar x - \hat x}
        = c
        = \lim_{n\to \infty} \iprod{M'(x^{k_m})}{\bar x - \hat x}
        = \dualprod{M'(\bar x)}{\bar x - \hat x}.
    \end{equation*}

    If $\{\thisx\}_{k \in \N}$ is bounded, it contains a \mbox{weakly-$*$} convergent subsequence by the Banach--Alaoglu theorem, so a limit point exists as claimed.

    Hence, if \eqref{eq:opial:mprime-zero} implies $\bar x =\hat x$, then every convergent subsequence of  $\{\thisx\}_{k \in \N}$ has the same weak limit.
    It lies in $\hat X$ by \cref{item:opial-limit}.
    The final claim now follows from a standard subsequence--subsequence argument: Assume to the contrary that there exists a subsequence of $\{x^k\}_{k\in\N}$ not convergent to $\hat x$. Then we can apply the above argument to obtain a further subsequence converging to $\hat x$.
    This contradicts the fact that any subsequence of a convergent sequences converges to the same limit.
\end{proof}

Finally, we can state and prove the weak convergence result.

\begin{theorem}[Weak convergence]
    \label{thm:fb:convergence:weak}
    Suppose \cref{ass:fb:all} holds with the step length condition \cref{item:fb:all:tau} strictly: $\tau L < 1$.
    Let $\{\this\mu\}_{k \ge 1}$ be generated by \cref{alg:fb:fb} for some $\mu^0 \in \DiscreteMasses$.
    Also suppose
    \begin{enumerate}[label=(\roman*)]
        \item\label{item:fb:convergence:weak:lowerbound}
        $F$ is bounded from below;
        \item\label{item:fb:convergence:weak:fprime-cont}
        $F'$ is continuous with respect to $\norm{\freevar}_\Wave$; and
        \item\label{item:fb:convergence:weak:cn}
        the tolerance sequence $\{\epsilon_{k+1}\}_{k \in \N} \subset (0, \infty)$ satisfies $\sum_{k=0}^{\infty} \epsilon_{k+1} < \infty$.
    \end{enumerate}
    Then there exists a weak-$*$ limit point $\bar\mu$ of $\{\this\mu\}_{k \in \N}$ satisfying $0 \in \subdiff[F + G](\bar\mu)$, and any two such points $\bar\mu, \hat\mu$ satisfy $\norm{\bar\mu-\hat\mu}_\Wave=0$.
    If a point satisfying both properties is unique, in particular when all the assumptions of \cref{thm:wave:constr:weak-new}\,\cref{item:wave:constr:wave-to-weak-new} hold for $\rho$, then $\this\mu \weaktostar \opt\mu$ weakly-$*$ in $\Meas(\Omega)$.
\end{theorem}

We recall that the assumptions of \cref{thm:wave:constr:weak-new}\,\cref{item:wave:constr:wave-to-weak-new} hold for the kernel $\rho$ of \cref{ex:convolution:fast}, but the $C_0$ assumption on $\rho^{1/2}$ presumably fails for \cref{ex:convolution:gaussian}.

\begin{proof}
    We use \cref{lemma:opial} with $\hat X = \inv{(\subdiff[F+G])}(0)$ and $M=\frac{1}{2}\norm{\freevar}_\Wave^2$, in which case $B_M(\mu,\nu)=\frac{1}{2}\norm{\mu-\nu}_\Wave^2$.
    We need to verify its assumptions.
    \Cref{alg:fb:fb,lemma:fb:step-problem:alg,ass:fb:all} ensure \eqref{eq:fb:step-problem} for all $k \in \N$.
    Therefore \cref{lemma:fb:epsilon} provides $\tilde\epsilon_{k+1} \in \subdiff E_k(\nexxt\mu)$ satisfying $-\epsilon_{k+1} \le \tilde\epsilon_{k+1} \le \epsilon_{k+1}$ as well as \eqref{eq:fb:tilde-epsilon-estimate}.
    The latter with the present assumption \cref{item:fb:convergence:weak:cn} establishes
    \begin{equation}
        \label{eq:fb:convergence:weak:c-def}
        \sum_{k=0}^{\infty} \max\{0, \dualprod{\tilde\epsilon_{k+1}}{\nexxt\mu-\opt\mu}\}
        \le
        (\kappa + \norm{\opt\mu}_{\Masses})\sum_{k=0}^{\infty} \epsilon_{k+1}
        =: C
        \quad\text{for all} \quad \opt\mu \in \hat X.
    \end{equation}
    Hence \cref{lemma:opial}\,\cref{item:opial-epsilon-bound} holds.
    The case $\lambda_k \equiv 1$ of \cref{lemma:fb:general-estimate} with $\tau_k \equiv \tau$ establishes the Féjer monotonicity “with error”
    \[
        \frac{1}{2}\norm{\nexxt\mu - \opt\mu}_\Wave^2
        \le
        \frac{1}{2}\norm{\this{\mu} - \opt\mu}_\Wave^2
        + \dualprod{\tilde\epsilon_{k+1}}{\nexxt\mu-\opt\mu}
    \]
    for every $\opt\mu \in \hat X$.
    This proves the assumption \cref{lemma:opial}\,\cref{item:opial-nonincreasing}.
    We still need to prove \cref{lemma:opial}\,\cref{item:opial-limit}.
    We do this together with showing the boundedness of $\{\this\mu\}_{k \in \N}$.

    The conditions of \cref{thm:fb:convergence:function-ergodic} readily hold for any $\opt\mu \in \hat X$ by our assumptions. Thus \eqref{eq:convergence:descent-estimate} and \eqref{eq:convergence:ergodic-first-estimate} in its proof hold for all $N \in \N$.
    Since $0 \in \subdiff[F+G](\opt\mu)$, we have $[F+G](\nexxt\mu) \ge [F+G](\opt\mu)$ for all $k \in \N$.
    Therefore, since $\tau L < 1$, \eqref{eq:convergence:ergodic-first-estimate} and \ref{item:fb:convergence:weak:cn} allow us to deduce both $\norm{\nexxt\mu-\this\mu}_\Wave \to 0$ and
    \begin{equation}
        \begin{aligned}[t]
        \frac{1}{2}\norm{\mu^0 - \opt\mu}^2
        + C
        &
        \ge
        \sum_{k=0}^{N-1}\left(\tau[F+G](\nexxt\mu) - \tau[F+G](\opt\mu)
            +\frac{1}{2}\norm{\nexxt\mu - \this\mu}^2
        \right)
        \\
        &
        \ge
        \tau[F+G](\mu^{j+1}) - \tau[F+G](\opt\mu)
        \\
        &
        \ge
        \tau\alpha\norm{\mu^{j+1}}_{\Meas}  + \tau \inf F - \tau[F+G](\opt\mu)
        \end{aligned}
    \end{equation}
    for $C$ defined in \eqref{eq:fb:convergence:weak:c-def}, for any $j \in \{0, \ldots, N-1\}$ and $\opt\mu \in \hat X$.
    By \ref{item:fb:convergence:weak:lowerbound}, this establishes the boundedness of $\{\this\mu\}_{k \in \N}$.

    For any $x \in \Omega$ and $k \in \N$, we have
    \begin{equation}
        \label{eq:fb:weak:x-test}
        \abs{\dualprod{\delta_x}{\Wave(\nexxt\mu-\this\mu)}}
        = \abs{\iprod{\delta_x}{\nexxt\mu-\this\mu}_\Wave}
        \le \norm{\delta_x}_{\Wave}\norm{\nexxt\mu-\this\mu}_\Wave
        = \abs{\rho(0)}\norm{\nexxt\mu-\this\mu}_\Wave.
    \end{equation}
    Since we have just shown that $\norm{\nexxt\mu-\this\mu}_\Wave \to 0$,
    it follows from \eqref{eq:fb:weak:x-test} that $\Wave(\nexxt\mu-\this\mu) \to 0$ uniformly, i.e., strongly in $\FullPredual$.
    By \cref{item:fb:convergence:weak:fprime-cont}, moreover, $F'(\this\mu)-F'(\nexxt\mu) \to 0$ strongly in $\FullPredual$.
    Let then $\eta_k \defeq F'(\nexxt\mu) +  \alpha \nexxt w \in \subdiff[F+G](\nexxt\mu)$.
    Since we have assumed $\epsilon_{k+1} \downto 0$, also $\tilde\epsilon_{k+1} \to 0$ strongly.
    Testing the approximate optimality conditions \eqref{eq:fb:implicit:expanded:convex} by $\delta_x$ we now have
    \[
       \tau\eta_k(x)
       = \tilde \epsilon_{k+1}(x) - \tau[F'(\this\mu)-F'(\nexxt\mu)](x) - \dualprod{\delta_x}{\Wave(\nexxt\mu-\this\mu)}
       \to 0
       \quad\text{uniformly in $\Omega$}.
    \]
    Since $\eta_k \in \subdiff[F+G](\nexxt\mu)$, the weak-$*$-to-strong outer semicontinuity of (pre)subdifferentials (see \cite{clasonvalkonen2020nonsmooth}) establishes that any weak-$*$ limit point $\hat\mu$ of $\{\this\mu\}_{k \in \N}$ satisfies $0 \in \subdiff[F+G](\hat\mu)$.
    This establishes the assumption \cref{lemma:opial}\,\cref{item:opial-limit}.

    We have now verified the main conditions of \cref{lemma:opial}, which shows that all weak-$*$ limit points $\hat\mu,\bar\mu$ of $\{\this\mu\}_{k \in \N}$ belong to $\inv{(\subdiff[F+G])}(0)$ and satisfy
    $
         0
         =
         \dualprod{M'(\hat\mu) - M'(\bar\mu)}{\hat\mu-\bar\mu}
         =
         \norm{\hat\mu-\bar\mu}_\Wave^2.
    $
    Since we have shown that $\{\this\mu\}_{k \in \N}$ is bounded, such a limit point exists.
    This establishes the first claim.

    If now $\norm{\hat\mu-\bar\mu}_\Wave=0$ implies $\hat\mu=\bar\nu$, then \cref{lemma:opial} establishes weak-$*$ convergence.
    When the assumptions of \cref{thm:wave:constr:weak-new}\,\cref{item:wave:constr:wave-to-weak-new} hold, \cref{cor:distances:norm} proves that $\norm{\freevar}_\Wave$ is a norm, in particular, so this property holds.
\end{proof}

\begin{remark}
    The condition \cref{thm:fb:convergence:weak}\,\cref{item:fb:convergence:weak:fprime-cont} along with \cref{ass:fb:all}\,\ref{item:fb:all:f} follow from
    \begin{equation}
        \label{eq:fb:convergence:lipschitz-diff}
        \norm{F'(\mu) - F'(\nu)}_\infty \le L \norm{\mu-\nu}_\Wave
        \quad (\mu, \nu \in \Meas(\Omega)).
    \end{equation}
    The former is immediate while the proof of the latter is standard; compare \cite[Lemma 7.1 (iii) $\Rightarrow$ (v)]{clasonvalkonen2020nonsmooth}.
\end{remark}

\begin{example}
    For $F(\mu)=\frac12\norm{A\mu-b}^2$ with pre-adjointable $A \in \linear(\Masses; Y)$ and $Y$ a Hilbert space, \cref{thm:fb:convergence:weak}\,\cref{item:fb:convergence:weak:fprime-cont} follows by slightly adapting the previous remark. Indeed, we take $L$ satisfying $A_*A \le L \Wave$ as in \cref{example:fb:convergence:standard-f-lipschitz}.
    Since $F'(\mu)-F'(\nu) = A_*A(\mu-\nu)$, \eqref{eq:fb:convergence:lipschitz-diff}, we then have
    \[
        \norm{F'(\mu) - F'(\nu)}_\infty
        \le
        \norm{A_*}_{\linear(\R^n; \FullPredual)} \norm{A(\mu-\nu)}_{\R^n}
        \le
        \norm{A_*}_{\linear(\R^n; \FullPredual)} \sqrt{L} \norm{\mu-\nu}_{\Wave}.
     \]
     This establishes the required continuity.
\end{example}

\section{Inertial forward-backward}
\label{sec:fista}

We now provide an inertial version of \cref{alg:fb:fb}, i.e., FISTA \cite{beck2009fista} on measures.
We first describe the algorithm in \cref{sec:fista:description}, and then prove its convergence in \cref{sec:fista:convegence}.
For this we need the full generality of \cref{lemma:fb:general-estimate}.

\subsection{Algorithm description}
\label{sec:fista:description}

With $\lambda_0=1$, and $\breve\mu^0=\mu^0$, we define the inertial measures and parameters
\begin{equation}
    \label{eq:inertia:update}
    \nexxt{\breve\mu} \defeq (1+\theta_{k+1})\nexxt\mu - \theta_{k+1}\this\mu,
    \quad
    \theta_{k+1} \defeq \lambda_{k+1}(\inv\lambda_k-1),
    \quad
    \lambda_{k+1} \defeq \frac{2\lambda_k}{\lambda_k + \sqrt{4 + \lambda_k^2}}
    \quad
    (k \in \N).
\end{equation}
Then we “rebase” \cref{alg:fb:fb} at $\this{\breve\mu}$ in place of $\this\mu$ to obtain the inertial algorithm \cref{alg:inertia:inertia}.

\begin{algorithm}[t]
    \caption{Inertial forward-backward for Radon norm regularisation ($\mu$FISTA)}
    \label{alg:inertia:inertia}
    \begin{algorithmic}[1]
        \Require Regularisation parameter $\alpha>0$; convex and pre-differentiable $F: \Masses \to \R$;
            self-adjoint particle-to-wave operator $\Wave \in \linear(\Meas(\Omega); \FullPredual)$.
        \State Choose tolerances $\{\epsilon_{k+1}\}_{k \in \N} \subset (0, \infty)$ and a step length parameter $\tau>0$
            subject to \cref{thm:inertia:convergence}.
        \State Choose a fractional tolerance $\kappa \in (0, 1)$ for finite-dimensional subproblems.
        \State Pick an initial iterate $\mu^0 \in \DiscreteMasses$
        \State Set $\lambda_0 \defeq 1$ and $\breve\mu^0 \defeq \mu^0$.
        \For{$k \in \N$}
            \State $\this v \defeq F'(\this{\breve\mu})$.
            \State $\nexxt\mu \defeq \textproc{insert\_and\_adjust}(\this{\breve\mu}, \tau \this v - \Wave\this{\breve\mu}, \tau\alpha, \epsilon_{k+1}, \kappa)$.
                \Comment{Solves \eqref{eq:fb:step-problem} with \cref{alg:fb:insert-and-remove}.}
            \State Calculate $\nexxt{\breve\mu}$ and $\lambda_{k+1}$ according to \eqref{eq:inertia:update}.
            \State Prune spikes with zero weight in \emph{both} $\mu^{k+1}$ and $\nexxt{\breve\mu}$.
                \Comment{Optionally also heuristic merging.}
        \EndFor
    \end{algorithmic}
\end{algorithm}

\subsection{Convergence}
\label{sec:fista:convegence}

\begin{theorem}[Inertial method convergence]
    \label{thm:inertia:convergence}
    Suppose \cref{ass:fb:all} holds and let $\opt\mu \in \Masses$ satisfy $0 \in \subdiff[F + G](\opt \mu)$.
    Let $\{\this\mu\}_{k \ge 1}$ be generated by \cref{alg:inertia:inertia} for some $\mu^0 \in \DiscreteMasses$ with the tolerance sequence $\{\epsilon_{k+1}\}_{k \in \N} \subset (0, \infty)$ satisfying
    \[
        \lim_{N \to \infty}  \lambda_N^2 \sum_{k=0}^{N-1} \inv\lambda_k\epsilon_{k+1} = 0.
    \]
    Then $[F+G](\mu^N) \to [F+G](\opt\mu)$, more precisely, $\lambda_N^2 = O(1/N^2)$ in
    \[
        [F+G](\mu^N)
        \le
        [F+G](\opt\mu)
        + \frac{\lambda_N^2}{\tau}\left(
            (\kappa + \norm{\opt\mu}_{\Masses}) \sum_{k=0}^{N-1} \inv\lambda_k\epsilon_{k+1}
            + C_0
        \right)
    \]
    for $C_0 = \tau([F+G](\mu^0)-[F+G](\opt\mu)) + \frac{1}{2}\norm{\mu^0-\opt\mu}_{\Wave}^2$.
\end{theorem}

\begin{proof}
    We have $\inv \lambda_N \ge  1+ N/2$, see \cite[Lemma 4.3]{beck2009fista} or \cite[Lemma 3.4]{tuomov-inertia}.
    This establishes that $\lambda_N^2 = O(1/N^2)$.
    Taking $\nextz \defeq \inv\lambda_k\nexxt\mu - (\inv\lambda_k - 1)\this\mu$, we have $\lambda_k(\nexxt z-\this\mu) = \nexxt\mu - \this\mu$ as well as $\lambda_k(\nextz-\thisz) = \nexxt\mu-\this{\breve \mu}$, as required by \cref{lemma:fb:general-estimate}; see \cite[(2.8)]{tuomov-inertia}.
    \Cref{alg:inertia:inertia,lemma:fb:step-problem:alg,ass:fb:all} ensure \eqref{eq:fb:step-problem} for all $k \in \N$.
    Therefore \cref{lemma:fb:epsilon} provides $\tilde\epsilon_{k+1} \in \subdiff E_k(\nexxt\mu)$ satisfying \eqref{eq:fb:tilde-epsilon-estimate}.
    Since $\tau L \le 1$ and $\lambda_k \in (0, 1]$, \cref{lemma:fb:general-estimate} thus establishes
    \begin{multline*}
        \frac{\lambda_k^2}{2}\norm{\nextz-\opt\mu}_{\Wave}^2
        + \tau([F+G](\nexxt\mu) - [F+G](\opt\mu))
        \\
        \le
        (1-\lambda_k)\tau([F+G](\this\mu)-[F+G](\opt\mu))
        +\frac{\lambda_k^2}{2}\norm{\thisz-\opt\mu}_{\Wave}^2
        +\lambda_k\dualprod{\tilde\epsilon_{k+1}}{\thisz - \opt\mu}.
    \end{multline*}
    Since $\lambda_k^{-2}(1-\lambda_k) = \lambda_{k-1}^{-2}$ for all $k \in \N$ when we set $\lambda_{-1} = 0$, see again, e.g., \cite[Lemma 3.4]{tuomov-inertia}, multiplying this expression by $\lambda_k^{-2}$ yields
    \begin{multline*}
        \frac{1}{2}\norm{\nextz-\opt\mu}_{\Wave}^2
        + \lambda_k^{-2}\tau([F+G](\nexxt\mu) - [F+G](\opt\mu))
        \\
        \le
        \lambda_{k-1}^{-2}\tau([F+G](\this\mu)-[F+G](\opt\mu))
        +\frac{1}{2}\norm{\thisz-\opt\mu}_{\Wave}^2
        +\inv\lambda_k\dualprod{\tilde\epsilon_{k+1}}{\thisz - \opt\mu}.
    \end{multline*}
    Summing over $k=0,\ldots,N-1$ establishes
    \[
        \lambda_N^{-2}\tau([F+G](\mu^N) - [F+G](\opt\mu))
        \\
        \le
        C_0
        + \sum_{k=0}^{N-1}\inv\lambda_k\dualprod{\tilde\epsilon_{k+1}}{\thisz - \opt\mu}.
    \]
    Estimating the sum with \eqref{eq:fb:tilde-epsilon-estimate}, and multiplying by $\lambda_N^2/\tau$, the claim follows.
\end{proof}

\noindent\begin{minipage}{\textwidth}\noindent
\begin{example}[Tolerance sequence for inertia]
    \label{ex:inertia:convergence}
    Since $\lambda_0=1$, we have
    \[
        \inv\lambda_{k+1} = \frac{1}{2}(1 + \sqrt{1 + 4\smash[b]{\lambda_k^{-2}}})
        \le
        1 + \inv\lambda_k
        \le \ldots \le
        k + 1.
    \]
    Therefore, taking $\epsilon_{k+1}=1/(k+1)^p$ for some $p > 2$, \cref{thm:inertia:convergence} shows $O(1/N^2)$ function value convergence.
    If $p \in (1,2]$ as in \cref{ex:fb:tolerance:ergodic}, the rate reduces to $O(1/N)$.
\end{example}
\end{minipage}

\section{Primal-dual proximal splitting}
\label{sec:pdps}

We now provide a version of the primal-dual proximal splitting (PDPS) of \cite{chambolle2010first}.
We first describe the algorithm in \cref{sec:pdps:description}, and then sketch its convergence in \cref{sec:pdps:convegence}.
We allow the function $F$ in \eqref{eq:fb:problem} to take the form
\[
    F(\mu) = F_0(A\mu),
\]
where $A \in \linear(\Meas(\Omega); Y)$ with $Y$ a Hilbert space, and $F_0: Y \to \extR$ is convex, proper, and lower semicontinuous, but possibly nonsmooth. This is a relaxation from \cref{sec:fb,sec:fista} that required $F$ to have Lipschitz Fréchet derivative. For \eqref{eq:intro:problem}, i.e., a Gaussian noise model, we would take $F_0(y) = \frac{1}{2}\norm{y-b}_2^2$ on $Y=\R^n$. For a salt-and-pepper noise model, we would take $F_0(y) = \norm{y-b}_1$.

\subsection{Algorithm description}
\label{sec:pdps:description}

The optimality conditions for \eqref{eq:fb:problem} may now be written \cite{clasonvalkonen2020nonsmooth}
\[
    0 \in H(\opt\mu, \opty)
    \quad\text{where}\quad
    H(\mu, y)
    \defeq
    \begin{pmatrix}
        \subdiff G(\mu) + A_* y \\
        \subdiff F_0^*(y) - A\mu
    \end{pmatrix},
\]
where we recall from \eqref{eq:fb:g} that $G(\mu) = \alpha \norm{\mu}_{\Meas} + \delta_{\ge 0}(\mu)$.
With the help of the step length and preconditioning operators $W_k \in \linear(\FullPredual \times Y; \FullPredual\times Y)$ and $M_k \in \linear(\Meas(\Omega) \times Y; \FullPredual \times Y)$ defined as
\[
    W_k \defeq \begin{pmatrix}
        \tau_k \Id & 0 \\
        0 & \sigma_{k+1} \Id
    \end{pmatrix}
    \quad\text{and}\quad
    M_k \defeq \begin{pmatrix}
        \Wave & - \tau_k A_* \\
        \omega_k\sigma_{k+1} A & \Id
    \end{pmatrix},
\]
for some step length and over-relaxation parameters $\{(\tau_k, \sigma_k, \omega_k)\}_{k \in \N}$, we introduce as in \cite{tuomov-proxtest,clasonvalkonen2020nonsmooth,he2012convergence} the implicit form algorithm
\begin{equation}
    \label{eq:pdps:implicit}
    \begin{pmatrix} \tilde\epsilon_{k+1} \\ 0 \end{pmatrix} \in W_k H(\nextu) + M_k(\nextu-\thisu),
\end{equation}
where $\thisu=(\this\mu, \thisy)$, and with exact steps, we would have $\tilde\epsilon_{k+1}=0$.
Following \cref{sec:fb} we, however, replace the first line primal update by \cref{eq:fb:step-problem} with $\this v = A_* \thisy$.
The second line dual update reads
\[
    0 \in \sigma_{k+1} \subdiff F_0^*(\nexty) - \sigma_{k+1}A[(1+\omega_k)\nexxt\mu-\omega_k\this\mu] + (\nexty-\thisy).
\]
Taking $(\tau_k, \sigma_k, \omega_k) \equiv (\tau, \sigma, 1)$, and rewriting the dual update in terms of an explicit proximal mapping, in analogy with \cref{alg:fb:fb}, we obtain \cref{alg:pdps:pdps}.

When $F_0^*$ is strongly convex with factor $\gamma_{F_0^*}$, we can also accelerate
\begin{equation}
    \label{eq:pdps:accel}
    \tau_{k+1} \defeq \tau_k/\omega_k
    \quad\text{and}\quad
    \sigma_{k+1} \defeq \sigma_k\omega_k
    \quad\text{for}\quad
    \omega_k \defeq 1/\sqrt{1+\gamma_{F_0^*}\sigma_k},
\end{equation}
replacing $2\nexxt\mu-\this\mu$ in the dual update by $(1+\omega_k)\nexxt\mu-\omega_k\this\mu$.
We can take $\omega_k \defeq 1/\sqrt{1+2\gamma_{F_0^*}\sigma_k}$ if only iterate convergence is desired, and no function value convergence.

\begin{algorithm}[t]
    \caption{Primal-dual proximal splitting for Radon norm regularisation ($\mu$PDPS)}
    \label{alg:pdps:pdps}
    \begin{algorithmic}[1]
        \Require Regularisation parameter $\alpha>0$; convex and Fréchet-differentiable $F_0: Y \to \R$; pre-adjointable $A \in \linear(\Meas(\Omega); Y)$ with $Y$ a Hilbert space; self-adjoint $\Wave \in \linear(\Meas(\Omega); \FullPredual)$.
        \State Choose tolerances $\{\epsilon_{k+1}\}_{k \in \N} \subset (0, \infty)$, and step length parameters $\tau,\sigma>0$ satisfying $\tau\sigma A_*A < \Wave$.
        \State Choose a fractional tolerance $\kappa \in (0, 1)$ for finite-dimensional subproblems.
        \State Pick initial iterates $(\mu^0, y^0) \in \DiscreteMasses \times Y$.
        \For{$k \in \N$}
            \State $\this v \defeq A_* \this y$.
            \State $\nexxt\mu \defeq \textproc{insert\_and\_adjust}(\this{\mu}, \tau \this v - \Wave\this\mu, \tau\alpha, \epsilon_{k+1}, \kappa)$.
                \Comment{Solves \eqref{eq:fb:step-problem} with \cref{alg:fb:insert-and-remove}.}
            \State Prune zero weight spikes from $\nexxt\mu$.
                \Comment{Optionally also apply a spike merging heuristic.}
            \State $\nexty \defeq \prox_{\sigma F_0^*}(\this y + \sigma A[2\nexxt\mu-\this\mu])$.
        \EndFor
    \end{algorithmic}
\end{algorithm}

\subsection{Sketch of convergence}
\label{sec:pdps:convegence}

Following \cite{tuomov-proxtest}, see also \cite{clasonvalkonen2020nonsmooth}, we introduce for some testing parameters $\phi_k,\psi_k >0$ the \term{testing operator} $Z_k \in \linear(\FullPredual \times Y; \FullPredual \times Y)$ defined by
\[
    Z_k \defeq \begin{pmatrix}
        \phi_k \Id & 0 \\
        0 & \psi_{k+1} \Id
    \end{pmatrix}.
\]
We need $Z_kM_k \in \linear(\Meas(\Omega) \times Y; \FullPredual \times Y)$ to be self-adjoint and positive semi-definite. By standard arguments \cite{tuomov-proxtest,clasonvalkonen2020nonsmooth} this holds when $\omega_k = \inv\sigma_{k+1}\inv\psi_{k+1}\phi_k\tau_k$ and $\psi_k\sigma_k=\phi_k\tau_k$ as well as $\tau_k\sigma_k A_*A < \Wave$ with $\{\psi_k\}_{k \in \N}$ non-decreasing.

We then test \eqref{eq:pdps:implicit} by the application of $\dualprod{Z_k\freevar}{\nextu-\optu}$ for some $\optu \in \inv H(0)$. This yields
\begin{equation}
    \label{eq:pdps:main-estimate0}
    \begin{aligned}[t]
    \phi_k\dualprod{\tilde\epsilon_{k+1}}{\nexxt\mu-\opt\mu} & \in
    \phi_k\tau_k\dualprod{\subdiff G(\nexxt\mu) + A_* \nexty}{\nexxt\mu-\opt\mu}
    \\
    \MoveEqLeft[-1]
    +\psi_{k+1}\sigma_{k+1}\dualprod{\subdiff F_0^*(\nexty) - A\nexxt\mu}{\nexty-\opty}
    + \dualprod{Z_kM_k(\nextu-\thisu)}{\nextu-\optu}.
    \end{aligned}
\end{equation}
Using the definition of the convex subdifferential and the three-point identity for $Z_kM_k$, now
\begin{equation}
    \label{eq:pdps:main-estimate}
    \begin{aligned}[t]
        \phi_k\dualprod{\tilde\epsilon_{k+1}}{\nexxt\mu-\opt\mu}
        &
        \ge
        \frac{1}{2}\norm{\nextu-\optu}_{Z_kM_k}^2
        - \frac{1}{2}\norm{\thisu-\optu}_{Z_kM_k}^2
        + \frac{1}{2}\norm{\nextu-\thisu}_{Z_kM_k}^2
        \\
        \MoveEqLeft[-1]
        + \phi_k\tau_k[G(\nexxt\mu)-G(\opt\mu)]
        + \psi_{k+1}\sigma_{k+1}[F_0^*(\nexty)-F_0^*(\opty)]
        \\
        \MoveEqLeft[-1]
        + \frac{\psi_{k+1}\sigma_{k+1}\gamma_{F_0^*}}{2}\norm{\nexty-\opty}_Y^2
        \\
        \MoveEqLeft[-1]
        + \phi_k\tau_k[\dualprod{A_*\nexty}{\nexxt\mu-\opt\mu}
        - \psi_{k+1}\sigma_{k+1}\dualprod{A\nexxt\mu}{\nexty-\opty}.
    \end{aligned}
\end{equation}

If now $(\tau_k, \sigma_k, \omega_k) \equiv (\tau, \sigma, 1)$, i.e., we consider the unaccelerated algorithm, we can take $\phi_k \equiv 1$ and $\psi_{k+1} = \inv\sigma\tau$ while satisfying the non-negativity and self-adjointness of $Z_kM_k = Z M$, which is now  independent of the iteration. Then \eqref{eq:pdps:main-estimate} gives
\begin{equation}
    \label{eq:pdps:first-estimate}
    \frac{1}{2}\norm{\nextu-\optu}_{ZM}^2
    + \frac{1}{2}\norm{\nextu-\thisu}_{ZM}^2
    + \tau \mathcal{G}(\nextu, \optu)
    \le
    \frac{1}{2}\norm{\thisu-\optu}_{ZM}^2
    + \dualprod{\tilde\epsilon_{k+1}}{\nexxt\mu-\opt\mu},
\end{equation}
where the \term{Lagrangian gap functional}
\[
    \begin{aligned}
    \mathcal{G}(u, \optu)
    &
    \defeq
    G(\mu)-G(\opt\mu)
    + F_0^*(y)-F_0^*(\opty)
    + \dualprod{A_*y}{\mu-\opt\mu}
    - \dualprod{A\mu}{y-\opty}
    \\
    &
    =
    [G(\mu) - F_0^*(\opty) + \dualprod{A\mu}{\opty}]
    -
    [G(\opt\mu) - F_0^*(y) + \dualprod{A\opt\mu}{y}].
    \end{aligned}
\]
We have $\mathcal{G}(u, \optu) \ge 0$ when $0 \in H(\optu)$; compare \cite{clasonvalkonen2020nonsmooth}.
Summing \eqref{eq:pdps:first-estimate} over $k=0,\ldots,N-1$, estimating $\dualprod{\tilde\epsilon_{k+1}}{\nexxt\mu-\opt\mu}$ with \cref{lemma:fb:epsilon}, and using Jensen's inequality to pass to ergodic variables, we get:

\begin{theorem}[Ergodic convergence of the Lagrangian gap functional]
    \label{thm:pdps:convergence:function-ergodic}
    Let $\Wave \in \linear(\Meas(\Omega); \FullPredual)$ be defined as $\Wave\mu = \rho * \mu$ on the compact set $\Omega$, where $0 \not\equiv \rho \in C_0(\R^n)$ is symmetric and positive definite, i.e., $\rho(-y)=\rho(y)$ for all $y$, and $\fourier[\rho] \ge 0$.
    Let $\optu \in \Masses \times Y$ with $Y$ a Hilbert space satisfy $0 \in H(\optu)$, and $\{(\this\mu,\this y)\}_{k \ge 1}$ be generated by \cref{alg:pdps:pdps} for some $(\mu^0, y^0) \in \DiscreteMasses \times Y$, step lengths $\tau,\sigma>0$ satisfying $\tau\sigma A_*A \le \Wave$, and the tolerance sequence $\{\epsilon_{k+1}\}_{k \in \N} \subset (0, \infty)$ satisfying
    \[
        \lim_{N \to \infty}  \frac{1}{N} \sum_{k=0}^{N-1} \epsilon_{k+1} = 0;
    \]
    see \cref{ex:fb:tolerance:ergodic}. Then $\mathcal{G}(\tilde u^N; \optu) \to 0$, more precisely,
    \[
        0 \le \mathcal{G}(\tilde u^N; \optu)
        \le
        \frac{1}{N\tau}\left( (\kappa+\norm{\opt\mu}_{\Masses}) \sum_{k=0}^{N-1} \epsilon_{k+1} + \frac{1}{2}\norm{u^0-\optu}_{ZM}^2\right).
    \]
\end{theorem}

Weak convergence of the iterates can also be obtained following the arguments of \cref{thm:fb:convergence:weak}.
When $F_0^*$ is strongly convex, we can with the step length parameter updates \eqref{eq:pdps:accel} obtain $O(1/N^2)$ convergence of the gap functionals and $O(1/N)$ convergence of the dual iterates $\{\thisy\}_{k \in \N}$, however, the treatment of the gap functional is very technical; we refer to \cite{tuomov-proxtest,clasonvalkonen2020nonsmooth}. We therefore skip the treatment of the gap functional, and directly use in \eqref{eq:pdps:main-estimate0} the fact that $0 \in H(\optu)$ together with the monotonicity of $\subdiff G$ and the $\gamma_{F_0^*}$-strong monotonicity of $F_0^*$.
This gives in place of \eqref{eq:pdps:main-estimate0} the estimate
\begin{equation}
    \label{eq:pdps:main-estimate1}
    \begin{aligned}[t]
        \phi_k\dualprod{\tilde\epsilon_{k+1}}{\nexxt\mu-\opt\mu}
        &
        \ge
        \frac{1}{2}\norm{\nextu-\optu}_{Z_kM_k}^2
        - \frac{1}{2}\norm{\thisu-\optu}_{Z_kM_k}^2
        + \frac{1}{2}\norm{\nextu-\thisu}_{Z_kM_k}^2
        \\
        \MoveEqLeft[-1]
        + \psi_{k+1}\sigma_{k+1}\gamma_{F_0^*}\norm{\nexty-\opty}_Y^2.
    \end{aligned}
\end{equation}
Taking
$
    \phi_k \equiv 1
$
and
$
    \psi_{k+1} = \psi_k(1+2\sigma_k\gamma_{F_0^*}),
$
we have
\[
    Z_{k+1}M_{k+1}
    =
    Z_kM_k
    + \begin{pmatrix} 0 & 0 \\ 0 & \psi_{k+1}\sigma_{k+1} \gamma_{F_0^*} \Id\end{pmatrix}
    + \Xi_k,
\]
for some skew-adjoint $\Xi_k$. Since operator-relative norms are invariant with respect skew-adjoint components, summing \eqref{eq:pdps:main-estimate1} over $k=0,\ldots,N-1$ yields
\begin{equation}
    \label{eq:pdps:main-estimate2}
    \frac{1}{2}\norm{u^N-\optu}_{Z_{N-1}M_{N-1}}^2
    \le
    \frac{1}{2}\norm{u^0-\optu}_{Z_0M_0}^2
    + \sum_{k=0}^{N-1} \phi_k\dualprod{\tilde\epsilon_{k+1}}{\nexxt\mu-\opt\mu}.
\end{equation}
For some $\kappa \in (0, 1)$, we have
\[
    Z_{N-1}M_{N-1} \ge \begin{pmatrix}
        0 & 0 \\
        0 & (1-\kappa)\psi_N \Id
    \end{pmatrix}
\]
Therefore, dividing \eqref{eq:pdps:main-estimate2} by $\psi_N$, which grows at the rate $\Omega(N^2)$ \cite[Lemma 10.7]{clasonvalkonen2020nonsmooth}, we obtain $O(1/N)$ convergence of the dual iterates $\{\thisy\}_{k \in \N}$.

\section{Numerical experience}
\label{sec:numerical}

We now discuss our implementations of the proposed methods, and their practical performance compared to conditional gradient methods from the literature.
We first describe the implementation and parametrisation details in \cref{sec:numerical:details}.
We then describe the sample problems that we solve in \cref{sec:numerical:experiments}.
We finish with a report and discussion on the performance in \cref{sec:numerical:results}.

\subsection{Algorithm implementation and parametrisation}
\label{sec:numerical:details}

We implemented \cref{alg:fb:fb,alg:inertia:inertia,alg:pdps:pdps} (“our methods” $\mu$FB, $\mu$FISTA, and $\mu$PDPS) as well as the “relaxed” and “fully corrective” conditional gradient methods \cite[Algorithm 5.1]{brediespikkarainen2013inverse} and \cite[Algorithm 2]{walter2019linear}, denoted FWr and FWf.
All are applicable to the squared data term $F(\mu) = \frac{1}{2}\norm{A\mu-b}^2$. $\mu$PDPS is also applicable to $F(\mu)=\norm{A\mu-b}_1$.
Our Rust implementation is available on Zenodo \cite{tuomov-pointsource-codes}.

\paragraph{Finite-dimensional subproblems}
Both conditional gradient methods find a maximiser of $-F'(\this\mu)(\freevar)$.\footnote{The maximisation of $\abs{F'(\this\mu)(\freevar)}$ in \cite{brediespikkarainen2013inverse} corresponds to a version of \eqref{eq:intro:problem} without the non-negativity constraint; compare \cref{rem:fb:unconstrained}. The more general cone constraints of \cite{walter2019linear} readily treat this constraint.}
They add the point to the support $S$ of $\this\mu$, and then adapt weights by solving the finite-dimensional subproblem
\[
    \min_w J(w) \defeq F(Pw) + G(Pw),
\]
where $P \in \linear(\R^n; \Meas(\Omega))$ maps weights $w \in \R^n$ to measures $\sum_{i=1}^n w_i \delta_{x_i}$ with support $S=\{x_1,\ldots,x_n\}$. In the finite-dimensional subproblems of \cref{alg:fb:fb,alg:inertia:inertia,alg:pdps:pdps}, by contrast, $F$ is linearised, and a quadratic penalty based on $\Wave$ is added:
\[
    \min_w J(w) \defeq \tau[F(\this\mu) + \iprod{\this v}{Pw-\this\mu} + G(Pw)] + \frac{1}{2}\norm{Pw-\this\mu}_\Wave^2.
\]

The FWr takes a single forward-backward step on each finite-dimensional subproblem after a specific initialisation \cite{brediespikkarainen2013inverse}.
For our methods and the FWf, we solve the subproblems to the accuracy $\inf_{q \in \subdiff J(w^\ell)} \norm{q}_\infty \le 0.1\epsilon_{k+1}$, where $\{\epsilon_{k+1}\}_{k \in \N}$ is the overall tolerance sequence.
For FWf we use forward-backward splitting capped at 2000 iterations.
For \cref{alg:fb:fb,alg:inertia:inertia,alg:pdps:pdps} we use a semismooth Newton method (SSN) modelled after \cite{hintermuller2002active}. It is based on the proximal optimality condition
\[
    w = \prox_{\tau_{\text{sub}} G \circ P}(w - \tau_{\text{sub}} P^* \this v - \tau_{\text{sub}} P^* \Wave P w).
\]
SSN can only be used with our methods because $A_*A$ is in general not positive definite, but $\Wave$ is subject to conditions in \cref{thm:wave:constr:weak-new}.\footnote{
    If the spikes of $\this\mu$ start grouping very close together, $P_*\Wave P$ may become ill-conditioned.
    If observed, this could be avoided by using a merging heuristic.
    Alternatively, a first-order method could be used for the subproblem.
}
We set the subproblem step length $\tau_{\text{sub}}$ in both methods as $0.99$ times an estimate of the Lipschitz factor of the gradient of the smooth part of $J$. For the conditional gradient methods this is $n$ times an estimate of $\norm{A_*A}$, and for our methods $n$ times an estimate of $\norm{\Wave}$.

\paragraph{Branch-and-bound subproblem}
To maximise $-F'(\this\mu)(\freevar)$ for the conditional gradient methods, or to do the corresponding step of  \cref{alg:fb:insert-and-remove} for our methods, we do branch-and-bound optimisation on a geometric bisection tree representation of weighted sums $\zeta(x) = \sum_{i=1}^n w_i \rho(x - x_i) + \sum_{j=1}^m \beta_j (\theta * \psi)(x - z_j)$ of functions with a small support.
Each node of the tree corresponds to a subcube of the domain, and maintains a list of active, non-zero components of the sum, and rough upper and lower bounds on it, within the subcube.
The sum is, of course, computationally easier if the number of active components is small. To achieve this, we construct the component functions to have the small supports that was one of our guiding principles in \cref{sec:sensorgrid:examples}.

The branch-and-bound maximisation starts by putting the domain $\Omega=[0,1]^n$ in a priority queue (a binary heap).
Based on upper bound estimates for each cube in the priority queue, it picks the most promising one for refinement. If the maximum of a second-order polynomial model of $\zeta$ within the cube is $0.1\epsilon_k$ of the rough upper bound, the cube is inserted back into the priority queue, marked as a candidate solution. Otherwise its $2^n$ subcubes are inserted into the priority queue, unmarked.
When a marked cube is taken from the priority queue, the local approximate maximiser is returned as an approximate maximiser of $\zeta$.
Cubes are pruned from the priority queue when it is apparent that they do not contain the solution.
For further details we refer to our Rust implementation and its documentation.

\paragraph{Merging heuristics}
In all algorithms we prune unneeded spikes with zero weights to improve computational performance and result visualisation. The conditional gradient methods are, however, \emph{very} weak at achieving completely zero weights. They instead depend on \term{merging heuristics} to reduce the number of spikes.
We use the following: if two spikes $w_i \delta_{x_i}$ and $w_j \delta_{x_j}$ of $\mu^{k+1}$ satisfy $\norm{x_i-x_j}_\infty \le 0.02$, we merge them as $(w_i+w_j)\delta_{(w_i x_i + w_j x_j)/(w_i+w_j)}$ if doing so does not increase the value of the data term $F$.
Then convergence is not affected.

For our \cref{alg:fb:fb,alg:inertia:inertia,alg:pdps:pdps} this is not true:
non-increase of the data term is \emph{not} sufficient to not affect convergence.
Instead, \eqref{eq:fb:step-problem} should be maintained. This is significantly more costly.
Therefore, we will \emph{not} use a merging in heuristic.
One is not needed: the methods tend to keep the spike count reasonable even without one.
This stems from the conservative insertion strategy in \cref{alg:fb:insert-and-remove}: \emph{if weight optimisation is sufficient to satisfy \eqref{eq:fb:step-problem}, a new point will not be inserted into the support of $\nexxt\mu$}.
However, to improve visual quality, we use the same merging heuristic as with FWr and FWf to postprocess the final iterate.

\paragraph{Parameters}

We always take the initial iterate $\mu^0=0$. For the PDPS, $y^0 \in \subdiff F_0(A\mu^0)$.
Based on trial and error, we take the tolerance sequence $\epsilon_k = 0.5 \tau \alpha/(1+0.2k)^{1.4}$, where $k$ is the iteration number.
This choice balances between fast initial convergence and not slowing down later iterations too much via excessive accuracy requirements.
Moreover, as a bootstrap heuristic, on the first 10 iterations, we insert in \cref{alg:fb:insert-and-remove} at most one point irrespective of the tolerance $\epsilon_k$.
This does not affect convergence, as the convergence theory can be applied starting from any fixed iteration number.
For $\mu$FB and $\mu$FISTA we take $\tau = 0.99 / L$ where $L$ satisfying $A_*A \le L \Wave$ is given by \cref{ex:convolution:gaussian} or \ref{ex:convolution:fast} based on the choice of the spread $\psi$ for $A$ and the corresponding kernel $\rho$ of $\Wave$.
For $\mu$PDPS we take $\tau_0 = 0.5 / \sqrt{L}$ and $\sigma_0 = 1.98 / \sqrt{L}$.
For the squared data term we use the acceleration scheme \eqref{eq:pdps:accel} with $\gamma_{F_0^*}=1$.
The implemented conditional gradient methods have no discretionary parameters \cite{brediespikkarainen2013inverse,walter2019linear}.

\subsection{Experiments}
\label{sec:numerical:experiments}

To compare the algorithms against one another, we use the squared data term $F(\mu)=\frac{1}{2}\norm{A\mu-b}^2$ in both $\Omega=[0, 1]$ and $\Omega=[0, 2]^2$.
For the forward operator $A$, we consider both the cut Gaussian spread of \cref{ex:convolution:gaussian}, and the “fast” spread of \cref{ex:convolution:fast}. We use the wave-to-particle operator $\Wave$ from the same example.
For the cut Gaussian the standard deviation $\sigma_u = \sigma_v = 0.05$, and the cut-off $a=0.15$.
For the “fast” case $\sigma = 0.16$.
For the radius $b$ of the rectangular sensor $\theta$ as in \cref{ex:convolution:rectangular-sensor} we take $0.4$ times the spacing between the sensors on a regular grid $\grid$. Thus 80\% of the domain is “sensed” by some sensor.
In $[0, 1]$ we use 100, and in $[0, 2]^2$ we use $16 \times 16$ equally spaced sensors.\footnote{%
    The small number of sensors along each axis in 2D is for visualisation purposes: our implementation scales well with computation using a $32 \times 32$ sensor grid with the same spread only requiring slightly over double the CPU time.
}

To generate the synthetic measurement data $b$, we apply $A$ to a ground-truth measure $\hat\mu$ with four spikes of distinct magnitudes, as depicted in our result figures. Then we add to each sensor reading independent Gaussian noise of standard deviation $0.2$ in 1D (both spreads), and standard deviation $0.1$ (cut Gaussian spread) or $0.15$ (“fast” spread) in 2D.
These choices produce in all cases 3.8--4.8 dB SSNR (57--65\% noise).
We use the trial-and-error regularisation parameter $\alpha=0.09$ (1D, cut Gaussian), $0.06$ (1D, “fast”), $0.19$ (2D, cut Gaussian), or $0.12$ (2D, “fast” spread).

We also demonstrate $\mu$PDPS on the $\ell^1$ data term $F(\mu)=\norm{A\mu-b}_1$ on the cut Gaussian spread in one dimension; the other combinations are available through our software implementation \cite{tuomov-pointsource-codes}. The parametrisation of $A$ is as above.
We apply salt-and-pepper noise of magnitude $m=0.6$ and probability $p=0.4$ in 1D.
Each sensor will therefore have noise values in $\{0,-m,m\}$ with corresponding probabilities $\{1-p, p/2, p/2\}$.
This gives again 4.8 dB SSNR (57\% noise). We set $\alpha = 0.1$.

\subsection{Results}
\label{sec:numerical:results}

\begin{figure}[t!]
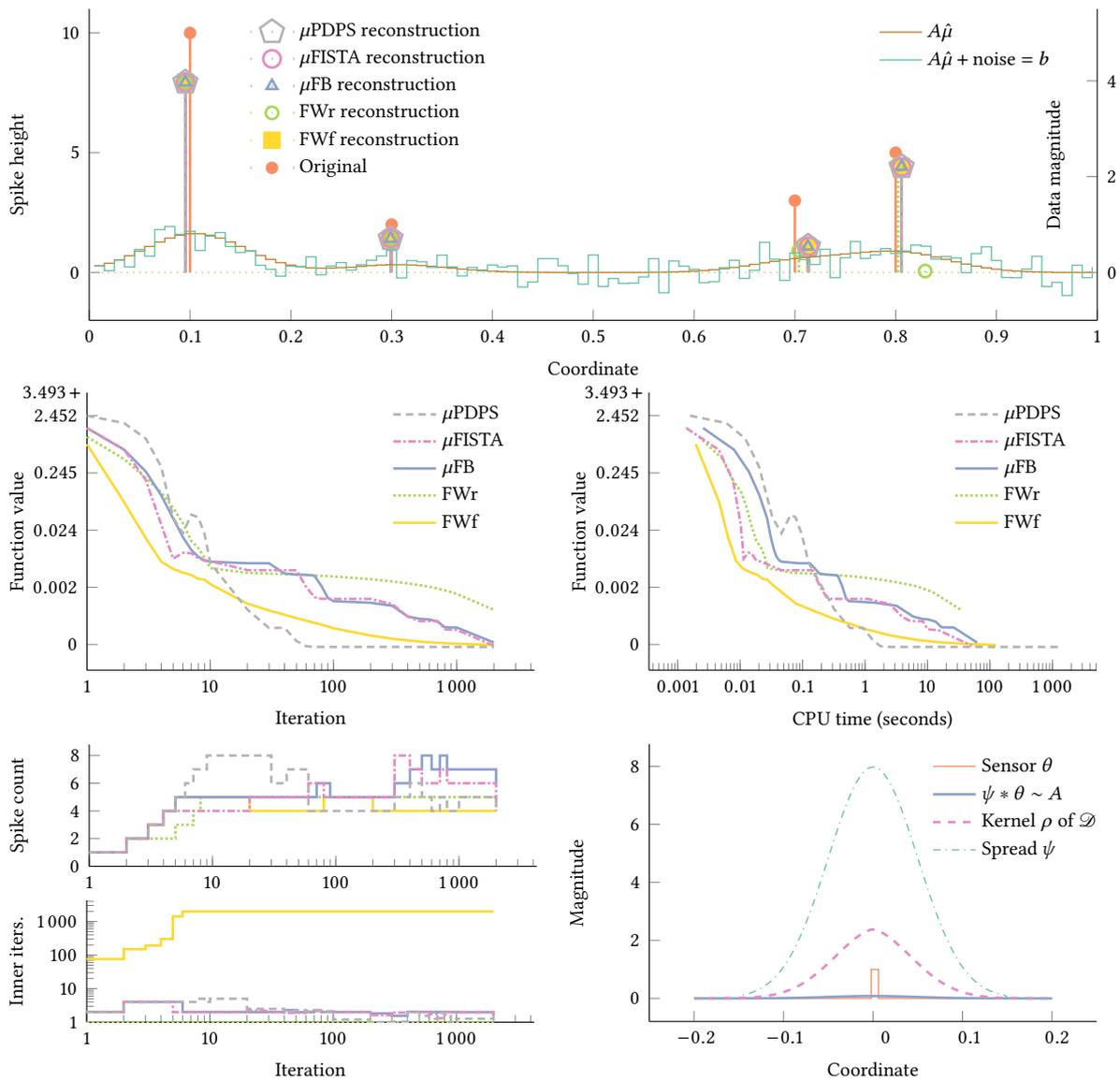

    \plotsONEd{pointsource1d}
    \caption{%
        Reconstructions and performance on 1D problem with cut Gaussian spread.
        \textbf{Top:} reconstruction and original data.
        The measurement data magnitude scale is on the right, spike magnitude on the left.
        \textbf{Middle:} Function value in terms of iteration count (left) and CPU time (right).
        The thin line indicates function value for $\mu$FB after postprocessing weight optimisation.
        \textbf{Bottom:} spike evolution, inner iteration count (left), and kernels (right).
        The thick lines indicate the spike count, and the thinner and dimmer lines the inner iteration count.
    }
    \label{fig:1dproblem}
\end{figure}

\begin{figure}[t!]
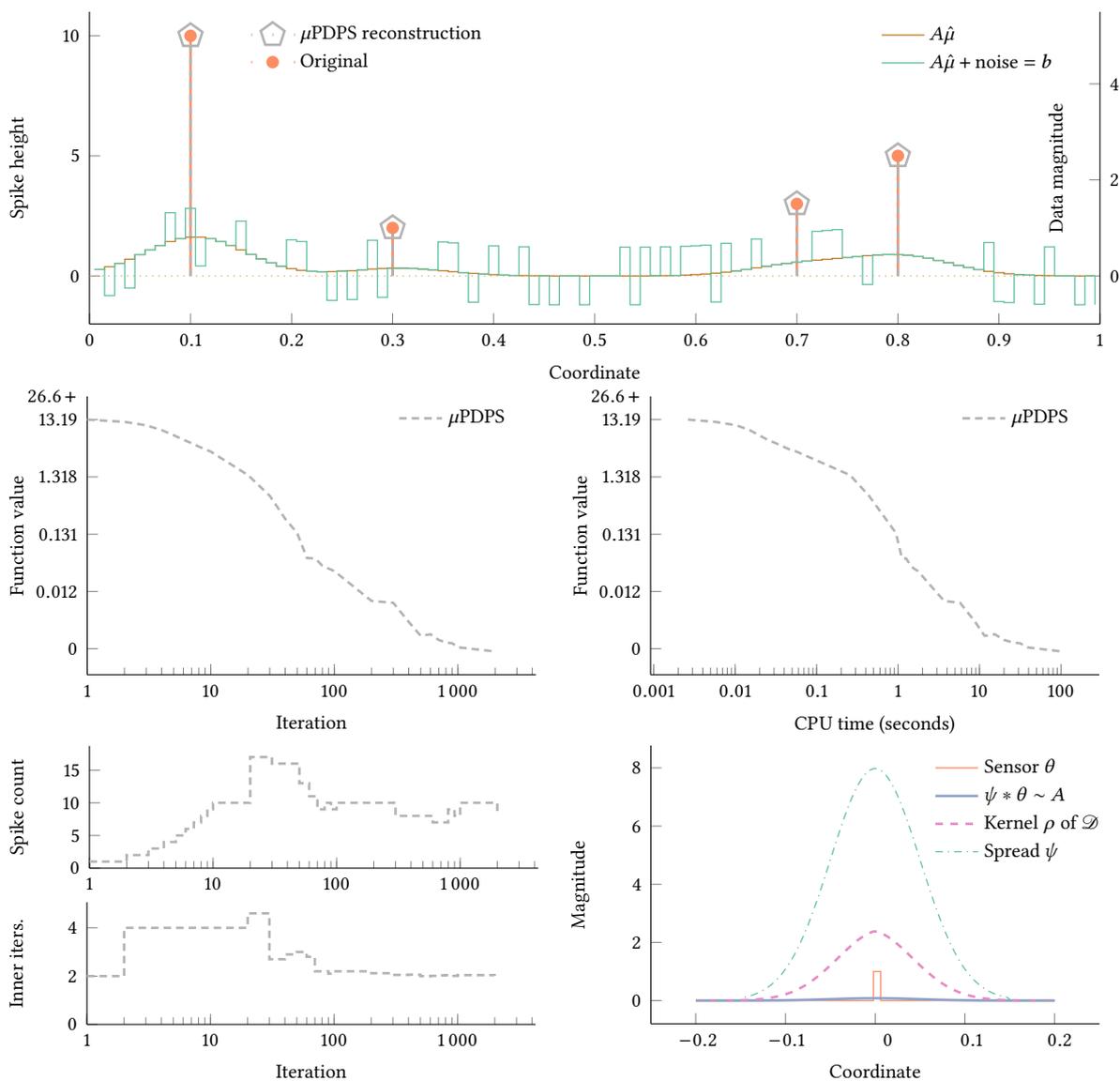

    \allalgsfalse
    \plotsONEd{pointsource1d_l1}
    \caption{%
        Reconstructions and performance on 1D problem with salt and pepper noise, $\ell_1$ data term, and gaussian spread.
        The plots are to be read the same way as \cref{fig:1dproblem}.
    }
    \label{fig:1dproblem_l1}
\end{figure}

We ran the experiments on a 2020 MacBook Air M1 with 16GB of memory.
We take advantage of the 4 high performance CPU cores of the 8-core machine by using 4 parallel computational threads to calculate $A^*z$ and for the branch-and-bound optimisation.
We report the performance of all of the algorithms for the cut Gaussian spread with squared data term and Gaussian noise in 1D in \cref{fig:1dproblem}, and in 2D in \cref{fig:2dproblem}. For the “fast” spread the results are in the corresponding \cref{fig:1dproblem_fast,fig:2dproblem_fast}. We also report the performance of $\mu$PDPS on salt-and-pepper noise with $\ell^1$ data term in one dimension in \cref{fig:1dproblem_l1}.
Each of the figures depicts the spread $\psi$, kernel $\rho$, and sensor $\theta$ involved in $A$ and $\Wave$. They also depict the noisy and noise-free data, the ground-truth measure $\hat\mu$, and the algorithmic reconstructions. The reconstructions or optimal solutions to \eqref{eq:fb:problem} cannot be expected to equal $\hat\mu$ due to noise and ill-conditioning of the inverse problem $A\mu=b$.

Each of the figures plots function value against both iteration count and CPU time spent.
The plots are logarithmic on both axes, and we sample the reported values logarithmically only on iterations $1,2,\ldots,10,20,\ldots,100,200,\ldots$.
We limit the number of iterations to 2000. The CPU time is a sum over the time spent by each computational thread, so several times the clock time requirement. Since the branch-and-bound optimisation stage gets increasingly more difficult closer to a solution, algorithms that converge the fastest generally require more CPU time to reach the final iteration.

The figures also indicate the spike count evolution and the number of iterations needed to solve the finite-dimensional subproblems. The subproblem iteration counts are averages over the corresponding period. FWr by design only ever takes one forward-backward iteration for the subproblems. FWf uses forward-backward splitting, as discussed in \cref{sec:numerical:details}. Recall that we cap the iteration count at 2000.
For $\mu$FB, $\mu$FISTA, and $\mu$PDPS the subproblem iteration count may include more than one spike insertion.

Since the weight optimisation subproblem of $\mu$FB, $\mu$FISTA, and $\mu$PDPS is not simply a finite-dimensional version of the original problem, as it is for FWf and FWr, the reported function value for the initial iterations can be suboptimal. An improved function value can be obtained by cheap postprocessing weight optimisation. This is graphed for $\mu$FB with a thin line in the iteration vs.~function value plot.

\subsection{Comparison}

\begin{figure}[t!]
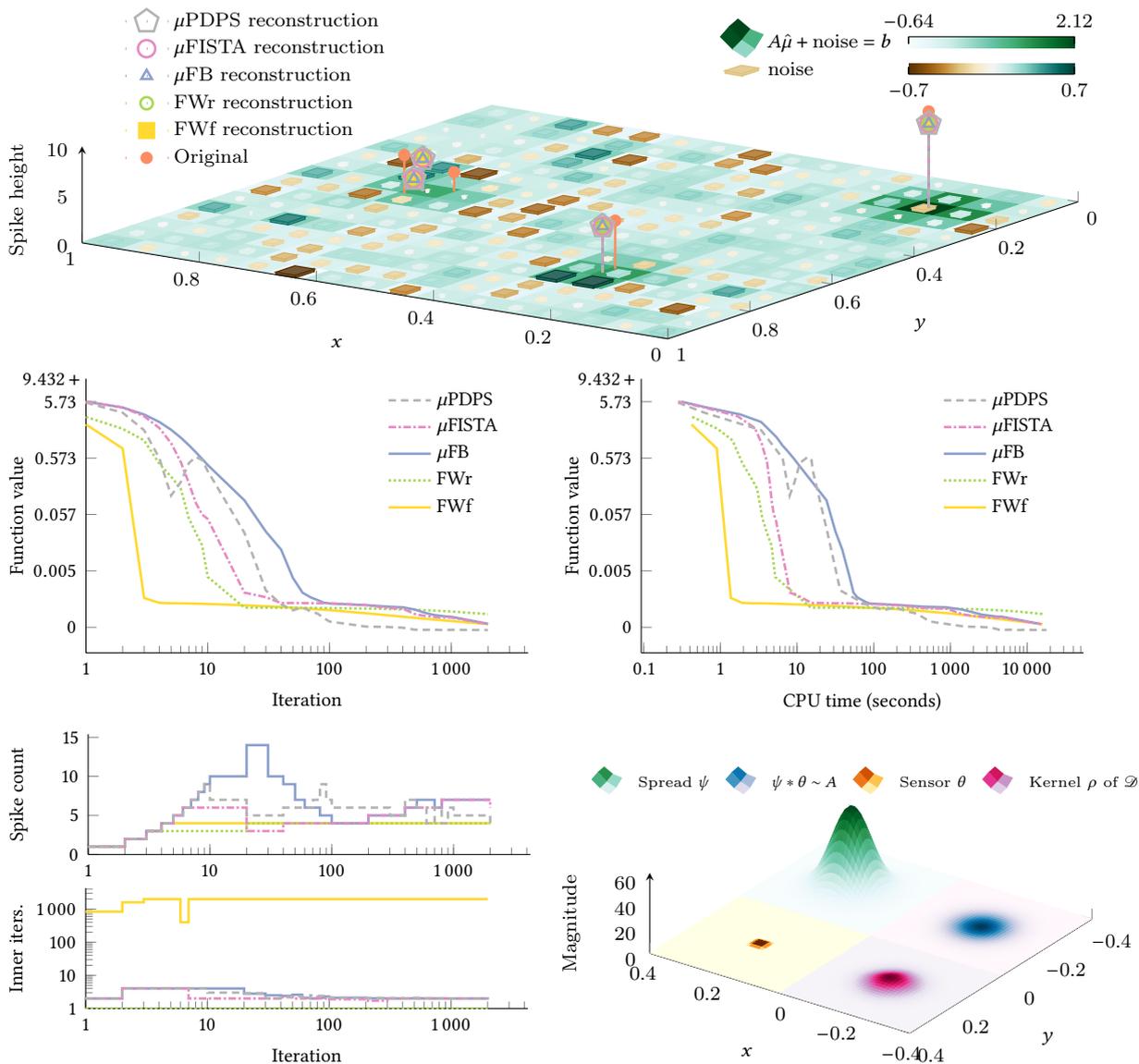

    \plotsTWOd{pointsource2d}
    \caption{%
        Reconstructions and performance on 2D problem with cut Gaussian spread.
        \textbf{Top:} reconstruction and original data. The area area of the top surface of the boxes is proportional the noise level of the underlying sensor, and their colour the sign of the noise.
        \textbf{Middle:} Function value in terms of iteration count and CPU time.
        The thin line indicates function value for $\mu$FB after postprocessing weight optimisation.
        \textbf{Bottom:} spike count evolution and kernels.
        The kernels have been shifted by $\pm 0.2$ in the $x$ and $y$ directions for visualisation-technical reasons.
    }
    \label{fig:2dproblem}
\end{figure}

Based on all the computations, for all algorithms, the “fast” spread generally has much lower CPU time requirements for convergence than the cut Gaussian. The bulk of the computational time is spent computing the $\erf$ of the expression for $\theta_z * \psi$ in \cref{lemma:convolution:fov-spread}.
The FWf and FWr seem quicker to start than the $\mu$FB and $\mu$FISTA, but eventually slow down, and the latter two overtake them. Sometimes, as in \cref{fig:2dproblem}, the conditional gradient methods are incredibly slow, and FWr even unstable.
Overall, however, it is difficult to rank the FWf, FWr, and $\mu$FISTA based on these experiments.
Generally $\mu$FISTA performs better than $\mu$FB, but not always significantly.
Perhaps surprisingly, as the squared data term does not require primal-dual structure for a prox-simple algorithm, the $\mu$PDPS is consistently a top-performer.
It is also the only one of the algorithms that can handle the $\ell^1$ data term.

\begin{figure}[t!]
    \plotsONEd{pointsource1d_fast}
    \caption{%
        Reconstructions and performance on 1D problem with the “fast” spread.
        The plots are to be read the same way as \cref{fig:1dproblem}.
    }
    \label{fig:1dproblem_fast}
\end{figure}

\begin{figure}[t!]
    \plotsTWOd{pointsource2d_fast}
    \caption{%
        Reconstructions and performance on 2D problem with the “fast” spread.
         The plots are to be read the same way as \cref{fig:2dproblem}.
    }
    \label{fig:2dproblem_fast}
\end{figure}

\input{pointsource0.bbl}



\end{document}

%% file: symbols.tex
\newcommand{\term}{\emph}

\newcommand{\field}[1]{\mathbb{#1}}
\newcommand{\N}{\mathbb{N}}

\newcommand{\R}{\field{R}}
\newcommand{\extR}{\overline \R}
\newcommand{\B}{B}

\newcommand{\norm}[1]{\|#1\|}

\newcommand{\abs}[1]{|#1|}
\newcommand{\adaptabs}[1]{\left|#1\right|}

\newcommand{\inv}[1]{#1^{-1}}

\newcommand{\freevar}{\,\boldsymbol\cdot\,}

\newcommand{\Union}\bigcup
\newcommand{\Isect}\bigcap
\newcommand{\union}\cup
\newcommand{\isect}\cap
\newcommand{\bigunion}\bigcup
\newcommand{\bigisect}\bigcap

\newcommand{\defeq}{:=}

\newcommand{\downto}{\searrow}
\newcommand{\upto}{\nearrow}
\newcommand{\subdiff}{\partial}

\DeclareMathOperator*{\argmin}{arg\,min}

\DeclareMathOperator{\closure}{cl}

\DeclareMathOperator{\sign}{sign}

\DeclareMathOperator{\Span}{span}

\DeclareMathOperator{\range}{ran}

\def \uminusSym{\setbox0=\hbox{$\cup$}\rlap{\hbox 
        to\wd0{\hss\raise0.5ex\hbox{$\scriptscriptstyle{-}$}\hss}}\box0}
    
\newcommand{\iprod}[2]{\langle #1,#2\rangle}

\def\llangle{\langle\kern-3pt\langle}
\def\rrangle{\rangle\kern-3pt\rangle}

\newcommand{\weakto}{\mathrel{\rightharpoonup}}
\def \weaktostarSym{\setbox0=\hbox{$\rightharpoonup$}\rlap{\hbox 
        to\wd0{\hss\raise1ex\hbox{$\scriptscriptstyle{*\,}$}\hss}}\box0}
    \def \weaktostar    {\mathrel{\weaktostarSym}}

\def\linear{\mathbb{L}}

\newcommand{\setto}{\rightrightarrows}

\def\extR{\overline \R}

\def\opt#1{\bar #1}

\def\this#1{#1^k}
\def\nexxt#1{#1^{k+1}}

\def\optu{{\opt{u}}}
\def\optx{{\opt{x}}}

\def\opty{{\opt{y}}}

\def\nextu{\nexxt{u}}
\def\nextx{\nexxt{x}}
\def\nextz{\nexxt{z}}
\def\nexty{\nexxt{y}}

\def\thisu{\this{u}}
\def\thisx{\this{x}}
\def\thisz{\this{z}}
\def\thisy{\this{y}}

\def\thisz{\this{z}}
\def\nextz{\nexxt{z}}

\DeclareMathOperator{\prox}{prox}

\DeclareMathOperator{\supp}{supp}

\def\d{\,d}

\def\dualprod#1#2{\langle #1|#2\rangle}

\newcommand{\Meas}{{\mathscr{M}}}

\newcommand{\Lebesgue}{{\mathscr{L}}}

\let\phi=\varphi
\let\epsilon=\varepsilon

\def\Id{\mathop{\mathrm{Id}}}

\def\Masses{\Meas(\Omega)}
\def\DiscreteMasses{\mathscr{Z}(\Omega)}
\def\FullPredual{C_0(\Omega)}
\def\Wave{{\mathscr{D}}}
\def\WaveSqr{\mathscr{Q}}
\def\fourier{\mathscr{F}}
\def\imag{\mathbb{i}}
\def\autoconvolution{\mathscr{A}}
\DeclareMathOperator{\rect}{rect}
\DeclareMathOperator{\sinc}{sinc}
\DeclareMathOperator{\tri}{tri}
\DeclareMathOperator{\erf}{erf}

%% file: symbols-texonly.tex
\renewrobustcmd{\downto}{{{\mathchoice%
            {\rotatebox[origin=c]{-20}{$\to$}}
            {\rotatebox[origin=c]{-20}{$\to$}}
            {\rotatebox[origin=c]{-20}{\scalebox{0.75}{$\to$}}}
            {\rotatebox[origin=c]{-20}{\scalebox{0.6}{$\to$}}}
}}}

\renewrobustcmd{\upto}{{{\mathchoice%
            {\rotatebox[origin=c]{20}{$\to$}}
            {\rotatebox[origin=c]{20}{$\to$}}
            {\rotatebox[origin=c]{20}{\scalebox{0.75}{$\to$}}}
            {\rotatebox[origin=c]{20}{\scalebox{0.6}{$\to$}}}
}}}

%% file: plotting.tex
\usepackage{tikz,pgfplots,pgfplotstable}
\usetikzlibrary{fpu}
\usepgfplotslibrary{colorbrewer,fillbetween}

\newlength{\figheight}
\newlength{\twofigheight}
\newlength{\dataplotheight}
\newlength{\twodkernelplotheight}
\pgfplotsset{
    compat = 1.17,
    scale only axis = false,
    scale mode = stretch to fill,
    every axis/.append style = {
        trim axis left,
        trim axis right,
    },
    plot coordinates/math parser = false,
    every axis label/.append style = {font = \scriptsize},
    every tick label/.append style = {font = \scriptsize},
    tick label style = {
        /pgf/number format/fixed,
        /pgf/number format/set thousands separator={\,}
    },
    legend style = {
        inner sep = 0pt,
        outer xsep = 5pt,
        outer ysep = 0pt,
        legend cell align = left,
        align = left,
        draw = none,
        fill = none,
        font = \scriptsize
    },
    legend pos = north east,
    width = \linewidth,
    height = \figheight,
    scaled x ticks = false,
    xminorticks = true,
    unbounded coords = jump,
    axis x line* = bottom,
    axis y line* = left,
    yminorticks = true,
    log ticks with fixed point,
    ylabel style = {
         at = {(yticklabel* cs:0.5)},
         anchor = center,
         yshift = 6ex,
    },
    y tick scale label style = {
         at = {(yticklabel* cs:1.0)},
         anchor = center,
         xshift = -3ex,
    }
}

\def\FB{$\mu$FB}
\def\FISTA{$\mu$FISTA}
\def\PDPS{$\mu$PDPS}

\def\datapath{results/}

\newif\ifallalgs
\allalgstrue

\pgfkeys{/matrixdim/.initial = 16}

\pgfplotsset{
    solid mark/.style = {
        every mark/.append style = {solid},
    }
}

\pgfplotstableset{
    create on use/inner_iters_mean/.style={
        create col/expr={\thisrow{inner_iters} / \thisrow{this_iters}},
    }
}

\pgfplotscreateplotcyclelist{reconstructions}{
    {color = Set2-B, line width = 1pt, mark = *, mark size = 2pt}, 
    {color = Set2-F, line width = 1pt, mark = square*}, 
    {color = Set2-E, line width = 1pt, mark = o, mark size = 2.5pt, densely dotted, solid mark}, 
    {color = Set2-C, line width = 1pt, mark = triangle, mark size = 2.5pt}, 
    {color = Set2-D, line width = 1pt, mark = o, mark size = 3.5pt, densely dashdotted, solid mark}, 
    {color = Set2-H, line width = 1pt, mark = pentagon, mark size = 5pt, densely dashed, solid mark}, 
}

\pgfplotsset{
    reconstructionplot/.style = {
        reverse legend,
        cycle list name = reconstructions,
        ycomb,
        mark size = 3pt,
    },
    performanceplot/.style = {
        reverse legend,
        no markers,
        cycle list shift = 1,
        cycle list name = reconstructions,
    }
}

\newcommand{\findminmax}[5]{
    \def\tempa{#3}
    \def\tempb{true}
    \ifx\tempa\tempb
        \pgfplotstablegetelem{0}{#2}\of{#1}%
        \pgfmathparse{\pgfplotsretval}
        \global\let#5\pgfmathresult
        \pgfplotstablegetelem{0}{#2}\of{#1}%
        \pgfmathparse{\pgfplotsretval}
        \global\let#4\pgfmathresult
    \fi
    \pgfplotstablegetrowsof#1
    \pgfmathparse{\pgfplotsretval-1}
    \foreach \i in {0,...,\pgfmathresult}{
        \pgfplotstablegetelem{\i}{#2}\of{#1}%
        \pgfmathparse{max(#5, \pgfplotsretval)}
        \global\let#5\pgfmathresult
        \pgfplotstablegetelem{\i}{#2}\of{#1}%
        \pgfmathparse{min(#4, \pgfplotsretval)}
        \global\let#4\pgfmathresult
    }
}

\newcommand{\SetMinMax}[3]{%
    \pgfplotstableread{#1}{#3}
    \findminmax{#3}{#2}{true}{\ExperimentMinValue}{\ExperimentMaxValue}
}

\newcommand{\UpdMinMax}[3]{
    \pgfplotstableread{#1}{#3}
    \findminmax{#3}{#2}{false}{\ExperimentMinValue}{\ExperimentMaxValue}
}

\pgfplotsset{
    fixedylog/.code 2 args={
        \pgfmathsetmacro{\base}{10^(-#1)}
        \pgfmathsetmacro{\logscale}{#1}
        \gdef\scaledlog##1{%
            (log10(\base + ((##1-\ExperimentMinValue)/(\ExperimentMaxValue-\ExperimentMinValue))) + \logscale)%
        }
        \pgfmathsetmacro{\logmax}{\scaledlog{\ExperimentMaxValue}}
        \pgfplotsset{
            y filter/.expression={\scaledlog{y}},
            ytick={0, 0.25*\logmax, 0.50*\logmax, 0.75*\logmax, \logmax},
            /pgf/number format/.cd,#2,
            scaled y ticks = manual:{$\pgfmathprintnumber{\ExperimentMinValue}\,+$}%
                {\pgfmathparse{10^((##1) - \logscale) - \base) * (\ExperimentMaxValue - \ExperimentMinValue)}}
        }
    }
}

\def\recodeps#1{%
    \tikzpicturedependsonfile{\datapath#1/orig.txt}%
    \tikzpicturedependsonfile{\datapath#1/pdps_reco.txt}%
    \ifallalgs%
    \tikzpicturedependsonfile{\datapath#1/fb_reco.txt}%
    \tikzpicturedependsonfile{\datapath#1/fw_reco.txt}%
    \tikzpicturedependsonfile{\datapath#1/fwrelax_reco.txt}%
    \tikzpicturedependsonfile{\datapath#1/fista_reco.txt}%
    \tikzpicturedependsonfile{\datapath#1/b_hat.txt}%
    \tikzpicturedependsonfile{\datapath#1/b_noisy.txt}%
    \fi%
}
\def\logdeps#1{%
    \tikzpicturedependsonfile{\datapath#1/pdps_log.txt}%
    \ifallalgs%
    \tikzpicturedependsonfile{\datapath#1/fb_log.txt}%
    \tikzpicturedependsonfile{\datapath#1/fista_log.txt}%
    \tikzpicturedependsonfile{\datapath#1/fw_log.txt}%
    \tikzpicturedependsonfile{\datapath#1/fwrelax_log.txt}%
    \fi%
}

\def\dataplotONEd#1{
    \tikzsetnextfilename{#1_data}
    \recodeps{#1}
    \begin{tikzpicture}
        \pgfplotsset{set layers}
        \begin{axis}[%
            scale only axis,
            width = 0.9\linewidth,
            height = \dataplotheight,
            xmin = 0, xmax = 1,
            ymax = 11, ymin = -2,
            legend pos = north west,
            legend style = {xshift=10ex},
            ylabel = {Spike height},
            xlabel = {Coordinate},
            reconstructionplot
            ]

            \draw [color = BrBG-E, line width = 0.5pt, dotted] (0,0)--(1,0);

            \addplot table [x = x0, y = weight]{\datapath#1/orig.txt};
            \addlegendentry{Original}

            \ifallalgs
            \addplot table [x = x0, y = weight]{\datapath#1/fw_reco.txt};
            \addlegendentry{FWf reconstruction}

            \addplot table [x = x0, y = weight]{\datapath#1/fwrelax_reco.txt};
            \addlegendentry{FWr reconstruction}

             \addplot table [x = x0, y = weight]{\datapath#1/fb_reco.txt};
            \addlegendentry{{\FB} reconstruction}

            \addplot table [x = x0, y = weight]{\datapath#1/fista_reco.txt};
            \addlegendentry{{\FISTA} reconstruction}
            \else
            \pgfplotsset{cycle list shift = 4}
            \fi

           \addplot table [x = x0, y = weight]{\datapath#1/pdps_reco.txt};
            \addlegendentry{{\PDPS} reconstruction}
        \end{axis}
        \begin{axis}[%
            scale only axis,
            width = 0.9\linewidth,
            height = \dataplotheight,
            xmin = 0, xmax = 1,
            ymax = 5.5, ymin = -1,
            legend pos = north east,
            axis y line* = right,
            axis x line = none,
            ylabel = {Data magnitude},
            ylabel style = {
                at = {(yticklabel cs:0.5)},
                anchor = center,
            },
            ]

            \addplot [color = BrBG-D, line width = 0.5pt, const plot]
                table [header=false]{\datapath#1/b_hat.txt};
            \addlegendentry{$A\hat\mu$}

            \addplot [color = Set2-A, line width = 0.5pt, const plot]
                table [header=false]{\datapath#1/b_noisy.txt};
            \addlegendentry{$A\hat\mu + \text{noise} = b$}

        \end{axis}
    \end{tikzpicture}%
    \pgfplotsset{set layers=none}%
}

\def\dataplotTWOd#1{
    \tikzsetnextfilename{#1_data}
    \recodeps{#1}
    \begin{tikzpicture}
        \pgfplotstableread{\datapath#1/b_hat.txt}{\bdata}
        \pgfplotstableread{\datapath#1/b_noisy.txt}{\bnoisy}
        \pgfplotstablecreatecol[copy column from table = {\bnoisy}{[index]2}]{noisy}{\bdata}
        \pgfplotstableclear{\bnoisy}
        \pgfplotstablecreatecol[expr = {\thisrow{noisy} - \thisrowno{2}}]{noise}{\bdata}
        \findminmax{\bdata}{noisy}{true}{\noisyMin}{\noisyMax}
        \findminmax{\bdata}{noise}{true}{\noiseMinX}{\noiseMaxX}
        \pgfmathsetmacro{\noiseMax}{max(abs(\noiseMinX),abs(\noiseMaxX)))}
        \pgfmathsetmacro{\noiseMin}{-\noiseMax}
        \pgfplotsset{set layers}
        \begin{axis}[%
            scale only axis,
            width = 0.9\linewidth,
            height = \dataplotheight,
            xmin = 0, xmax = 1,
            ymin = 0, ymax = 1,
            zmax = 20, zmin = 0,
            legend pos = north east,
            legend style = {
                name = datalegend,
                yshift = 0ex,
                xshift = -12ex,
            },
            axis x line = none,
            axis y line = none,
            axis z line = none,
            view/az = 215,
            view/el = 60,
            point meta = explicit,
            point meta rel = per plot,
            mesh/rows = \pgfkeysvalueof{/matrixdim},
            mesh/cols = \pgfkeysvalueof{/matrixdim},
            mesh/ordering = x varies,
            z buffer = sort,
            ]

            \addplot3[
                matrix plot*, 
                colormap/BuGn,
                line width = 0.5pt
            ] table [
                header = false,
                x index = 0,
                y index = 1,
                z expr = {0},
                meta = noisy,
            ]{\bdata};
            \addlegendentry{$A\hat\mu + \text{noise} = b$}

            \addplot3 [
                scatter,
                only marks,
                mark = cube*,
                cube/size z = 1pt,
                shader = flat,
                colormap/BrBG,
                line width = 0.5pt,
                scatter/use mapped color = {
                    draw = mapped color,
                    fill = mapped color!80,
                },
                visualization depends on = {
                    10*sqrt(abs(meta)) \as \noisew
                },
                scatter/@pre marker code/.append style = {
                    /pgfplots/cube/size x = \noisew,
                    /pgfplots/cube/size y = \noisew,
                },
                legend image post style = {
                    /pgfplots/cube/size z = 1pt,
                    /pgfplots/cube/size x = 7pt,
                    /pgfplots/cube/size y = 7pt,
                    mark options = {
                        fill = BrBG-F!80,
                        draw = BrBG-F,
                    }
                },
                x filter/.expression = { abs(\thisrow{noise}) < 0.00000001 ? nan : x },
                unbounded coords =discard,
            ]  table [
                header = false,
                x index = 0,
                y index = 1,
                z expr = {0},
                meta = noise,
            ]{\bdata};
            \addlegendentry{noise}
        \end{axis}
        \begin{axis}[%
            scale only axis,
            width = 0.9\linewidth,
            height = \dataplotheight,
            xmin = 0, xmax = 1,
            ymin = 0, ymax = 1,
            zmax = 10, zmin = 0,
            legend pos = north west,
            legend style = { yshift = 1ex },
            zlabel = {Spike height},
            xlabel = {$x$},
            ylabel = {$y$},
            axis z line = left,
            view/az = 215,
            view/el = 60,
            reconstructionplot
            ]

            \addplot3 table [x = x0, y = x1, z = weight,]{\datapath#1/orig.txt};
            \addlegendentry{Original}

            \ifallalgs
            \addplot3 table [x = x0, y = x1, z = weight]{\datapath#1/fw_reco.txt};
            \addlegendentry{FWf reconstruction}

            \addplot3 table [x = x0, y = x1, z = weight]{\datapath#1/fwrelax_reco.txt};
            \addlegendentry{FWr reconstruction}

            \addplot3 table [x = x0, y = x1, z = weight]{\datapath#1/fb_reco.txt};
            \addlegendentry{{\FB} reconstruction}

            \addplot3 table [x = x0, y = x1, z = weight]{\datapath#1/fista_reco.txt};
            \addlegendentry{{\FISTA} reconstruction}
            \else
            \pgfplotsset{cycle list shift = 4}
            \fi

            \addplot3 table [x = x0, y = x1, z = weight]{\datapath#1/pdps_reco.txt};
            \addlegendentry{{\PDPS} reconstruction}

        \end{axis}
        \pgfplotstableclear{\bdata}
        \begin{axis}[
            hide axis,
            scale only axis,
            height = 0pt,
            width = 0pt,
            colormap/BuGn,
            colorbar horizontal,
            point meta min = \noisyMin,
            point meta max = \noisyMax,
            colorbar style = {
                name = noisybar,
                at = (datalegend.north east),
                yshift = 0.3ex,
                anchor = north west,
                axis x line* = top,
                xtick = {\noisyMin,\noisyMax},
                font = {\tiny},
                width = 0.15\linewidth,
                height = 0.01\linewidth,
                axis y line = none,
                tick style = { major tick length = 3pt, semithick, color = black },
            }
            ]
            \addplot [draw=none] coordinates { (0, 0) (1, 1)};
        \end{axis}%
        \begin{axis}[
            hide axis,
            scale only axis,
            height = 0pt,
            width = 0pt,
            colormap/BrBG,
            colorbar horizontal,
            point meta min = \noiseMin,
            point meta max = \noiseMax,
            colorbar style = {
                at = (noisybar.south west), 
                yshift = -2ex,
                anchor = north west,
                xtick = {\noiseMin,\noiseMax},
                font = {\tiny},
                width = 0.15\linewidth,
                height = 0.01\linewidth,
                axis y line = none,
                tick style = { major tick length = 3pt, semithick, color = black },
                anchor = south west,
            }
            ]
            \addplot [draw=none] coordinates { (0, 0) (1, 1)};
        \end{axis}%
    \end{tikzpicture}%
    \pgfplotsset{set layers=none}%
}

\def\performanceplots#1{%
    \SetMinMax{\datapath#1/pdps_log.txt}{value}{\pdpsPerformance}%
    \ifallalgs
    \UpdMinMax{\datapath#1/fb_log.txt}{value}{\fbPerformance}%
    \UpdMinMax{\datapath#1/fista_log.txt}{value}{\fistaPerformance}%
    \UpdMinMax{\datapath#1/fw_log.txt}{value}{\fwPerformance}%
    \UpdMinMax{\datapath#1/fwrelax_log.txt}{value}{\fwrelaxPerformance}%
    \fi
    \tikzsetnextfilename{#1_performance}
    \logdeps{#1}%
    \begin{tikzpicture}[baseline]
        \begin{axis}[%
            xmode = log,
            xmin = 1,
            width = 0.5\linewidth,
            height = \figheight,
            axis x line* = bottom,
            axis y line* = left,
            ylabel = {Function value},
            xlabel = {Iteration},
            fixedylog = {4}{fixed,precision=3},
            performanceplot,
            ]

            \ifallalgs
            \addplot table [x = iter, y = value]{\fwPerformance};
            \addlegendentry{FWf}

            \addplot table [x = iter, y = value]{\fwrelaxPerformance};
            \addlegendentry{FWr}

            \addplot+[forget plot, line width = 0.5pt, opacity = 0.5] table [x = iter, y = post_value]{\fbPerformance};
            \addplot table [x = iter, y = value]{\fbPerformance};
            \addlegendentry{\FB}

            \addplot table [x = iter, y = value]{\fistaPerformance};
            \addlegendentry{\FISTA}
            \else
            \pgfplotsset{cycle list shift = 5}
            \fi

            \addplot table [x = iter, y = value]{\pdpsPerformance};
            \addlegendentry{\PDPS}
        \end{axis}
    \end{tikzpicture}%
    \quad%
    \tikzsetnextfilename{#1_cputime}%
    \logdeps{#1}%
    \begin{tikzpicture}[baseline]
        \begin{axis}[%
            width = 0.5\linewidth,
            height = \figheight,
            xmode = log,
            axis x line* = bottom,
            axis y line* = left,
            ylabel = {Function value},
            xlabel = {CPU time (seconds)},
            fixedylog = {4}{fixed,precision=3},
            performanceplot,
            ]

            \ifallalgs
            \addplot table [x = cpu_time, y=value]{\fwPerformance};
            \addlegendentry{FWf}

            \addplot table [x = cpu_time, y=value]{\fwrelaxPerformance};
            \addlegendentry{FWr}

            \addplot table [x = cpu_time, y=value]{\fbPerformance};
            \addlegendentry{\FB}

            \addplot table [x = cpu_time, y=value]{\fistaPerformance};
            \addlegendentry{\FISTA}
            \else
            \pgfplotsset{cycle list shift = 5}
            \fi

            \addplot table [x = cpu_time, y=value]{\pdpsPerformance};
            \addlegendentry{\PDPS}
        \end{axis}
    \end{tikzpicture}%
}

\def\evolutionplot#1{%
    \begin{minipage}[b]{0.5\linewidth}%
    \tikzsetnextfilename{#1_evolution}
    \logdeps{#1}%
    \begin{tikzpicture}%
        \begin{axis}[%
            width = \linewidth,
            height = \twofigheight,
            name = spike count,
            xmode = log,
            xmin = 1,
            ymin = 0,
            ylabel = {Spike count},
            const plot,  
            performanceplot,
            x tick label style = { draw = none },
            ]

            \ifallalgs
            \addplot table [x = iter, y = n_spikes]{\fwPerformance};
            \addplot table [x = iter, y = n_spikes]{\fwrelaxPerformance};
            \addplot table [x = iter, y = n_spikes]{\fbPerformance};
            \addplot table [x = iter, y = n_spikes]{\fistaPerformance};
            \else
            \pgfplotsset{cycle list shift = 5}
            \fi
            \addplot table [x = iter, y = n_spikes]{\pdpsPerformance};

        \end{axis}
    \end{tikzpicture}%
    \\%
    \tikzsetnextfilename{#1_spikes}%
    \logdeps{#1}%
    \begin{tikzpicture}
        \begin{axis}[%
            at = {(spike count.south west)},
            anchor = north west,
            width = \linewidth,
            height = \twofigheight,
            xmode = log,
            xmin = 1,
            ymin = {\ifallalgs 1\else 0\fi},
            ylabel = {Inner iters.},
            xlabel = {Iteration},
            ymode = {\ifallalgs log\else normal\fi},
            const plot,  
            performanceplot,
            ]

            \ifallalgs
            \addplot table [x = iter, y = inner_iters_mean]{\fwPerformance};
            \addplot table [x = iter, y = inner_iters_mean]{\fwrelaxPerformance};
            \addplot table [x = iter, y = inner_iters_mean]{\fbPerformance};
            \addplot table [x = iter, y = inner_iters_mean]{\fistaPerformance};
            \else
            \pgfplotsset{cycle list shift = 5}
            \fi
            \addplot table [x = iter, y = inner_iters_mean]{\pdpsPerformance};
        \end{axis}
    \end{tikzpicture}%
    \end{minipage}%
}
\def\kernelplotONEd#1{%
    \tikzsetnextfilename{#1_kernel}%
    \tikzpicturedependsonfile{\datapath#1/spread.txt}%
    \tikzpicturedependsonfile{\datapath#1/sensor.txt}%
    \tikzpicturedependsonfile{\datapath#1/base_sensor.txt}%
    \tikzpicturedependsonfile{\datapath#1/kernel.txt}%
    \begin{tikzpicture}
        \begin{axis}[%
            width = 0.5\linewidth,
            height = \figheight,
            xmin = -0.25, xmax = 0.25,
            axis x line* = bottom,
            axis y line* = left,
            ylabel = {Magnitude},
            xlabel = {Coordinate},
            legend pos = north east,
            legend style = { xshift = 2ex, outer xsep = 0pt },
            tick label style = { text width = 3ex, align = right },
            ]

            \addplot [color = Set2-B, line width = 0.5pt, const plot]
                table [header=false]{\datapath#1/sensor.txt};
            \addlegendentry{Sensor $\theta$}

            \addplot [color = Set2-C, line width = 1pt]
                table [header=false]{\datapath#1/base_sensor.txt};
            \addlegendentry{$\psi * \theta \sim A$}

            \addplot [color = Set2-D, line width = 1pt, dashed]
                table [header=false]{\datapath#1/kernel.txt};
            \addlegendentry{Kernel $\rho$ of $\Wave$ }

            \addplot [color = Set2-A, line width = 0.5pt, dashdotted]
                table [header=false]{\datapath#1/spread.txt};
            \addlegendentry{Spread $\psi$}

        \end{axis}
    \end{tikzpicture}%
}

\def\kernelplotTWOd#1{%
    \tikzsetnextfilename{#1_kernel}%
    \tikzpicturedependsonfile{\datapath#1/spread.txt}%
    \tikzpicturedependsonfile{\datapath#1/sensor.txt}%
    \tikzpicturedependsonfile{\datapath#1/base_sensor.txt}%
    \tikzpicturedependsonfile{\datapath#1/kernel.txt}%
    \begin{tikzpicture}
        \begin{axis}[%
            width = 0.5\linewidth,
            height = \twodkernelplotheight,
            xmin = -0.4, xmax = 0.4,
            ymin = -0.4, ymax = 0.4,
            axis z line = left,
            view/az = 215,
            view/el = 60,
            ylabel = {$y$},
            xlabel = {$x$},
            zlabel = {Magnitude},
            mesh/rows = 40,
            mesh/cols = 40,
            mesh/ordering = x varies,
            point meta rel = per plot,
            legend columns = 4,
            legend style = {
                at = {(0.5, 0.97)},
                anchor = north,
                column sep = 1ex,
                font = \tiny,
            },
            tick label style = { text width = 3ex, align = right },
            legend image post style = { scale = 0.7 },
            ]

            \addplot3 [
                surf, shader = flat, colormap/BuGn,
                restrict x to domain = 0:0.4,
                restrict y to domain = -0.4:0,
            ] table [
                header = false,
                x expr = {\thisrowno{0} + 0.2},
                y expr = {\thisrowno{1} - 0.2},
                z index = 2,
            ]{\datapath#1/spread.txt};
            \addlegendentry{Spread $\psi$}

            \addplot3 [
                surf, shader = flat, colormap/PuBu,
                restrict x to domain = -0.4:0,
                restrict y to domain = -0.4:0,
            ] table [
                header = false,
                x expr = {\thisrowno{0} - 0.2},
                y expr = {\thisrowno{1} - 0.2},
                z index = 2,
            ]{\datapath#1/base_sensor.txt};
            \addlegendentry{$\psi * \theta \sim A$}

            \addplot3 [
                surf, shader = flat, colormap/YlOrBr,
                restrict x to domain = 0:0.4,
                restrict y to domain = 0:0.4,
            ] table [
                header = false,
                x expr = {\thisrowno{0} + 0.2},
                y expr = {\thisrowno{1} + 0.2},
                z index = 2,
            ]{\datapath#1/sensor.txt};
            \addlegendentry{Sensor $\theta$}

            \addplot3 [
                surf, shader = flat, colormap/PuRd,
                restrict x to domain = -0.4:0,
                restrict y to domain = 0:0.4,
            ] table [
                header = false,
                x expr = {\thisrowno{0} - 0.2},
                y expr = {\thisrowno{1} + 0.2},
                z index = 2,
            ]{\datapath#1/kernel.txt};
            \addlegendentry{Kernel $\rho$ of $\Wave$ }

        \end{axis}
    \end{tikzpicture}%
}

\def\plotsONEd#1{%
    \tikzexternaldisable
    \dataplotONEd{#1}%
    \\%
    \performanceplots{#1}%
    \\%
    \evolutionplot{#1}%
    \kernelplotONEd{#1}%
}

\def\plotsTWOd#1{%
    \tikzexternalenable%
    \dataplotTWOd{#1}%
    \\%
    \tikzexternaldisable
    \performanceplots{#1}%
    \\%
    \evolutionplot{#1}%
    \tikzexternalenable
    \kernelplotTWOd{#1}%
    \tikzexternaldisable%
}

%% file: pointsource0.bbl
 \providecommand{\eprint}[1]{\href{http://arxiv.org/abs/#1}{arXiv:#1}}
  \providecommand{\noopsort}[1]{}
  \providecommand{\eprint}[1]{\href{http://arxiv.org/abs/#1}{arXiv:#1}}

%% file: pointsource0.bbl
\begin{thebibliography}{10}

\bibitem{ambrosio2000fbv}
L{.\nobreak\kern 0.33333em}Ambrosio, N{.\nobreak\kern 0.33333em}Fusco, and
  D{.\nobreak\kern 0.33333em}Pallara, \emph{Functions of Bounded Variation and
  Free Discontinuity Problems}, Oxford University Press, 2000.

\bibitem{beck2009fista}
A{.\nobreak\kern 0.33333em}Beck and M{.\nobreak\kern 0.33333em}Teboulle, A fast
  iterative shrinkage-thresholding algorithm for linear inverse problems,
  \emph{SIAM Journal on Imaging Sciences} 2 (2009),  183--202,
  \href{https://dx.doi.org/10.1137/080716542}{\nolinkurl{doi:10.1137/080716542}}.

\bibitem{blank2017extension}
L{.\nobreak\kern 0.33333em}Blank and C{.\nobreak\kern 0.33333em}Rupprecht, An
  extension of the projected gradient method to a {B}anach space setting with
  application in structural topology optimization, \emph{SIAM Journal on
  Control And Optimization} 55 (2017),  1481–1499,
  \href{https://dx.doi.org/10.1137/16m1092301}{\nolinkurl{doi:10.1137/16m1092301}}.

\bibitem{bredies2022generalized}
K{.\nobreak\kern 0.33333em}Bredies, M{.\nobreak\kern 0.33333em}Carioni,
  S{.\nobreak\kern 0.33333em}Fanzon, and F{.\nobreak\kern 0.33333em}Romero, A
  generalized conditional gradient method for dynamic inverse problems with
  optimal transport regularization, \emph{Foundations of Computational
  Mathematics}  (2022),
  \href{https://dx.doi.org/10.1007/s10208-022-09561-z}{\nolinkurl{doi:10.1007/s10208-022-09561-z}}.

\bibitem{brediespikkarainen2013inverse}
K{.\nobreak\kern 0.33333em}Bredies and H.\,K{.\nobreak\kern
  0.33333em}Pikkarainen, Inverse problems in spaces of measures, \emph{ESAIM:
  Control, Optimization and Calculus of Variations} 19 (2013),  190--218,
  \href{https://dx.doi.org/10.1051/cocv/2011205}{\nolinkurl{doi:10.1051/cocv/2011205}}.

\bibitem{bregman1967relaxation}
L.\,M{.\nobreak\kern 0.33333em}Bregman, The relaxation method of finding the
  common point of convex sets and its application to the solution of problems
  in convex programming, \emph{USSR Computational Mathematics and Mathematical
  Physics} 7 (1967),  200–217,
  \href{https://dx.doi.org/10.1016/0041-5553(67)90040-7}{\nolinkurl{doi:10.1016/0041-5553(67)90040-7}}.

\bibitem{candes2014towards}
E.\,J{.\nobreak\kern 0.33333em}Candès and C{.\nobreak\kern
  0.33333em}Fernandez-Granda, Towards a mathematical theory of
  super-resolution, \emph{Communications on Pure and Applied Mathematics} 67
  (2014),  906--956,
  \href{https://dx.doi.org/10.1002/cpa.21455}{\nolinkurl{doi:10.1002/cpa.21455}}.

\bibitem{casas2012approximation}
E{.\nobreak\kern 0.33333em}Casas, C{.\nobreak\kern 0.33333em}Clason, and
  K{.\nobreak\kern 0.33333em}Kunisch, Approximation of elliptic control
  problems in measure spaces with sparse solutions, \emph{SIAM Journal on
  Control And Optimization} 50 (2012),  1735--1752,
  \href{https://dx.doi.org/10.1137/110843216}{\nolinkurl{doi:10.1137/110843216}}.

\bibitem{casas2013parabolic}
E{.\nobreak\kern 0.33333em}Casas, C{.\nobreak\kern 0.33333em}Clason, and
  K{.\nobreak\kern 0.33333em}Kunisch, Parabolic control problems in measure
  spaces with sparse solutions, \emph{SIAM Journal on Optimization} 51 (2013),
  28--63,
  \href{https://dx.doi.org/10.1137/120872395}{\nolinkurl{doi:10.1137/120872395}}.

\bibitem{chambolle2010first}
A{.\nobreak\kern 0.33333em}Chambolle and T{.\nobreak\kern 0.33333em}Pock, A
  first-order primal-dual algorithm for convex problems with applications to
  imaging, \emph{Journal of Mathematical Imaging and Vision} 40 (2011),
  120--145,
  \href{https://dx.doi.org/10.1007/s10851-010-0251-1}{\nolinkurl{doi:10.1007/s10851-010-0251-1}}.

\bibitem{cheney2000approximation}
W{.\nobreak\kern 0.33333em}Cheney and W{.\nobreak\kern 0.33333em}Light, \emph{A
  Course in Approximation Theory}, volume 101 of Graduate Studies in
  Mathematics, American Mathematical Society, 2000.

\bibitem{chizat2021sparse}
L{.\nobreak\kern 0.33333em}Chizat, Sparse optimization on measures with
  over-parameterized gradient descent, \emph{Mathematical Programming}  (2021),
  \href{https://dx.doi.org/10.1007/s10107-021-01636-z}{\nolinkurl{doi:10.1007/s10107-021-01636-z}}.

\bibitem{chizat2023convergence}
L{.\nobreak\kern 0.33333em}Chizat, Convergence rates of gradient methods for
  convex optimization in the space of measures, 2023,
  \href{https://arxiv.org/abs/2105.08368}{\nolinkurl{arXiv:2105.08368}}.

\bibitem{chizat2018global}
L{.\nobreak\kern 0.33333em}Chizat and F{.\nobreak\kern 0.33333em}Bach, On the
  global convergence of gradient descent for over-parameterized models using
  optimal transport, in \emph{Advances in Neural Information Processing
  Systems}, S{.\nobreak\kern 0.33333em}Bengio, H{.\nobreak\kern
  0.33333em}Wallach, H{.\nobreak\kern 0.33333em}Larochelle, K{.\nobreak\kern
  0.33333em}Grauman, N{.\nobreak\kern 0.33333em}Cesa-Bianchi, and
  R{.\nobreak\kern 0.33333em}Garnett (eds.), volume~31, Curran Associates,
  Inc., 2018,
  \url{https://proceedings.neurips.cc/paper/2018/file/a1afc58c6ca9540d057299ec3016d726-Paper.pdf}.

\bibitem{clasonvalkonen2020nonsmooth}
C{.\nobreak\kern 0.33333em}Clason and T{.\nobreak\kern 0.33333em}Valkonen,
  Introduction to Nonsmooth Analysis and Optimization, 2020,
  \href{https://arxiv.org/abs/2001.00216}{\nolinkurl{arXiv:2001.00216}},
  \url{https://tuomov.iki.fi/m/nonsmoothbook_part.pdf}.
\newblock Work in progress.

\bibitem{courbot2021fast}
J.\,B{.\nobreak\kern 0.33333em}Courbot and B{.\nobreak\kern
  0.33333em}Colicchio, A fast homotopy algorithm for gridless sparse recovery,
  \emph{Inverse Problems} 37 (2021),  025002,
  \href{https://dx.doi.org/10.1088/1361-6420/abd29c}{\nolinkurl{doi:10.1088/1361-6420/abd29c}}.

\bibitem{denoyelle2019sliding}
Q{.\nobreak\kern 0.33333em}Denoyelle, V{.\nobreak\kern 0.33333em}Duval,
  G{.\nobreak\kern 0.33333em}Peyr{\'e}, and E{.\nobreak\kern 0.33333em}Soubies,
  The sliding {F}rank--{W}olfe algorithm and its application to
  super-resolution microscopy, \emph{Inverse Problems}  (2019),
  \href{https://dx.doi.org/10.1088/1361-6420/ab2a29}{\nolinkurl{doi:10.1088/1361-6420/ab2a29}}.

\bibitem{duval2017sparse}
V{.\nobreak\kern 0.33333em}Duval and G{.\nobreak\kern 0.33333em}Peyr{\'{e}},
  Sparse regularization on thin grids I: the Lasso, \emph{Inverse Problems} 33
  (2017),  055008,
  \href{https://dx.doi.org/10.1088/1361-6420/aa5e12}{\nolinkurl{doi:10.1088/1361-6420/aa5e12}}.

\bibitem{flinth2020linear}
A{.\nobreak\kern 0.33333em}Flinth, F{.\nobreak\kern 0.33333em}de~Gournay, and
  P{.\nobreak\kern 0.33333em}Weiss, On the linear convergence rates of exchange
  and continuous methods for total variation minimization, \emph{Mathematical
  Programming} 190 (2020),  221–257,
  \href{https://dx.doi.org/10.1007/s10107-020-01530-0}{\nolinkurl{doi:10.1007/s10107-020-01530-0}}.

\bibitem{fonseca2007mmc}
I{.\nobreak\kern 0.33333em}Fonseca and G{.\nobreak\kern 0.33333em}Leoni,
  \emph{Modern Methods in the Calculus of Variations: {$L^p$} Spaces},
  Springer, 2007.

\bibitem{frank1956algorithm}
M{.\nobreak\kern 0.33333em}Frank and P{.\nobreak\kern 0.33333em}Wolfe, An
  algorithm for quadratic programming, \emph{Naval Research Logistics
  Quarterly} 3 (1956),  95--110,
  \href{https://dx.doi.org/10.1002/nav.3800030109}{\nolinkurl{doi:10.1002/nav.3800030109}}.

\bibitem{he2012convergence}
B{.\nobreak\kern 0.33333em}He and X{.\nobreak\kern 0.33333em}Yuan, Convergence
  analysis of primal-dual algorithms for a saddle-point problem: from
  contraction perspective, \emph{SIAM Journal on Imaging Sciences} 5 (2012),
  119--149,
  \href{https://dx.doi.org/10.1137/100814494}{\nolinkurl{doi:10.1137/100814494}}.

\bibitem{hintermuller2002active}
M{.\nobreak\kern 0.33333em}Hinterm{\"u}ller, K{.\nobreak\kern 0.33333em}Ito,
  and K{.\nobreak\kern 0.33333em}Kunisch, The primal-dual active set strategy
  as a semismooth {N}ewton method, \emph{SIAM Journal on Optimization} 13
  (2002),  865--888 (2003),
  \href{https://dx.doi.org/10.1137/s1052623401383558}{\nolinkurl{doi:10.1137/s1052623401383558}}.

\bibitem{hormander2003analysis}
L{.\nobreak\kern 0.33333em}Hörmander, \emph{The Analysis of Linear Partial
  Differential Operators: I. Distribution Theory and Fourier Analysis},
  Classics in Mathematics, Springer Berlin Heidelberg, 2003,
  \href{https://dx.doi.org/10.1007/978-3-642-61497-2}{\nolinkurl{doi:10.1007/978-3-642-61497-2}}.

\bibitem{kammler2008first}
D.\,W{.\nobreak\kern 0.33333em}Kammler, \emph{A First Course in Fourier
  Analysis}, Cambridge University Press, 2008,
  \href{https://dx.doi.org/10.1017/CBO9780511619700}{\nolinkurl{doi:10.1017/cbo9780511619700}}.

\bibitem{lindberg2012mathematical}
J{.\nobreak\kern 0.33333em}Lindberg, Mathematical concepts of optical
  superresolution, \emph{Journal of Optics} 14 (2012),  083001,
  \href{https://dx.doi.org/10.1088/2040-8978/14/8/083001}{\nolinkurl{doi:10.1088/2040-8978/14/8/083001}}.

\bibitem{mattila1999geometry}
P{.\nobreak\kern 0.33333em}Mattila, \emph{Geometry of sets and measures in
  {E}uclidean spaces: Fractals and rectifiability}, Cambridge University Press,
  1999,
  \href{https://dx.doi.org/10.1017/CBO9780511623813}{\nolinkurl{doi:10.1017/cbo9780511623813}}.

\bibitem{opial1967weak}
Z{.\nobreak\kern 0.33333em}Opial, Weak convergence of the sequence of
  successive approximations for nonexpansive mappings, \emph{Bulletin of the
  American Mathematical Society} 73 (1967),  591--597,
  \href{https://dx.doi.org/10.1090/S0002-9904-1967-11761-0}{\nolinkurl{doi:10.1090/s0002-9904-1967-11761-0}}.

\bibitem{walter2019linear}
K{.\nobreak\kern 0.33333em}Pieper and D{.\nobreak\kern 0.33333em}Walter, Linear
  convergence of accelerated conditional gradient algorithms in spaces of
  measures, \emph{ESAIM: Control, Optimization and Calculus of Variations} 27
  (2021), ~38,
  \href{https://dx.doi.org/10.1051/cocv/2021042}{\nolinkurl{doi:10.1051/cocv/2021042}},
  \href{https://arxiv.org/abs/1904.09218}{\nolinkurl{arXiv:1904.09218}}.

\bibitem{robbins1971convergence}
H{.\nobreak\kern 0.33333em}Robbins and D{.\nobreak\kern 0.33333em}Siegmund, A
  convergence theorem for non negative almost supermartingales and some
  applications, \emph{Optimizing Methods in Statistics}  (1971),  233–257,
  \href{https://dx.doi.org/10.1016/b978-0-12-604550-5.50015-8}{\nolinkurl{doi:10.1016/b978-0-12-604550-5.50015-8}}.

\bibitem{rudin2006functional}
W{.\nobreak\kern 0.33333em}Rudin, \emph{Functional Analysis}, McGraw--Hill,
  2006.

\bibitem{tuomov-inertia}
T{.\nobreak\kern 0.33333em}Valkonen, Inertial, corrected, primal-dual proximal
  splitting, \emph{SIAM Journal on Optimization} 30 (2020),  1391--1420,
  \href{https://dx.doi.org/10.1137/18M1182851}{\nolinkurl{doi:10.1137/18m1182851}},
  \href{https://arxiv.org/abs/1804.08736}{\nolinkurl{arXiv:1804.08736}},
  \url{https://tuomov.iki.fi/m/inertia.pdf}.

\bibitem{tuomov-proxtest}
T{.\nobreak\kern 0.33333em}Valkonen, Testing and non-linear preconditioning of
  the proximal point method, \emph{Applied Mathematics and Optimization} 82
  (2020),
  \href{https://dx.doi.org/10.1007/s00245-018-9541-6}{\nolinkurl{doi:10.1007/s00245-018-9541-6}},
  \href{https://arxiv.org/abs/1703.05705}{\nolinkurl{arXiv:1703.05705}},
  \url{https://tuomov.iki.fi/m/proxtest.pdf}.

\bibitem{tuomov-firstorder}
T{.\nobreak\kern 0.33333em}Valkonen, First-order primal-dual methods for
  nonsmooth nonconvex optimisation, in \emph{Handbook of Mathematical Models
  and Algorithms in Computer Vision and Imaging}, K{.\nobreak\kern
  0.33333em}Chen, C.\,B{.\nobreak\kern 0.33333em}Schönlieb,
  X.\,C{.\nobreak\kern 0.33333em}Tai, and L{.\nobreak\kern 0.33333em}Younes
  (eds.), Springer, Cham, 2021,
  \href{https://dx.doi.org/10.1007/978-3-030-03009-4_93-1}{\nolinkurl{doi:10.1007/978-3-030-03009-4_93-1}},
  \href{https://arxiv.org/abs/1910.00115}{\nolinkurl{arXiv:1910.00115}},
  \url{https://tuomov.iki.fi/m/firstorder.pdf}.

\bibitem{tuomov-subreg}
T{.\nobreak\kern 0.33333em}Valkonen, Preconditioned proximal point methods and
  notions of partial subregularity, \emph{Journal of Convex Analysis} 28
  (2021),  251--278,
  \href{https://arxiv.org/abs/1711.05123}{\nolinkurl{arXiv:1711.05123}},
  \url{https://tuomov.iki.fi/m/subreg.pdf}.

\bibitem{tuomov-regtheory}
T{.\nobreak\kern 0.33333em}Valkonen, Regularisation, optimisation,
  subregularity, \emph{Inverse Problems} 37 (2021),  045010,
  \href{https://dx.doi.org/10.1088/1361-6420/abe4aa}{\nolinkurl{doi:10.1088/1361-6420/abe4aa}},
  \href{https://arxiv.org/abs/2011.07575}{\nolinkurl{arXiv:2011.07575}},
  \url{https://tuomov.iki.fi/m/regtheory.pdf}.

\bibitem{tuomov-pointsource-codes}
T{.\nobreak\kern 0.33333em}Valkonen, Proximal methods for point source
  localisation: the implementation, Software on Zenodo, 2022,
  \href{https://dx.doi.org/10.5281/zenodo.7402054}{\nolinkurl{doi:10.5281/zenodo.7402054}}.

\end{thebibliography}
